\documentclass{article}

\usepackage[hidelinks]{hyperref}       
\usepackage{url}            
\usepackage{booktabs}       
\usepackage{amsfonts}       
\usepackage{nicefrac}       
\usepackage{microtype}      
\usepackage{enumerate}
\usepackage{wrapfig}
\usepackage{authblk}
\usepackage[margin=1in]{geometry}

\usepackage{amsmath} 
\usepackage{amsfonts,amssymb,amsthm,amsmath,bm,mathtools}
\usepackage{algorithm}
\usepackage{algpseudocode}
\usepackage{systeme}
\usepackage{float}
\usepackage{tikz}
\usepackage{caption}
\usepackage{soul,xpatch}
\usepackage[normalem]{ulem}
\usepackage[english]{babel} 
\usepackage{placeins}
\usepackage{multirow}
\usepackage{bbm}
\usepackage{listings}
\usepackage{color} 
\usepackage{comment}

\usepackage{amsopn}  
\usepackage{lipsum}
\usepackage{graphicx}
 
\usepackage{color}
\usepackage{ amssymb }
\usepackage{bbm}
\usepackage{bm}

\usepackage[normalem]{ulem}

\usepackage{xcolor}
\usepackage{longtable}
\usepackage[utf8]{inputenc}
\usepackage{booktabs}
\usepackage{tabularx}

\usepackage{tablefootnote}
\usepackage{threeparttable}
\usepackage{bbding}
\usepackage{multirow}
\usepackage{booktabs}
\usepackage{array}
\usepackage{ragged2e}

\newcolumntype{P}[1]{>{\Centering\hspace{0pt}}p{#1}}
\newcolumntype{Z}{>{\centering\let\newline\\\arraybackslash\hspace{0pt}}X}

\usepackage{amsthm}

\usepackage{algorithm}

\makeatletter

\newcommand{\Eb}{\mathbb{E}}

\newcommand{\Rmnum}[1]{\expandafter\@slowromancap\romannumeral #1@}
\makeatother

\DeclareMathOperator*{\argmin}{arg\,min}

 \newcommand{\RN}[1]{\textup{\uppercase\expandafter{\romannumeral#1}}
}

\newtheorem{assumption}{Assumption}

 \newcommand{\vg}{\textbf{g}}
\newcommand{\cF}{\mathcal{F}}

\newcommand{\norm}[1]{\left\lVert#1\right\rVert}
\newcommand\inn[2]{\left\langle 
#1,#2\right\rangle}

\newcommand{\overm}{\frac{1}{m}}
\newcommand{\summ}{\sum_{i=1}^m}

\newcommand{\ep}{\epsilon}
\newcommand{\bx}{\bar{x}}
\newcommand{\bz}{\bar{z}}

\newcommand{\tx}{\tilde{x}}

\newcommand{\vx}{\textbf{x}}
\newcommand{\vy}{\textbf{y}}

\newcommand{\vv}{\textbf{v}}

\newcommand{\vz}{\textbf{z}}

\newcommand{\tz}{\tilde{z}}

\newcommand{\vs}{\bm{s}}

\newcommand{\tg}{\tilde{g}}

\newcommand{\cM}{\mathcal{M}}
  
\def \Gh {\mathcal{G}}

\def \cL {\mathcal{L}}

\newcommand{\Mud}{\texttt{M}_{\frac{1}{\delta}\,u}}
\newcommand{\Muds}{\texttt{M}_{\frac{1}{\delta}\,u}^{\star}}
\newcommand{\muM}{\mu_{\texttt{M}}}

\newtheorem{example}{Example} 
\newtheorem{theorem}{Theorem}
\newtheorem{lemma}[theorem]{Lemma} 
\newtheorem{proposition}[theorem]{Proposition} 
\newtheorem{remark}[theorem]{Remark}
\newtheorem{corollary}[theorem]{Corollary}
\newtheorem{definition}[theorem]{Definition}

 \usepackage{cite}

\title{DCatalyst: A Unified Accelerated Framework for Decentralized Optimization}
\author{Tianyu Cao\thanks{School of Industrial Engineering, Purdue University, West Lafayette, IN 47906 (email:  {\tt cao357@purdue.edu}). }\qquad Xiaokai Chen
\thanks{School of Industrial Engineering, Purdue University, West Lafayette, IN 47906 (email: {\tt chen4373@purdue.edu}).}\qquad Gesualdo Scutari\thanks{School of Industrial Engineering, Purdue University, West Lafayette, IN 47906 (email {\tt gscutari@purdue.edu}). \\
This work has been supported by the ORN Grant N. N000142412751. 
}
	} 

\begin{document}

\maketitle

\begin{abstract} We study decentralized optimization over a network of agents, modeled as graphs, with no central server. The goal is to minimize   $f+r$, where  $f$ represents a (strongly) convex function averaging the local agents' losses, and $r$ is a convex, extended-value function. 

We introduce DCatalyst, a unified black-box framework that integrates Nesterov acceleration into decentralized optimization algorithms. 
At its core,   DCatalyst operates as an \textit{inexact}, \textit{momentum-accelerated} proximal method (forming the outer loop) that seamlessly incorporates any selected decentralized algorithm (as the inner loop).  
We demonstrate that DCatalyst achieves optimal communication and computational complexity (up to log-factors) across various decentralized algorithms and problem instances.  Notably, it extends acceleration capabilities to problem classes  previously lacking accelerated solution methods, thereby broadening the effectiveness of decentralized methods. 

On the technical side, our framework introduce the {\it inexact estimating sequences}--a novel extension of the well-known Nesterov's estimating sequences, tailored for the minimization of  composite losses in decentralized settings. This method adeptly handles consensus errors and inexact solutions of agents' subproblems, challenges not addressed by existing models. 
 
\end{abstract}

\section{Introduction}

   We explore decentralized optimization across a network comprising $m>1$ agents. The network is modeled as an undirected, static graph, possibly with no centralized nodes (e.g., servers).   The agents aim to cooperatively solve the following optimization problem:   
\begin{equation}\label{eq:problem}
\min_{x\in\mathbb{R}^{d}}\,\,u(x) 
:= f(x)+r(x),\quad \text{with}\quad f(x) :=\overm\summ f_i(x) \tag{P}. 
\end{equation}
Here $f_i: \mathbb{R}^d\rightarrow \mathbb{R}$ is  the  cost function associated with agent $i$, assumed to be smooth and (strongly) convex, and known only to that agent. Additionally,    $r:\mathbb{R}^d \rightarrow \mathbb{R}\cup\{+\infty\}$ is an extended-value, convex function known to all agents, instrumental to enforce    constraints or promote some desirable structure on the solution,  such as  sparsity or low rank. 

   The problem formulation    \eqref{eq:problem}  is notably versatile, encompassing a variety of applications across diverse fields. This paper places particular emphasis on (supervised) decentralized machine learning problems, specifically framed as Empirical Risk Minimization (ERM).   In such formulations, each   $f_i$  has a finite-sum structure, resulting in the average of  a predefined loss $\ell(\bullet;\xi):\mathbb{R}^d\to \mathbb{R}$  over   data samples $\xi\in \mathcal D$, these being distributed among the agents in the network.  Specifically,  denoting by  $\{\xi_{i,1},\cdots,\xi_{i,n}\}\subset \mathcal{D}$ the set of $n$   samples accessed by agent $i$, drawn from an unknown distribution $\mathbb{P}$,   the empirical risk  $f_i$  in \eqref{eq:problem} reads 
\begin{equation}\label{eq:f_ij}
    f_i(x) := \frac{1}{n}\sum_{j=1}^n\underset{:=f_{ij(x)}}{\underbrace{\ell(x; \xi_{i,j})}}, 
\end{equation}
 where  $\ell(x;\xi_{i,j})$   
 quantifies the discrepancy between the hypothesized model, parameterized by $x$, and the data linked to $\xi_{i,j}$. Under mild assumptions on the loss function $\ell$ and i.i.d data distribution across the nodes,   agents' cost functions $f_i$ in (\ref{eq:f_ij})   reflect statistical similarities in  the data residing at different nodes, quantified by the following property   
 \begin{equation}\label{eq:similarity}
   \big\|\big(\nabla f(x) - \nabla f_i(x)\big) - \big(\nabla f(y) - \nabla f_i(y)\big)\big\|\leq \beta \|x-y\|,\quad \forall x,y\in \texttt{dom}\, r, 
\end{equation}
which holds with high probability \cite{DISCO,pmlr-v119-hendrikx20a}. Here, $\beta = \widetilde{\mathcal{O}}(1/\sqrt{n})$     measures   similarity  of $f_i$'s ($\widetilde{\mathcal{O}}$ hides log-factors and dependence on $d$).  
 
 In the postulated decentralized setting, lacking of a centralized node aggregating other agents' information, solving (\ref{eq:problem}) calls for the design of decentralized algorithms. The prototype decentralized scheme  employs an iterative process that interleaves local computations  with inter-node (possibly multiple) communication rounds. These communication rounds are instrumental to acquire global information, which subsequently informs the computational tasks in the ensuing iterations. Given that communication overhead often represents the primary constraint in distributed systems, overshadowing the cost of local parallel computations (e.g., \cite{Bekkerman_book11,Lian17}), recent  research efforts  have been directed toward crafting algorithms that are  {\it communication-efficient}.

 Employing Nesterov’s acceleration techniques offers a strategic pathway to reduce the number of iterations required for convergence, potentially alleviating the communication overhead of the algorithms. While  acceleration techniques have been 
 successfully implemented in centralized settings, including server-client systems  
  {\cite{nesterov1983method, beck2009fast,Alexandre2021AccMethods}},  their application to decentralized networks  
  remains  unsatisfactory.  

  Deferring to   Sec.~\ref{sec:related-works}  a  comprehensive review of the pertinent literature,   it is sufficient to highlight here  that current decentralized approaches  employing acceleration are  predominantly bespoke, as summarized in  Tables~\ref{table:scvx}-\ref{table:finite-sum}. Notably,  {\bf (i)} they are designed for {\it  special instances} of the general formulation \eqref{eq:problem}, lacking guarantees/applicability to wider  classes of problems, such as     \textit{composite formulations} ($r\neq 0$). These models are at the core of machine learning applications, which rely  on nonsmooth regularization to control overfitting; \textbf{(ii)}   
  they often overlook exploiting inherent structural function properties,  such as similarity (\ref{eq:similarity}), which could be leveraged for enhanced communication efficiency, especially  when the objective function is ill-conditioned; \textbf{(iii)}  their  convergence analyses   lack the generalizability  to new algorithmic classes addressing these broader problem domains. 
  
   This shows a clear research gap:  a {\it universal} analytical framework that encapsulates a broader class of accelerated decentralized algorithms is missing.   
  This paper   bridges precisely this gap, proposing a   {\it unified black-box} acceleration framework {\it à la Catalyst} {\cite{guler1992new,lin2015universal, lin2018catalyst, Alexandre2021AccMethods}} (notice that the Catalyst method is confined to centralized settings). 
  The innovative approach,     termed {\it Decentralized Catalyst} (DCatalyst),  brings Nesterov-style acceleration to a wide array of existing decentralized algorithms that were not initially designed for acceleration, making them applicable to Problem~\eqref{eq:problem} in its full generality, yet unsolved.  Notably,   DCatalyst achieves {\it optimal} communication complexity (up to logarithmic factors) for the majority of these decentralized algorithms.   Furthermore, it is distinguished by its design versatility, allowing for the explicit incorporation of functional structures, such as similarity or  nonsmoothness in composite form, and control of the  communication overhead versus local computational complexity.

\subsection{Related Works}\label{sec:related-works}
Given the focus of this study on decentralized algorithms, our literature review is tailored to methods employing acceleration  and suitable for mesh networks, excluding the extensive body of research dedicated to centralized settings (including server-client systems).
We organize the literature of relevant decentralized algorithms in three groups, based upon the specific instances of Problem~\eqref{eq:problem} they are applicable to.     The categorization is as follows. 

\textbf{(i) Strongly convex    $f$ (Table~\ref{table:scvx}):}   The literature documents numerous accelerated decentralized first-order methods tailored for this instance of \eqref{eq:problem}, primarily targeting smooth (unconstrained) problems, i.e., where $r=0$. These schemes  all achieve linear convergence in terms of communication complexity, albeit exhibiting diverse scaling behaviors with the functional and network parameters.     Specifically,  early algorithms  leveraging a primal-dual reformulation of \eqref{eq:problem} (with $r=0$) 
 {  \cite{ scaman2017optimal, kovalev2020optimal, li2020revisiting,  rogozin2020towards, dvinskikh2021decentralized, kovalev2021adom}  } or invoking penalty strategies 
 {\cite{li2020decentralized, uribe2020dual}}, demonstrate a convergence rate scaling as $\mathcal{O}(\sqrt{\kappa_{\ell}}/\sqrt{1-\rho})$, where $\kappa_{\ell} := {L_{\max}}/{\mu_{\min}}$ denotes the \textit{local} condition number of the agents' losses $f_i$. Here, $L_{\max}$ and $\mu_{\min}$ denote  the maximal smoothness constant and the minimal strong convexity constant across $f_i$'s, and  $\rho\in [0,1)$ represents  the connectivity of the network  (detailed definitions provided in Sec.~\ref{section:prelim}). The dependence on the local condition number $\kappa_\ell$ is not desirable--in fact  $\kappa_\ell$ can far exceed  the {\it global} condition number  $\kappa_g:=L/\mu$,  where $L$ and $\mu$   are the smoothness constant and the strong convexity constant of the {\it centralized} loss $f$.  Recent works   {\cite{ye2020decentralized,ye2023multi}} have introduced  decentralized  accelerations that improve  communication complexity to scale with $\mathcal{O}(\sqrt{\kappa_{g}}/\sqrt{1-\rho})$. These few methods are  applicable also to instances of \eqref{eq:problem} with $r\neq 0$. \\\indent  Given their nature as first-order methods, these algorithms do not leverage potential function similarity, e.g., in the form \eqref{eq:similarity}. This shortcoming is notably impactful for ill-conditioned functions, whose  global condition number, $\kappa_g$, is substantially large. Consequently, the communication complexity scaling with $\sqrt{\kappa_g}$  becomes less than ideal.  This situation is prevalent in numerous  ERM  problems, characterized by an optimal regularization parameter for test predictive performance that is exceedingly small.  For instance, in ridge regression scenarios, typical values read
     $\mu=\mathcal{O}(1/\sqrt{mn})$ and $L=\mathcal{O}(1)$ while  $\beta=\mathcal{O}(1/\sqrt{n})$  \cite{DISCO}.  This yields a global condition number $\kappa_g=\mathcal{O}(\sqrt{mn})$, in contrast with the ratio  $\beta/\mu=\mathcal{O}(\sqrt{m})$--the former growing with the local sample size $n$, while the latter remains independent.   These situations have spurred a wave of research focused on harnessing function similarity through statistical preconditioning \cite{DANE,reddi2016aide,yuan2019convergence,DISCO,fan2023communication}, some in conjunction with acceleration techniques \cite{pmlr-v119-hendrikx20a,dvurechensky2021hyperfast}, to enhance communication efficiency while accepting increased local computation costs. Such accelerated methods have achieved linear convergence with the number of communication steps scaling, albeit asymptotically \cite{pmlr-v119-hendrikx20a}, as $1+\mathcal{O}(\sqrt{\beta/\mu})$. This scaling  can be substantially more favorable  than $\sqrt{\kappa_g}$, aligning closely with the lower complexity bounds (modulo logarithmic factors) \cite{arjevani2015communication}. Nevertheless, these methods are inherently \textit{centralized}, as   they   rely on  a server node connected directly to all the other clients; hence, they are \textit{infeasible} for deployment over mesh networks. To the best of our knowledge, Network-DANE \cite{NetDane} and SONATA \cite{sun2020distributed}  stand out as the sole algorithms that exploit statistical similarity over mesh networks. However, neither method incorporates acceleration, resulting in a communication complexity of $\widetilde{\mathcal{O}}((1/\sqrt{1-\rho})\cdot\beta/\mu\cdot \log(1/\varepsilon))$, in contrast to the anticipated $\mathcal{O}(\sqrt{\beta/\mu})$ dependence. Currently, there lacks a theoretical framework that supports the acceleration of these methods while simultaneously harnessing the benefits of function similarity. The framework proposed in this paper closes this gap (see Sec.~\ref{sec:major-contributions}).\\
 
\indent \textbf{(ii) (non-strongly) convex $f$ (Table~\ref{table:cvx}):}  The literature of accelerated distributed algorithms tailored to (non-strongly) convex instances of \eqref{eq:problem} remains   sparse. Existing schemes  {\cite{li2020decentralized,li2020revisiting, uribe2020dual, dvinskikh2021decentralized}} are designed for smooth objectives ($r=0$). Among these, \cite{li2020decentralized} stands out by achieving  acceleration rates for both gradient and communication oracles (albeit dependent on $L_{\max}$, instead of the more favorable $L$), namely:  $\mathcal{O}(\sqrt{L_{\max}/\epsilon})$ and  $\mathcal{O}(\sqrt{L_{\max}/((1-\rho)\cdot\epsilon)})$, respectively.  In contrast, the realm of composite functions ($r \neq 0$) is notably underexplored. The sole distributed method in this category, \cite{li2020decentralized}, falls short of achieving full acceleration, with gradient and communication complexities scaling at $\mathcal{O}(L_{\max}/\epsilon)$ and  $\mathcal{O}(L_{\max}/(\sqrt{1-\rho}\cdot \epsilon))$, rather than the anticipated $\mathcal{O}(\sqrt{L/\epsilon})$ and  $\mathcal{O}(\sqrt{L/((1-\rho)\cdot\epsilon)})$, respectively.  Moreover, current methodologies do not leverage potential similarities across functions.   The proposed framework addresses these gaps, offering optimal acceleration (modulo logarithmic factors), e.g.,  to the SONATA algorithm \cite{sun2020distributed}, as well as harnessing function similarity, if any.     
\\\indent\textbf{(iii) $f_i$'s with finite-sum structure (Table~\ref{table:finite-sum}):}  In {\it centralized} settings, the synergy of Variance Reduction (VR) techniques with acceleration has been well-documented   to exploit the finite-sum structure of $f_i$'s while maintaining rapid convergence rates  {\cite{lin2015universal, katyusha2018, zhou18simple, zhou2019direct, lan2019unified, driggs2022accelerating}}. VR methods  enable the use of more cost-effective (stochastic) gradient iterations by sampling instances of the full gradient.     

  In contrast, {\it decentralized} acceleration techniques incorporating VR  are   markedly sparse, and confined to scenarios with {\it smooth} objectives $u$ ($r=0$)  {\cite{ hendrikx2020dual, hendrikx2021anOptimal,li2022variance}}. To our knowledge, the only VR-based decentralized method tailored for composite objectives is  {\cite{ye2021pmgtvr}}, yet this algorithm does not offer any form of acceleration. The proposed framework introduces the first systematic design of decentralized VR-based  methods  for  composite functions that provably achieve acceleration, thereby filling the gap in the current decentralized optimization landscape.
 
\begin{table*}[t!]\Large
\caption{\scriptsize  \textbf{Decentralized accelerated algorithms for strongly convex $f$ in \eqref{eq:problem}}: The most efficient decentralized accelerated algorithms in the literature for different instances of  \eqref{eq:problem}. 
Key parameters   include (definitions in   Sec.~\ref{section:prelim}): $\kappa_g :=  {L}/{\mu}$ (global condition number of $f$), 
$\kappa_{\ell} :=  {L_{\max}}/{ \mu_{\min}}$ (local condition number of $f_i$'s), 
$\beta$ (function similarity,  \eqref{eq:similarity}), and  $\rho\in [0,1)$ (network connectivity). Notation $\widetilde{\mathcal{O}}$ omits  poly-log factors of problem parameters, independent of  $\epsilon$. Notably, the proposed methods  stands as the first to achieve   acceleration either in terms of gradient or communication complexity, with the latter   harnessing similarity. }
\resizebox{1.0\columnwidth}{!}{\begin{tabular}{c|c|c|c}
\toprule 
 
\textbf{Algorithm}  & \textbf{Problem}& \textbf{$\#$ (Proximal) gradient} & \textbf{$\#$ Communications}\\

\midrule
OPAPC \cite{kovalev2020optimal}    & \multirow{2}{*}{  $f_i$ scvx, $r\equiv 0$} &  $\mathcal{O}\left(\sqrt{\kappa_{\ell}}\log\frac{1}{\epsilon}\right)$  &   $\mathcal{O}\left(\sqrt{\frac{\kappa_{\ell}}{1-\rho}}\log\frac{1}{\epsilon}\right)$  \\
   
\cline{1-1}\cline{3-4}
 Mudag \cite{ye2023multi} &                   &    $\mathcal{O}\left(\sqrt{\kappa_g}\log\frac{1}{\epsilon}\right)$               &  $\widetilde{\mathcal{O}}\left(\sqrt{\frac{\kappa_g}{1-\rho}}\log\frac{1}{\epsilon}\right)$ \\

 \midrule
DPAG \cite{ye2020decentralized}   & \multirow{3}{*}[-7pt]{   $f$ scvx, $r \neq 0$} &     $\mathcal{O}\left(\sqrt{ \kappa_g }\log\frac{1}{\epsilon}\right)$               & $\widetilde{\mathcal{O}}\left(\sqrt{\frac{\kappa_g}{1-\rho}}\log\frac{1}{\epsilon}\right)$  \\

\cline{1-1}\cline{3-4}
\textbf{Proposed:} DCatalyst-SONATA-L  &                   &   $\widetilde{\mathcal{O}}\left(\sqrt{ \kappa_g }\log\frac{1}{\epsilon}\right)$ &   
 $\widetilde{\mathcal{O}}\left( \sqrt{\frac{\kappa_g}{1-\rho}}\log\frac{1}{\epsilon} \right)$                     \\

 \cline{1-1}\cline{3-4}
\textbf{Proposed:} DCatalyst-SONATA-F   &   & $\widetilde{\mathcal{O}}\left(\sqrt{\frac{L+\beta}{\beta}}\cdot  \frac{\beta}{\mu} \log^2\frac{1}{\epsilon}\right)$ &  $\widetilde{\mathcal{O}}\left(\sqrt{\frac{\beta/\mu}{1-\rho}}\log\frac{1}{\epsilon}\right)$ \\

 \bottomrule 
 
\end{tabular}}\label{table:scvx} 
\end{table*}

\begin{table*}[t!]\footnotesize
\centering
\caption{\scriptsize \textbf{Decentralized accelerated  algorithms for (nonstrongly) convex $f$}: The most efficient decentralized accelerated algorithm for each reported instance of  \eqref{eq:problem}, based either on gradient or communication complexity. Definitions are as in Table~\ref{table:scvx}. For nonsmooth  $u$, $\epsilon$-optimality is measured in terms of the  gradient of the  Moreau envelope of $u$.  
 The proposed framework enables acceleration for nonsmooth $u$ as well as rate dependence on   $L$ (rather than the less favorable $L_{\max}$). }
\resizebox{1.0\columnwidth}{!}{\begin{tabular}{c|c|c|c}
\toprule 
\textbf{Algorithm}& \textbf{Problem}  &  \textbf{$\#$ (Proximal) gradient}& \textbf{$\#$ Communications} \\
 
  \midrule 
 APM-C \cite{li2020decentralized} &  $f$ cvx, $r\equiv 0$ &  $\mathcal{O}\left(\sqrt{\frac{L_{\max}}{\epsilon}}\right)$  &   $\mathcal{O}\left(\sqrt{\frac{L_{\max}}{(1-\rho)\epsilon}}\log\frac{1}{\epsilon}\right)$ \\

 \midrule 
  APM \cite{li2020decentralized}   &   \multirow{3}{*}[-8pt] {\centering $f$ cvx, $r\neq 0$} & $\mathcal{O}\left(\frac{L_{\max}}{\epsilon}\right)$  & $\mathcal{O}\left(\frac{L_{\max}}{\epsilon\sqrt{1-\rho}}\right)$  \\

\cline{1-1}\cline{3-4} 
 {\bf Proposed: }  DCatalyst-SONATA-L  &                      &   $\widetilde{\mathcal{O}}\left(\sqrt{\frac{L}{\epsilon }}\log\frac{1}{\epsilon}\right)$ 
 
 &   $\widetilde{\mathcal{O}}\left(\sqrt{\frac{L}{(1-\rho)\epsilon}}\log\frac{1}{\epsilon}\right)    $  \\
\cline{1-1} \cline{3-4}
{\bf Proposed: } DCatalyst-SONATA-F &                   &  $\widetilde{\mathcal{O}}\left(\sqrt{\frac{L+\beta}{ \epsilon}} \log^2\frac{1}{\epsilon}\right)$  
  &   $\widetilde{\mathcal{O}}\left(\sqrt{\frac{{\beta}}{(1-\rho)\epsilon}}\log\frac{1}{\epsilon}\right)$  \\
  \bottomrule
\end{tabular}} \label{table:cvx} 
\end{table*}

\begin{table*}[t!]\Large
\centering
\caption{\scriptsize  \textbf{Decentralized accelerated algorithms for finite-sum $f_i$'s}: Current decentralized accelerated VR algorithms   for   \eqref{eq:problem} with $f_i$'s having a finite-sum structure. 
Same quantities as defined in   Table~\ref{table:scvx}; in addition,  $\tilde{\kappa}_s := \left(\max_{ij}L_{ij}\right)/\mu_{\min}$, $\kappa_s := (\max_{i} ({1}/{n})\sum_{j=1}^n L_{ij})/\mu_{\min}$, and  $b\in\{1,\cdots,n\}$  represents  the batch-size used in the (stochastic) estimation of each $\nabla f_i$. Here,   $L_{ij}$ is smoothness parameter of $f_{ij}$.  
Notice that  $\kappa_{\ell} \leq \kappa_s \leq n\kappa_{\ell}$ and  $\kappa_s \leq \tilde{\kappa}_{s}   \leq mn \kappa_s $    (details in Sec.~\ref{section:prelim}).   The proposed framework achieves acceleration  also when $u$ is nonsmooth.}   
\resizebox{1.0\columnwidth}{!}{\begin{tabular}{c|c|c|c}
\toprule 
 
\textbf{Algorithm}& \textbf{Problem} &  \textbf{$\#$ (Proximal) Gradient}  & \textbf{$\#$ Communications} \\
 \midrule 
ADFS\cite{hendrikx2021anOptimal}  & \multirow{4}{*}[-7pt]{$f_i$ scvx, $r\equiv 0$} &  NA & \multirow{2}{*}{ $\mathcal{O}\left(\sqrt{\frac{\kappa_{\ell}}{1-\rho}}\log\frac{1}{\epsilon}\right)$} \\

\cline{1-1}\cline{3-3}
Acc-VR-EXTRA-CA \cite{li2022variance}  &                   & $\mathcal{O}\left(\left(\sqrt{n\kappa_s}+n\right)\log\frac{1}{\epsilon}\right)$ &                   \\

\cline{1-1}\cline{3-4}
DVR-Catalyst \cite{hendrikx2020dual}  &                   &  $\widetilde{\mathcal{O}}\left((\sqrt{n\kappa_s}+n)\log\frac{1}{\epsilon}\right)$  &   $\widetilde{\mathcal{O}}\left(\sqrt{\frac{\kappa_{\ell}}{1-\rho}}\sqrt{\frac{n \kappa_{\ell}}{\kappa_s}}\log\frac{1}{\epsilon}\right)$                   \\

 \cline{1-1}\cline{3-4}
{\bf Proposed:} DCatalyst-VR-EXTRA &        &  $\widetilde{\mathcal{O}}\left(\sqrt{\kappa_{\ell} nb}\log\frac{1}{\epsilon}\right)$ &  $\widetilde{\mathcal{O}}\left(\sqrt{\frac{\kappa_{\ell}}{1-\rho}}\sqrt{\frac{n}{b}}\log\frac{1}{\epsilon}\right)$  \\
 \midrule
PMGT-LSVRG \cite{ye2021pmgtvr} & \multirow{2}{*}[-5pt]{$f$ scvx, $r\neq 0$} & $\mathcal{O}\left((n+\tilde{\kappa}_s)\log\frac{1}{\epsilon}\right)$& 
$\mathcal{O}\left( \frac{ \tilde{\kappa}_s\log\tilde{\kappa}_s+n\log n }{\sqrt{1-\rho}}\log\frac{1}{\epsilon}\right)$                  \\

 \cline{1-1}\cline{3-4}
{\bf Proposed:}  DCatalyst-PMGT-LSVRG &                   &  $
\widetilde{\mathcal{O}}\left(\sqrt{\kappa_s n } \log\frac{1}{\epsilon}\right)  
$
& 
$
 \widetilde{\mathcal{O}}\left(\sqrt{\frac{\kappa_s n}{1-\rho}} \log\frac{1}{\epsilon}\right)$                       \\
\bottomrule

\end{tabular}} 
\label{table:finite-sum} 
\end{table*}

\subsection{Main Contributions}\label{sec:major-contributions} Our   technical contributions can be summarized as follows. 
 
 \textbf{(i)  DCatalyst--A unified decentralized acceleration framework:}
    We introduce a unified black-box framework engineered to  integrate   Nesterov acceleration into  {\it decentralized} optimization algorithms,  thereby broadening their applicability to   effectively tackle the full spectrum of   Problem \eqref{eq:problem}. DCatalyst ensures compatibility with the majority of  decentralized optimization algorithms--specifically capable of linear convergence in solving strongly convex instances of \eqref{eq:problem}.  
    At its core,   DCatalyst operates as an \textit{inexact}, \textit{momentum-accelerated} proximal method (forming the outer loop) that seamlessly incorporates any selected decentralized algorithm (as the inner loop). This chosen inner algorithm  solves {\it approximately}   the auxiliary, proximal subproblems related to \eqref{eq:problem}, 
    adhering to a newly developed {\it criterion of inexactness} and {\it warm-start strategy}   tailored for decentralized settings. 
    
    We prove that DCatalyst attains  optimal communication and/or computational complexity (up to log-factors) across a broad spectrum of embedded decentralized algorithms  and  instances of \eqref{eq:problem}.   
Our numerical results corroborate the theoretical advancements posited by DCatalyst, demonstrating tangible acceleration in all simulated scenarios, with significant performance improvement  with respect to  nonaccelerated schemes, especially when dealing with  ill-conditioned problems. Notably, DCatalyst exhibits superior performance when benchmarked against existing {\it single-loop} decentralized, accelerated algorithms specifically designed for particular subclasses of Problem~\eqref{eq:problem}, showcasing its robustness and versatility across various optimization scenarios. \smallskip

     \textbf{(ii) New decentralized accelerated algorithms for unsolved instances of \eqref{eq:problem}:} The framework's dual-loop architecture is pivotal for enabling acceleration while retaining the essential attributes of the chosen inner-loop decentralized algorithm. This feature yields  {\it new decentralized accelerated} algorithms capable of solving  Problem~\eqref{eq:problem} instances that were hitherto unapproachable. Few notable examples follow. \begin{itemize}\item \textit{Leveraging function  similarity:}  In the case of strongly convex functions $u$ (possibly with $r\neq 0$) satisfying   (\ref{eq:similarity}),  DCatalyst, when integrated with the SONATA-F decentralized algorithm \cite{sun2020distributed}, which harnesses function similarity through local preconditioning for enhanced convergence,  achieves   for the first time near optimal (accelerated) communication complexity (Table 1),  $\widetilde{\mathcal{O}}((1/\sqrt{1-\rho})\cdot\sqrt{\beta/\mu}\cdot \log(1/\varepsilon))$. This improves over  the   $\mathcal{O}({\beta/\mu})$ nonaccelerated 
 scaling  observed with  SONATA-F.   \item \textit{Acceleration for nonstrongly convex, composite functions $u$:} DCatalyst facilitates acceleration (while possibly exploring function similarity) for such a class of problems. It adeptly extends any decentralized algorithm, originally achieving linear convergence for strongly convex, composite functions but lacking assurances for merely convex scenarios, into the accelerated domain. This   gives rise to new algorithms such as DCatalyst-SONATA-L   and DCatalyst-SONATA-F (Table 2). 
 \item \textit{Acceleration for composite $u$  with finite-sum structure $f_i$'s:} DCatalyst serves as an ideal platform for integrating VR techniques within decentralized schemes, thereby facilitating acceleration in particular  for composite functions $u$ whose $f_i$'s have a finite-sum structure. For instance, hinging on the PMGT-LSVRG   decentralized algorithm   {\cite{ye2021pmgtvr}}, we propose DCatalyst-PMGT-LSVRG,  which achieves for the first time acceleration in the aforementioned setting on both computation and communication complexities. This is particularly crucial for tackling ill-conditioned problems, where the number of data points $n$ is significantly smaller than the stochastic condition number     $\kappa_s$   (see Table 3).
 \end{itemize}

    \textbf{(iii)   New analysis: the notion of {\it \textbf{inexact}} estimating sequences:} 
Our convergence analysis   hinges on the notion of {\it inexact estimating sequences}, a  novel  extension of the well-known estimating sequences \cite{nesterov2013introductory}, tailored to  the minimization of {\it composite} losses in {\it decentralized} settings. This machinery rigorously incorporates consensus errors and inexact solutions of decentralized subproblems for each agent--challenges that traditional estimating sequence models,  such as   \cite{nesterov2013introductory,baes2009estimate,lin2015universal,lin2018catalyst,Alexandre2021AccMethods}, cannot directly incorporate. It provides  a unified framework for both the convergence analysis and tuning   of decentralized optimization algorithms, employing acceleration, and  represents a standalone contribution of independent interest.

    Guided by this novel framework, we have  established   new termination criteria for  inner-loop   decentralized (non-accelerated) algorithms, ensuring overall convergence at accelerated rates.  Unlike the traditional Catalyst framework in the centralized setting \cite{guler1992new,lin2015universal, lin2018catalyst, Alexandre2021AccMethods}, or its application to     decentralized algorithms, Catalyst-EXTRA \cite{li2020revisiting} and DVR-Catalyst  \cite{hendrikx2020dual}, the proposed framework   seamlessly integrates with composite (and hence nonsmooth) losses and supports the use of generic stochastic decentralized algorithms within the inner loop. This adaptability makes it uniquely suited to a broader range of applications.

\subsection{Notation and paper structure}

We use the following notation. Let $[m] := \{1,\cdots,m\}$,  with  $m$ being any positive integer.  
  We use $\norm{.}_p$ to denote the $\ell_p$ norm; $\norm{.}$ is by  default the $\ell_2$-norm when applied to vectors and the   Frobenius norm when we input   matrices.   We denote by $1_m$ the the $m$-dimensional (column) vector of all ones.  Bold letters is used to denote matrices resulting from stacking a set of given vectors row-wise; for example, given $x_1,\ldots, x_m\in \mathbb R^d$, we define $\vx := [x_1,\cdots, x_m]^{\top} \in \mathbb{R}^{m\times d}$ and $\bx := (1/m)\summ x_i$. For a square matrix $A$, we use $\sigma_{\max}(A)$ and  $\sigma_{\min}^{+}(A)$ respectively to denote its maximum eigenvalue and smallest positive eigenvalue.
  For a sequence of sets $\{B_k\}_{k=0}^{\infty}$, we use $\prod_{k=0}^{\infty}B_k$ to denote their Cartesian product. 
  
 Given $\eta>0$   and a proper, closed convex function $r:\mathbb{R}^d\rightarrow \mathbb{R}\cup\{+\infty\}$, the  Moreau envelope of $r$, $\texttt{M}_{\eta\, r}:\mathbb{R}^{d}\rightarrow\mathbb{R}$,   and the proximal operator, $\texttt{prox}_{\eta r}(x):\mathbb{R}^d\to \mathbb R^d$, are 
 $$\texttt{M}_{\eta\,r}(x) :=
 \min_{y\in\mathbb{R}^d} \eta\, r(y)+\frac{1}{2}\norm{y-x}^2\quad \text{and}\quad \texttt{prox}_{\eta\, r}(x):=  \arg\min_{y\in\mathbb{R}^d}  \eta \,r(z)+\frac{1}{2}\norm{y-x}^2,$$  respectively.  
Notice that  
 $\texttt{M}_{\eta\, r}$ is differentiable, with gradient  

$$
\nabla \texttt{M}_{\eta\, r}(x) = \frac{1}{\eta}\left(x-\texttt{prox}_{\eta\,r}(x)\right),  $$
which is thus $1/\eta$-Lipschitz continuous.  
Denote  by $\texttt{M}_{\eta\,r}^\star$  the minimum of $\texttt{M}_{\eta\,r}$ over $\texttt{dom}\,r$.
With a slight abuse of notation, given $\vx\in\mathbb{R}^{m\times d}$, we will write   $\vy :=$ {\rm $\texttt{prox}_{\eta r}$}$\left(\vx\right)\in\mathbb{R}^{m\times d}$, to denote in short row-wise   equalities, that is,   $y_i =\texttt{prox}_{\eta r}(x_i)$, for all    $i\in[m]$.  
 
The rest  of the paper is organized as follows.  
In Sec.~\ref{section:prelim}, we introduce the main assumptions underlying Problem~\eqref{eq:problem};   Sec.~\ref{section:catalyst-framework} formally introduces the proposed  framework, decentralized Catalyst, along with its convergence properties for strongly convex and convex losses.   Sec.~\ref{sec:Applications} applies the proposed framework to a variety of decentralized algorithms. Sec.~\ref{section:analysis} provides  the proofs of the main results. Finally, some numerical results validating the effectiveness of our framework  are discussed  in Sec.~\ref{section:simulations}.

\section{Assumptions and Preliminaries}\label{section:prelim}

  We study \eqref{eq:problem} 
  over a mesh network. The network is modeled as  a static undirected graph $\mathcal{G} := (\mathcal{V},\mathcal{E})$, where $\mathcal{V} :=  [m]$ denotes the set of vertices (agents) and $\mathcal{E} := \{(i,j) | i,j\in\mathcal{V}\}$ is the set of communication links: agents $i$ and $j$ can communicate if $(i,j)\in \mathcal E$. We make the blanket assumption that   $\mathcal{G}$ is connected, and $(i,i)\in\mathcal{E}$ for all $i\in [m]$.

  The assumptions on Problem~\eqref{eq:problem} are  the following.   
\begin{assumption}\label{assump:class}
For Problem~\eqref{eq:problem}, the following hold:  
\begin{enumerate} 
  \item[(i)]   Each $f_i: \mathbb{R}^{d}\rightarrow \mathbb{R}$ is continuously differentiable, $\mu_i$-strongly convex and $L_{i}$-smooth, with $0\leq \mu_i \leq L_i < \infty$;
    \item[(ii)] $f: \mathbb{R}^{d}\rightarrow \mathbb{R}$ is   $\mu$-strongly convex and $L$-smooth, with $0  \leq  \mu \leq L < \infty$;
      \item[(iii)]  $r:\mathbb{R}^d\to \mathbb{R}\cup\{+\infty\}$ is a proper, closed, and convex function.
\end{enumerate} 
\end{assumption}

This assumption encompasses both    strongly convex   or convex functions $f$, with the latter corresponding to  $\mu=0$. Associated to Assumption~\ref{assump:class}, there are the following quantities: 
 
 $$
L_{\max} := \max_{i\in [m]} L_i,\quad {\mu_{\min}} := \min_{i\in [m]} \mu_i,\quad \kappa_g := \frac{L}{\mu},\quad {\kappa_{\ell}} := \frac{{L_{\max}}}{\mu_{\min}}.
$$

  Here,   by writing $\kappa_g$ and $\kappa_{\ell}$ we implicitly assume $\mu>0$ and $\mu_{\min}>0$, respectively.   These quantities are   pivotal in assessing the efficiency of  decentralized algorithms. It is noteworthy that   $\kappa_g\leq \kappa_{\ell}$. The gap can be significantly large, as shown in the following example. \smallskip

 \begin{example}
     
 Let $f_1, f_2:\mathbb{R}^2\rightarrow\mathbb{R}$ such that 
 $$f_1(x) = \frac{1}{2} x^{\top}
 \left(\begin{matrix}
    1 & 0\\
    0 & \epsilon
 \end{matrix}\right)x, \quad f_2(x) = \frac{1}{2} x^{\top}
 \left(\begin{matrix}
    \epsilon & 0\\
    0 & 1
 \end{matrix}\right)x.$$ Hence, $f(x) =  ((1+\epsilon)/{4})\norm{x}^2.$ It follows that $\kappa_g = 1$ while  $\kappa_{\ell} =  \epsilon^{-1}$, implying that  $\kappa_{\ell}$ can be arbitrarily large.$\hfill \square$
 \end{example} \smallskip

 When dealing with ERM  problems with $f$ in the form \eqref{eq:f_ij}, we further postulate  the following standard property for the $f_{ij}$'s.     
\begin{assumption}\label{assump:class-finite-sum}
Each $f_{ij}:\mathbb{R}^d\to \mathbb{R}$ in (\ref{eq:f_ij})  is continuously differentiable, convex, and $L_{ij}$-smooth, with $0<L_{ij}<\infty$.  
\end{assumption}

In the context of Assumption~\ref{assump:class-finite-sum}, we define several critical quantities that are instrumental in the analysis of first-order methods solving in  particular ERM problems:  
 $$
\bar{L}_i := \frac{1}{n}\sum_{j=1}^n L_{ij},\quad \bar{L}_{\max} := \max_{i\in [m]}\bar{L}_i,\quad \kappa_s = \frac{\bar{L}_{\max}}{ \mu_{\min}}.  
$$
Since $L_i \leq \bar{L}_i \leq n L_i$, we have \begin{equation}\label{eq:stochastic_kappa_rel}\kappa_{\ell} \leq \kappa_s \leq n\kappa_{\ell}\quad \text{and}\quad  L_{\max}\leq \bar{L}_{\max} \leq n L_{\max}.\end{equation}  
Here, $\kappa_s$  represents the  stochastic condition number (in contrast to the ``batch'' condition numbers, $\kappa_g$ and $\kappa_\ell$), which appears in the convergence rates of   VR-based algorithms for finite-sum optimization problems. The relationship (\ref{eq:stochastic_kappa_rel}) implies that, in the worst-case scenario where all $f_{ij}$'s are independent, VR-based methods do not offer any advantage over their batch counterparts. However, ERM problems featuring correlated data samples often exhibit $\kappa_s \ll n\kappa_{\ell}$, indicating a potential for significant improvements in the computational efficiency of VR-methods compared to batch algorithms.

 In the   exploration of  ERM  problems, beyond the settings outlined in Assumptions \ref{assump:class} and \ref{assump:class-finite-sum}, we anticipated in (\ref{eq:similarity}) additional structure   in loss functions $f_i$, termed  function similarity. We   encapsulate this property in the subsequent assumption for ease of reference. 
 \begin{assumption}\label{assump:class-similarity}
Each $f_i$ in Assumption~\ref{assump:class}, $i\in [m]$,  satisfies: 
$$
   \big\|\big(\nabla f(x) - \nabla f_i(x)\big) - \big(\nabla f(y) - \nabla f_i(y)\big)\big\|\leq \beta \|x-y\|,\quad \forall x,y\in \texttt{dom}\, r,
 $$
for some  $\beta > 0$. It is assumed $\beta>\mu$ without loss of generality. 
\end{assumption}
The smaller $\beta$, the more similar $f_i$'s. Under Assumption 1, the following bounds holds for $\beta$: $\beta\leq \max_{i\in [m]}\max\big\{|L_i-\mu|,|\mu-L_i|\big\}$. The example below demonstrates the typical, more favorable scaling of $\beta$ in ERM scenarios, due to statistical similarity. 

\textbf{Example:} Consider the ERM setting (\ref{eq:f_ij}), where data are  i.i.d., locally and across the agents. Assume $\nabla^2 f_i$ is $M$-Lipschitz and under extra (mild) assumptions--see \cite{Zhang-Xiao-chapter18}--the following bounds hold for $\beta$, with high probability:  
$$\beta =   \widetilde{\mathcal{O}}\left(\sqrt{\frac{L^2}{n}}\right),\,\,\,\text{if}\,\,\,\text{$M = 0$ ($f_i$'s are quadratic)};\quad \text{or} \quad \beta =\widetilde{\mathcal{O}}\left(\sqrt{\frac{L^2 d}{n}}\right), \,\,\,\text{otherwise}. 
$$ It turns out that  
$$
\frac{\beta}{\mu} \bigg/ \frac{L}{\mu} = \min\left\{1, \widetilde{\mathcal{O}}\left(\sqrt{\frac{d}{n}}\right)\right\}.
$$
Hence, the ratio ${\beta}/{\mu}$ may be considerably smaller than $\kappa_g$, particularly in cases where $n$ is large, a common condition in many applications. This indicates potential for improved efficiency in exploiting function similarity within decentralized optimization.

 \section{The Decentralized Catalyst Framework}\label{section:catalyst-framework}  DCatalyst inputs a given  decentralized algorithm  and enhances it  through a black-box acceleration procedure.   Our initial step involves    identifying a spectrum of candidate algorithms via a broad set of assumptions, which   encapsulate the vast majority of current decentralized  schemes--see Sec.~\ref{sec:dec-alg}.    Sec.~\ref{sec:Catalyst-framework} proceeds to  introduce the black-box procedure enabling acceleration in the aforementioned decentralized schemes.   Convergence analysis of DCatalyst is presented in  Sec.~\ref{sec:convergence-scvx} (strongly convex $u$)  and Sec.~\ref{sec:convergence-cvx} (convex $u$).

\subsection{Distributed algorithms}\label{sec:dec-alg}
  We model  decentralized  algorithms through    operators that iteratively act on   variables held by the agents. To this end, let us  denote by  $s_i = [x_i^\top, y_i^\top]^\top$ the set of variables controlled by agent $i$, where $x_i\in\mathbb{R}^d$ is the local  estimate of the optimization variable $x$ of the problem under consideration, and   $y_i\in\mathbb{R}^{d^\prime}$ serves as some auxiliary  vector-variable,   instrumental to track locally some additional information used by the algorithm to generate agents' updates, such as  dual (e.g., \cite{alghunaim2020decentralized,ling2015dlm,shi2015extra})  or gradient-tracking estimates (e.g., \cite{di2016next,DIGing,qu2017harnessing,sun2020distributed}). Denote   by $s_i^t$ the values  of these variables at iteration $t$, and  by $\vs^{t} := (s_i^t)_{i\in[m]} $   the stack of the agents' tuples at that iteration. The decentralized algorithm under consideration is then modeled as a (possibly random) mapping  $\cM:\mathbb{R}^{m\times(d+d^\prime)  }\rightarrow \mathbb{R}^{m\times(d+d^\prime)}$: given   $\vs^t$ as input at iteration $t$, the decentralized algorithm $\cM$ produced    $\vs^{t+1} = \cM(\vs^{t})$ as    new iterate at time $t+1$.
    
We focus on accelerating decentralized   algorithms $\mathcal M$ enjoying the following   properties. 
\begin{assumption}\label{Lyapunov} Let the decentralized  algorithm $\cM$ be employed to solve a strongly convex instance of  Problem~\eqref{eq:problem}  (that is,  under Assumption~\ref{assump:class}, with $\mu > 0$)  over a mesh network. 
Let      $\{\mathbf{s}^{t}\}_{t}$ be  the sequence of state variables generated by $\cM$,   starting from a given, feasible $\mathbf{s}^{0}$. There exists  a merit  function   $\cL: \mathbf{s}\mapsto \cL(\mathbf{s})\in \mathbb{R}_{+}$  associated with $\cM$,    such that   
 \begin{enumerate} 
 \item[(i)] (\textit{linear convergence}):   
\begin{equation}\label{assum3}
    \Eb[\cL(\vs^{t})|\vs^0] \leq c_{\cL}\left(1-\frac{1}{r_{\cM}}\right)^{t}\cL(\vs^{0}),  \quad \forall t\geq 0,
\end{equation}
where $c_{\cL} > 0$ and $r_{\cM} > 1$ are  problem/algorithm-dependent constants; and    \item[(ii)] (\textit{optimality}):
    Given any $\vs = [\vx^\top, \vy^\top]^\top\in\mathbb{R}^{m\times(d+d^\prime)}$, 
 $$
    \frac{1}{m}\|\vx -\vx^{\star}\|^{2} \leq  \cL(\vs) , 
$$
    where $\vx^{\star} := 1_{m}(x^{\star})^{\top}$ and $x^{\star}$ is the unique solution  of the optimization problem.   
\end{enumerate}
\end{assumption}

Assumption~\ref{Lyapunov} is quite standard in the literature of decentralized optimization algorithms--Sec.~\ref{sec:Applications} specializes (i) and (ii) to the most common  schemes. Specifically, (i) states that, when minimizing a strongly convex loss on a network,   $\cM$ exhibits a linear convergence rate. Condition (ii) guaranteed that convergence of the merit function  $\cL$ implies that of the iterates to the optimal solution $x^\star$ (at the same rate). 
 
 \subsection{\texorpdfstring{Accelerating  $\cM$   via   decentralized Catalyst}{Accelerating  \cM   via   decentralized Catalyst}}\label{sec:Catalyst-framework}
 
 The DCatalyst framework is formally introduced in Algorithm~\ref{alg:framework} and commented next. 
At its core, it can be viewed as the successive application of the chosen decentralized algorithm $\mathcal M$  for the inexact minimization of the function $u^k(x)$,  defined in (\ref{eq:def_u}). The quadratic term $(\delta/2)\|x - z_i^k\|^2$ therein functions as an extrapolation mechanism akin to Nesterov's approach, to gain acceleration. The momentum variables $z_i^k$ are updated locally by the respective  agents at the conclusion of each outer loop (during the Extrapolation step (S.2)), following the   update rule in (\ref{eq:extrapolation-DCat}). 
Crucially, at the start of every inner loop, the decentralized algorithm $\cM$ undergoes a strategic \textit{warm-restart}, and runs  until a predefined \textit{termination criterion} is met. These  strategies are informed directly from our  convergence analysis, with their specifics being dependent on both the optimization problem and the  decentralized algorithm $\cM$  at hand. Details on these criteria along with   tuning  recommendations    are provided in Sec.~\ref{sec:convergence-scvx} and Sec.~\ref{sec:convergence-cvx}  for strongly convex and non-strongly convex functions $u$, respectively.

\textbf{A Bird's-eye view: }   DCatalyst can be viewed as an acceleration of a {\it centralized}, inexact proximal method (the outer loop), where the proximal subproblems are approximately solved in a distributed manner using the decentralized algorithm $\cM$. This setup adheres to a novel criterion of inexactness for solving the subproblems \eqref{eq:subproblem_DCat} distributively, as proposed in this paper.   
Assuming  one can absorb consensus errors on  the agents' variables $x_i$'s and  momentum vectors $z_i$'s into this   criterion of solution approximation (to be introduced), we can approximate \texttt{(S.1)} and \texttt{(S.2)} as  
  \begin{align*}
    &  x_i^k \approx \overline{x}^k:= \frac{1}{m}\sum_{i=1}^m x_i^k\quad \text{and}\quad  z_i^k \approx \bar{z}^k:=\frac{1}{m}\sum_{i=1}^m z_i^k,\\
   & \texttt{(S.1)$^\prime$:}\quad \bar{x}^{k+1} \approx   \argmin_{x\in \mathbb{R}^d}   u(x)+   \frac{\delta}{2} \|x-\bar{z}^k\|^2,\\
    & \texttt{(S.2)$^\prime$:}\quad \bar{z}^{k+1} = \bar{x}^{k+1} + \frac{\alpha^{k}(1-\alpha^{k})}{(\alpha^{k})^{2}+\alpha^{k+1}}  \,(\bar{x}^{k+1} - \bar{x}^k),
  \end{align*}
  where we used the fact that the minimization of $u^k(x)$ in (\ref{eq:def_u})  and that of the function in \texttt{(S.1)$^\prime$} have the same solution. 
  The  above dynamics  are a resemble of an inexact proximal acceleration  \cite{lin2015universal, lin2018catalyst, Alexandre2021AccMethods}.
  
 However, the conventional convergence analysis for centralized accelerated methods does not directly apply to DCatalyst due to several critical distinctions. 
    {\bf (i)}  The approximate solution criteria for proximal problems as  discussed in the literature  \cite{lin2015universal, lin2018catalyst, Alexandre2021AccMethods}, are not  applicable here,  due to consensus errors and the challenge of verifying such conditions in a distributed context, especially in the presence of a nonsmooth function $r$.    {\bf  (ii)}  The potential functions used in the aforementioned works are not suitable to    guarantee the accelerated convergence rate of DCatalyst's outer loop. They fail to account for the additional errors inherent in the inexact, distributed minimization of $u^k$ via $\mathcal{M}$.  
 Our novel approach   addresses exactly  these challenges.

\begin{algorithm}[!ht]
\caption{\texttt{Decentralized Catalyst}}\label{alg:framework} 
 {\bf Input}:  $x_i^0=z_i^0 =z_i^{-1}\in \texttt{dom}\,\, r,\,\,\,i\in [m]$;   $\delta>0$; \\ $K\in \mathbb{N}_{++}$ ($\#$ of outer loops), $ \{\alpha^{k}\in (0,1)\}_{k=0}^K$;   Decentralized algorithm $\cM$. \medskip 
   
 {\bf Output}: $\vx^K=(x_i^K)_{i\in [m]}$ \smallskip \\
\textbf{for} {$k=0,1,2,\ldots,K-1$} \textbf{do}

\texttt{Set}  
\begin{equation}\label{eq:def_u}
    u^{k}(x) = \frac{1}{m}\sum_{i=1}^m  f_i^k(x) + r(x),\quad \text{with}\quad f_i^k(x) = f_i(x) + \frac{\delta}{2}\|x-z_i^k\|^2;
\end{equation}  
 
 \texttt{(S.1) Minimization step:} Minimize approximately $u^{k}$ over the network $\Gh$  running     $\cM$, with proper initialization $s_i^k=[(x_i^k)^\top,(y_i^k)^\top]^\top$, $i\in [m]$, and up to a suitable   termination: 
    \begin{align}\label{eq:subproblem_DCat} 
     x_{i}^{k+1}     \approx_{\mathcal{M}}  \arg\min_{x\in\mathbb{R}^{d}}\,u^{k}(x),\quad \forall i\in [m];  \tag{P$^k$} 
    \end{align}

 \texttt{(S.2)     Extrapolation step:}
\begin{equation}\label{eq:extrapolation-DCat} z_i^{k+1} = x_i^{k+1} +\frac{\alpha^{k}(1-\alpha^{k})}{(\alpha^{k})^{2}+\alpha^{k+1}} \,\left(x_i^{k+1} - x_i^k\right),\quad \forall i \in[m].
\end{equation}
 
\textbf{end for}
\end{algorithm}

 \subsection{Convergence of DCatalyst: Strongly convex losses}\label{sec:convergence-scvx}

 In this section,   we delve into the convergence properties of DCatalyst under the condition that  the objective function $u$ is  strongly convex. Our analysis  begins assessing the      convergence dynamics of DCatalyst's outer loop (Proposition \ref{prop:outer-loop-scvx}). Here,     a suitable warm-restart mechanism and  termination criterion for $\cM$ are provided, which maintain the inexactness of solving the proximal subproblems in (S.1)  within the threshold ensuring  convergence of the outer loop.   Our second result  outlines specific tuning recommendations for the decentralized algorithm $\cM$ used in the inner loop, designed   to meet the   established termination criterion (Proposition~\ref{prop:inner-loop-scvx}). Finally, combining    both the outer- and inner-loop convergence, we obtain  the overall computation/communication complexity of  DCatalyst (Theorem~\ref{theorem:total-iter-scvx}). The proofs of the results in this section can be found in Sec.~\ref{section:analysis}.
  
\smallskip

\textbf{Notation and probability space:} 
 We denote  the state variables generated by the decentralized algorithm $\cM$ at the inner iteration $t$ of  the outer iteration $k$ by $\vs^{k,t}$. For the subproblem \eqref{eq:subproblem_DCat}, we denote its solution as $x^{k,\star}$ and $\vx^{k,\star} := 1_m(x^{k,\star})^{\top}$. 
 When solving the subproblem \eqref{eq:subproblem_DCat}  using the chosen algorithm $\mathcal M$,   we refer to  $\cL^k:\mathbb{R}^{m\times(d+d^{\prime}) } \rightarrow \mathbb{R}_{+}$ as   the merit function associated with $\cM$ when applied to \eqref{eq:subproblem_DCat} and satisfying   Assumption~\ref{Lyapunov}.  Also, with a slight abuse of notation,   we will write $r_{\cM,\delta}$ instead of  $r_{\cM}$ in (\ref{assum3}),  to explicitly highlight the dependency of the convergence rate on the design  parameter  $\delta$.

When $\cM$ is a randomized algorithm, the   sequences generated by Algorithm~\ref{alg:framework} are random, and the following   probability space is subsumed.
Define $T_k\in (0,\infty)$ as the number of (inner) iterations performed by the decentralized algorithm $\cM$ at the outer iteration $k$. Considering  $\{\{\vs^{k,t}\}_{t=0}^{T_k}\}_{k\geq 0} \in \prod_{0}^{\infty}\mathbb{R}^{m\times(d+d^{\prime})}$  as the sequence of random variables generated by Algorithm~\ref{alg:framework},  let  $\Omega$  denote  the sample space of all the sample paths $\omega  := \{\{\uwave{\vs}^{k,t}\}_{t=0}^{T_k}\}_{k\geq 0} \in \prod_{0}^{\infty}\mathbb{R}^{m\times(d+d^{\prime})}$, where $\uwave{\vs}^{k,t}$ is a realization of $\vs^{k,t}$. For each $k$, define  $\Omega_{k} := \prod_{0}^{T_k}\mathbb{R}^{m\times(d+d^{\prime})}$, and let $\mathcal{B}_k$ be the  Borel $\sigma$-algebra of $\prod_{\ell=0}^{k}\Omega_{\ell}$. 
Suppose $\cM$ induces  a sequence of  Borel probability measures  $\{P_k\}_{k=0}^{\infty
}$, where each $P_k$ is defined on $\mathcal{B}_k$, satisfying the  property:   for each  $E\in \mathcal{B}_k$ and any $n\geq0$, $P_k(E) = P_{k+n}\left(E\times \prod_{\ell=k+1}^{k+n}\Omega_{\ell}\right)$. By the Ionescu–Tulcea’s theorem \cite{Klenke2008ProbabilityTA}, there exists a unique probability measure $\mathbb{P}$ on $\mathcal{\cF} := \sigma\left(\cup_{k=0}^{\infty}\cup_{E_k\in\mathcal{B}_k} \left(E_k\times\prod_{\ell=k+1}^{\infty}\Omega_{\ell}\right)\right)$, which agrees with all $P_k$'s:  for any $E\in\mathcal{B}_k$, $ \mathbb{P}\left(E\times\prod_{\ell=k+1}^{\infty}\Omega_{\ell}\right) = P_k(E)$. The probability space of interest is then   $\left(\Omega, \mathbb{P},\cF\right)$. Further, we define the filtration $\{\cF_k\}_{k\geq 0}$, where   $\cF_{k+1} := \sigma\left(\{\{\vs^{\ell,t}\}_{t=0}^{ T_{\ell}}\}_{\ell=0}^{k}\right)= \mathcal{B}_k\times \prod_{\ell=k+1}^{\infty}\Omega_{\ell}$, and   $\cF_0 := \{\emptyset,\Omega\}$. 
All inequalities involving conditional expectations with respect to  $\cF_k$   are understood to hold almost surely, though this specification is  omitted for brevity.

    Our first result is the  convergence of outer loop of  DCatalyst, as given next.   
\begin{proposition}[Convergence of the outer loop]\label{prop:outer-loop-scvx}
Consider  Problem \eqref{eq:problem} under Assumption~\ref{assump:class} with $\mu > 0$. Let $\{\mathbf{x}^k\}$ be the sequence generated by Algorithm~\ref{alg:framework} in  the following setting:  {\bf (i)} in Step  \texttt{(S.1)}, each subproblem \eqref{eq:subproblem_DCat} is (approximately) solved by  a decentralized algorithm $\cM$, which   satisfies Assumption~\ref{Lyapunov}, with merit function denoted by $\cL^{k}:\mathbb{R}^{m\times(d+d^{\prime}) } \rightarrow \mathbb{R}_{+}$; {\bf (ii)} in  Step $(\texttt{S.2})$,   $\alpha^{k}$  is chosen as   $\alpha^{k} =  \alpha := \sqrt{\mu/(\mu+\delta)}$, for all $k\geq 0$; {\bf (iii) }     for each outer loop $k\geq 0$,   $\cM$  is terminated after     $T_k$ (inner) iterations,   such that 
\begin{equation}\label{ineq:prop-out-scvx-Tk}
    \mathbb{E}[\cL^{k}(\vs^{k,T_k}) ] \leq \epsilon_0(1-c\alpha)^{k+1},   
\end{equation}
where $c\in(0,1)$ is some universal constant, and   $\epsilon_{0}>0$ is arbitrarily chosen; and {\bf (iv)} the warm-restart    
    $\vx^{k} := \vx^{k-1,T_{k-1}}$ is used. 
    
    Then,  for all  $k\geq 0$,  
\begin{equation}\label{Prop1:result-scvx}
    \begin{aligned}
        &\overm\Eb[\|\vx^{k}-\vx^{\star}\|^{2}] \leq c_{scvx}\left(1-c\alpha\right)^{k},  \,\text{with}\,\, c_{scvx}=\mathcal{O}\left(\left(\frac{\delta^2}{\mu^2}+1\right)(\|\vx^{\star}-\vx^0\|^2+ \epsilon_0)\right).
    \end{aligned}
\end{equation} 

 The explicit expression of  $c_{scvx}$  
can be found in the proof [see  
  eq.~\eqref{cal:c_scvx},  Sec.~\ref{section:analysis}].  
\end{proposition}
Notice that, under  Assumption~\ref{Lyapunov},  $T_k$   is well-defined for any $k\geq 0$, that is, $T_k < \infty$.  

Proposition~\ref{prop:outer-loop-scvx} establishes linear convergence of the outer loop:   For any given  $\ep>0$,   $  (1/m)\Eb [\|\vx^{k}-\vx^{\star}\|^{2}] \leq \ep$, in $k\geq K$ (outer) iterations, where  

 $$
K  = c\alpha \log\frac{c_{scvx}}{\ep} = \mathcal{O}\left(\sqrt{\frac{\mu+\delta}{\mu}}\log\frac{(1+\delta^2/\mu^2)(\|\vx^{\star}-\vx^0\|^2+\epsilon_0)}{\ep}\right).  $$

The total number of inner-plus-outer iterations, is then 
$$
N_{\texttt{it}} = \sum_{k=0}^K T_k.
$$

We provide next sufficient conditions for $T_k$ to be {\it uniformly} bounded with respect to $k$.

\begin{assumption}\label{assump:warm-start}
Given $\mu > 0$ (resp. $\mu=0)$, instate the setting of Proposition~\ref{prop:outer-loop-scvx} (resp. Proposition~\ref{prop:outer-loop-cvx}). 
The following holds for each $k\geq 0$:  
 
 $$
\cL^{k+1}(\vs^{k+1,0})\leq c_{\cM}\,\cL^{k}(\vs^{k,T_k })+\frac{d_{\cM}}{m}\|\vz^{k}-\vz^{k+1}\|^{2}, 
 $$
where $c_{\cM}, d_{\cM} > 0 $ are some problem-dependent constants independent of $k$. 
\end{assumption}

At a high level, the assumption necessitates that the merit function at the start of each new inner loop does not significantly deviate from its value at the end of the preceding inner loop. In  Sec.~\ref{sec:Applications}, we demonstrate through several examples that this condition can be satisfied by properly designing a warm-start initialization of the state variables (in particular the $y$-variables) in the inner-loop algorithm 
 $\mathcal M$. 
  
   Equipped with  Assumption~\ref{assump:warm-start}, we prove  the following uniform upper bound for $T_k$.

\begin{proposition}\label{prop:inner-loop-scvx} Instate the  setting of Proposition~\ref{prop:outer-loop-scvx}, and  
additionally suppose that  
Assumption $\ref{assump:warm-start}$ holds. Then, the minimal $T_k$ for condition  \eqref{ineq:prop-out-scvx-Tk}   to hold is bounded by   
 \begin{equation*} 
      T_k  \leq \left\lceil r_{\cM,\delta}\cdot \max\left\{   \log\frac{c_{\cL}\cL^{0}(\vs^{0,0})}{\epsilon_{0}(1-c\alpha)},  
 \,\, \log\frac{  c_{\cL} c_{\cM}  + 36 d_{\cM}c_{scvx}\frac{1}{\epsilon_{0}(1-c\alpha)^{2}} }{1-c\alpha}\right\}\right\rceil  =  \widetilde{\mathcal{O}}(r_{\cM,\delta}).
 \end{equation*}

Here,  $r_{\cM,\delta}$ corresponds to   $r_{\cM}$ in (\ref{assum3}) with $\mathcal M$ now applied to minimize $u^k$ in \eqref{eq:def_u} (hence depending also on  $\delta$).
\end{proposition}
 
Decentralized algorithms employed in the inner loop may require multiple rounds of communications/iteration. Denote such a number by  $N_{\texttt{com},\cM}$,    
  the total number  $N_{\texttt{com}}$ of communication rounds  is 
$$
N_{\texttt{com}} = N_{\texttt{it}} \cdot N_{\texttt{com}, \cM}. 
$$

Combining Proposition~\ref{prop:outer-loop-scvx} and Proposition~\ref{prop:inner-loop-scvx}, we obtain the desired result for the  total number of  iterations and communication cost of Algorithm~\ref{alg:framework}.
 
\begin{theorem}\label{theorem:total-iter-scvx} Given    Problem (\ref{eq:problem}) under Assumption~\ref{assump:class} with $\mu > 0$, let $\{\mathbf{x}^k\}$ be the sequence generated by Algorithm~\ref{alg:framework} in  the following setting:  {\bf (i)} in Step  \texttt{(S.1)}, a decentralized algorithm $\cM$ is employed,    satisfying Assumption~\ref{Lyapunov} with merit function    denoted by $\cL^{k}:\mathbb{R}^{m\times(d+d^{\prime}) } \rightarrow \mathbb{R}_{+}$;  {\bf (ii)} in  Step $(\texttt{S.2})$,   $\alpha^{k}$  is chosen as   $\alpha^{k} =  \alpha := \sqrt{\mu/(\mu+\delta)}$, for all $k\geq 0$; {\bf (iii) }     for each outer loop $k\geq 0$,   $\cM$  is terminated after     $T_k$ inner iterations,   such that (\ref{ineq:prop-out-scvx-Tk}) holds with $\cL^{k}$ as defined in (i); and {\bf (iv)} the warm-restart    
    $\vx^{k} := \vx^{k-1,T_{k-1}}$ is used and   Assumption~\ref{assump:warm-start} holds. 
Then,    
    $$ \overm\Eb[\|\vx^{k}-\vx^{\star}\|^{2}]  \leq \ep,$$   after  {$N_{\texttt{it}}$} inner-plus-outer iterations       given by   \begin{equation}\label{eq:scvx_final_rate} N_{\texttt{it}} = \widetilde{\mathcal{O}}\left(r_{\cM,\delta}\,\sqrt{\frac{\mu+\delta}{\mu}}\log\frac{(1+\delta^2/\mu^2)(\|\vx^{\star}-\vx^0\|^2+\epsilon_0)}{\ep}\right).  
    \end{equation}  {The total number of communications is  $N_{\texttt{com}} =  N_{\texttt{it}}\cdot N_{\texttt{com},\cM}$.}  
\end{theorem}

Theorem~\ref{theorem:total-iter-scvx} certifies   linear convergence of  DCatalyst to  the  solution of (\ref{eq:problem}). The iteration/communication complexity, as outlined in (\ref{eq:scvx_final_rate}), is influenced by the tunable  parameter
 $\delta$ and the convergence rate  $r_{\cM,\delta}$ of the chosen  algorithm  $\cM$.   A crucial aspect of this theorem is providing tuning recommendations for $\delta$,      potentially achieving(near-)optimal communication complexity (up to log-factors). Section~\ref{sec:Applications} provides specific choices for  $\delta$   to enhance convergence of a range of existing, non-accelerated, decentralized algorithms.

\subsection{Convergence of DCatalyst: Convex losses}\label{sec:convergence-cvx} 
This section studies  convergence of Algorithm~\ref{alg:framework} when applied to (non-strongly) convex functions $u$. We mirror the structure of Sec.~\ref{sec:convergence-scvx}.

Our first result establishes convergence of the outer loop. 
\begin{proposition}[Convergence of the outer loop]\label{prop:outer-loop-cvx}  
Consider  Problem (\ref{eq:problem}) under Assumption~\ref{assump:class} with $\mu = 0$. Let $\{\mathbf{x}^k\}$ be the sequence generated by Algorithm~\ref{alg:framework} in  the following setting:  {\bf (i)}   in Step  \texttt{(S.1)}, a decentralized algorithm $\cM$   satisfying Assumption~\ref{Lyapunov} with merit function  $\cL^{k}:\mathbb{R}^{m\times(d+d^{\prime}) } \rightarrow \mathbb{R}_{+}$ is employed, to (approximately) solve the subproblem~(\ref{eq:subproblem_DCat}); {\bf (ii)} in  Step $(\texttt{S.2})$,   $\alpha^{k}$  is chosen recursively according to \begin{equation}\label{eq:alpha-rec}(\alpha^{k})^{2} = (1-\alpha^{k})(\alpha^{k-1})^2, \,\,\forall k\geq 1,\quad\text{and}\quad \alpha^{0} = \frac{\sqrt{5}-1}{2};\end{equation}  {\bf (iii) }     for each outer loop $k\geq 0$,   $\cM$  is terminated after     $T_k$  iterations,   such that

\begin{equation}\label{ineq:prop-out-cvx-Tk}
     \Eb[\cL^{k}\big(\vs^{k,T_k}\big)]\leq \left(\frac{1}{k+3}\right)^{4+2r_0}\epsilon_0,  
\end{equation}
where  $r_0 > 0$ is some universal constant, and $\epsilon_{0} > 0$ is arbitrarily chosen; and {\bf (iv)} the warm-restart    
    $\vx^{k} := \vx^{k-1,T_{k-1}}$ is used,  for all $k\geq 1$. 
    
    Then,  for all $k\geq 0$,   
 
 $$
    \frac{1}{2\delta m}\summ\Eb[\|\nabla \Mud(x_i^{k})\|^{2}] \leq \frac{c_{cvx}}{(k+2)^{2}}, \quad \text{with}\quad c_{cvx}=\mathcal{O}(\delta(\|\vx^{\star}-\vx^0\|^2+\epsilon_0)).
 $$ 
The explicit expression of $c_{cvx}$ is given in (\ref{cal:c_cvx}),  Sec.~\ref{section:analysis}.  

 Moreover, if $r\equiv 0$, for all $k\geq0$, we have  
$$   \frac{1}{2\delta m}\summ\Eb[\|\nabla u(x_i^{k})\|^{2}] \leq \frac{c_{cvx}(L+\delta)^2}{\delta^2(k+2)^{2}}. $$
\end{proposition}

 The following result  provides the growth rate of $T_k$.  
\begin{proposition}\label{prop:inner-loop-cvx} Instate the setting of Proposition~\ref{prop:outer-loop-cvx};   suppose that {\bf (i)} Assumption~\ref{assump:warm-start} holds, and {\bf (ii)} the iterates $\{\mathbf{x}^k\}$   are bounded,    $\|x_i^k\| \leq B$, for all $i\in [m]$ and $k\geq 0$, and some $B\in (0,\infty)$.    Then, the minimal   $T_k$ for  (\ref{ineq:prop-out-cvx-Tk})  to hold is bounded:   
 
\begin{equation}\label{inner-lp-cvx}
\begin{aligned}
   T_{k} &\leq \left\lceil r_{\cM,\delta}\cdot\max\left\{  
         \log\frac{ c_{\cL}9^{r_0+2}\cL^{0}(\vs^{0,0})}{\epsilon_0}, 
        \log \left(c_{\cL} c_{\cM}2^{4+2r_{0}} +\frac{36c_{\cL} d_{\cM}B^2(k+3)^{4+2r_{0}}}{\epsilon_{0}}\right)\right\}\right\rceil\\
        & = \widetilde{\mathcal{O}}\left(r_{\cM,\delta}\log k\right). 
\end{aligned}
\end{equation}

\end{proposition}

 Combining Proposition~\ref{prop:outer-loop-cvx} and Proposition~\ref{prop:inner-loop-cvx}, we obtain the  total number of inner-plus-outer iterations and communications to $\epsilon$-optimality.

\begin{theorem}\label{theorem:total-iter-cvx} 
Given    Problem (\ref{eq:problem}) under Assumption~\ref{assump:class} with $\mu = 0$, let $\{\mathbf{x}^k\}$ be the sequence generated by Algorithm~\ref{alg:framework} in  the following setting:  {\bf (i)} in Step  \texttt{(S.1)}, a decentralized algorithm $\cM$   satisfying Assumption~\ref{Lyapunov} is employed; {\bf (ii)} in  Step $(\texttt{S.2})$,   $\alpha^{k}$  is chosen recursively according to (\ref{eq:alpha-rec}); {\bf (iii) }     for each outer loop $k\geq 0$,   $\cM$  is terminated after     $T_k$ inner iterations,   such that  (\ref{ineq:prop-out-cvx-Tk}) holds; and {\bf (iv)} the warm-restart    
    $\vx^{k} := \vx^{k-1,T_{k-1}}$ is used, Assumption~\ref{assump:warm-start} holds, and in addition $\{\vx^k\}$ are bounded. 
Then, for all $k\geq0$,   
$$
 \frac{1}{2\delta m}\summ\Eb[\|\nabla \texttt{M}_{\frac{1}{\delta}\,r}(x_i^{k})\|^{2}] \leq \epsilon,
$$ 
or when $r \equiv 0$,  
\begin{equation*}
    \frac{\delta}{2  m(L+\delta)^2}\summ\Eb[\|\nabla u(x_i^{k})\|^{2}] \leq  \epsilon, 
\end{equation*}  
after $N_{\texttt{it}}$ inner-plus-outer iterations  given by $$N_{\texttt{it}} = \widetilde{\mathcal{O}}\left(  r_{\cM,\delta}\sqrt{\frac{\delta(\|\vx^{\star}-\vx^0\|^2+\epsilon_0)}{\epsilon}}\log  \frac{\delta(\|\vx^{\star}-\vx^0\|^2+\epsilon_0)}{\epsilon}\right). 
$$ The total number of communications is $N_{\texttt{com}} =  N_{\texttt{it}}\cdot N_{\texttt{com},\cM}$.
\end{theorem}

We note that our findings diverge from those of the Catalyst framework in the centralized case, as   in works like
\cite{lin2015universal, lin2018catalyst}. 
Beyond extending to the decentralized setting, a critical distinction lies in the nature of the termination time $T_k$, which is deterministic in our approach, unlike the stochastic nature of those  in the literature, which predict termination in expectation. This makes our termination more practical.   

The boundedness of the sequence is a standard assumption for (nonstrongly) convex losses, which is implied e.g. by the compactness of   $\texttt{dom}\,r$--a frequent constraint in machine learning applications that helps manage solution complexity.  
Alternatively, Theorem~\ref{theorem:total-iter-cvx} may  be reformulated under the coerciveness of  $u$ rather than boundedness of the iterates--another conditions widely used  in the machine learning community (e.g.,  \cite{lin2015universal, lin2018catalyst, li2020revisiting}). In this   case,    $T_k$  would be  defined  such that $\cL^k(\vs^{k,T_k})\leq \epsilon_{k+1}$, making it a random variable. We can then  demonstrate that   $\mathbb{E}[T_k]$ is of the same order of the RHS of   (\ref{inner-lp-cvx}). We omit further details because of  page limits.  

   \section{Applications of the DCatalyst  Framework} \label{sec:Applications} 
This section applies DCatalyst  to three representative non-accelerated decentralized algorithms: SONATA \cite{sun2020distributed}, PUDA~\cite{alghunaim2020decentralized}, and PMGT-LSVRG~\cite{ye2021pmgtvr}. This  selection is motivated by the distinctive features of each of these algorithms.  {\bf (i)} SONATA is applicable to both strongly convex and convex (possibly composed) functions $u$, and can exploit function similarity to enhance convergence rates. By integrating it with DCatalyst, we aim to inherit these features while achieving accelerated rates. {\bf (ii)}  PUDA is representative of a variety of (primal-dual) (proximal) decentralized algorithms,   which includes others such as EXTRA \cite{shi2015extra}, DIGing \cite{DIGing}, and Exact Diffusion \cite{ExactDiffusion}--employing DCatalyst provides a unified approach to accelerating all these methods.  {\bf (iii)}  For scenarios where each $f_i$ has a finite-sum structure (satisfying Assumption~\ref{assump:class-finite-sum}), we apply DCatalyst  to accelerate PMGT-LSVRG~\cite{ye2021pmgtvr},   marking the first accelerated variance-reduction decentralized algorithm applicable to composite functions   $u$.    Another by-product of our framework  is DCatalyst-VR-EXTRA, which is obtained by applying  DCatalyst   to VR-EXTRA   \cite{li2022variance};  we omit the details because of space constraints.

To   evaluate the performance of DCatalyst, in addition to the total (i.e., inner-plus-outer) communications towards $\varepsilon$-optimality,   we report  also  the total computational complexity, measured by the total  number of (proximal) gradient  updates conducted by each   agent.

\subsection{\texorpdfstring{DCatalyst SONATA \cite{sun2020distributed}}{DCatalyst SONATA} }\label{subsec:application-sonata}
The SONATA algorithm \cite{sun2020distributed} applied to Problem (\ref{eq:problem})
reads: for each $i\in [m]$ and $t\geq 0$ (and proper initialization), 
  \begin{align*} 
     x^{t+1/2}_i    & = 
      \argmin_{x\in \mathbb{R}^{d}}\, \tilde{f}_i(x;x_i^t)+\inn{y_i^t-\nabla f_i(x_i^t)}{x-x_i^t}+r(x),\\
          x_i^{t+1} &= \sum_{j=1}^m w_{ij}x_j^{t+1/2},\\
         y_i^{t+1}& = \sum_{j=1}^m w_{ij}\left(y_j^t+\nabla f_j(x_j^{t+1})-\nabla f_j(x_j^t)\right).
    \end{align*}
 Here,  $\tilde f_i$ is a local surrogate of $f_i$, chosen to exploit the potential structure of $f_i$; we will consider explicitly the following two instances of $\tilde f_i$:\begin{equation*}
\tilde{f}_i(x;x_i^t) = \left\{\begin{aligned}
    & f_i(x)+\frac{\beta}{2}\norm{x-x_i^t}^2,& (\text{SONATA-F, under Ass.~\ref{assump:class-similarity}}); \\
 & f_i(x_i^t)+\inn{\nabla f_i(x_i^t)}{x-x_i^t}+\frac{L}{2}\norm{x-x_i^t},    &(\text{SONATA-L}).
\end{aligned}\right. 
\end{equation*}
The first choice retains the entire function $f_i$, which will be shown being particularly suitable to exploit function similarity among the $f_i$ (specifically, under Assumption~\ref{assump:class-similarity}), while the second one retains only the first-order information, generally leading to prox-friendly subproblems. The gossip weights, $w_{ij}$, are chosen such that (i)  $w_{ij} > 0$ if  $(i,j) \in \mathcal{G}$, and $w_{ij} = 0$ otherwise; and (ii) $W=[w_{ij}]_{i,j=1}^m\in \mathbb{R}^{m\times m}$ is doubly stochastic, i.e., $ {1}_m^{\top}W =  {1}_m^{\top}$ and $W {1}_m =  {1}_m$. We denote   $\rho := \norm{W-(1/m){1}_m {1}_m^{\top}}<1$.  Notice that several rules have been proposed for the choice of such a $W$; example include the Laplacian, the Metropolis-Hasting, and the maximum-degree weight rules \cite{network-rule,Nedic_Olshevsky_Rabbat2018}.

 We apply the DCatalyst framework to  SONATA-F and SONATA-L, which serve as the algorithm 
$\cM$
  in the inner loop of Algorithm~\ref{alg:framework}.  In the following, we outline the specific tuning of these algorithms to ensure they satisfy all the assumptions required for the application of Theorem~\ref{theorem:total-iter-scvx}  and Theorem~\ref{theorem:total-iter-cvx}. Throughout this section, when using  SONATA-F we will implicitly assume similarity among the $f_i$'s, as stipulated by Assumption~\ref{assump:class-similarity}.

  \smallskip

\noindent  $\bullet$ {\bf On Assumption~\ref{Lyapunov}:} The result below summarizes the convergence property of  SONATA-F and SONATA-L, in agreement with  Assumption~\ref{Lyapunov}.  

\begin{lemma}\label{application:sonata-assump4}
   Consider SONATA applied to  Problem (\ref{eq:problem}) under  Assumption~\ref{assump:class}  (with $\mu>0$) with the following initialization:  
$x^0_i\in \texttt{dom}\, r$ and   $y_i^0=\nabla f_i(x_i^0)$,  for all $i\in [m]$. Suppose 
 $$
\rho^2 \leq \left\{\begin{aligned}
 & \frac{\mu^2\beta^2}{5712(L+\beta)^2\left(9\beta^2+4\left(L+\beta\right)^2\right)}= \mathcal{O}\left(\left(\frac{\beta/\mu}{\kappa}+1\right)^{4}\right), &(\text{SONATA-F});\\
    & \min\left\{\frac{\mu^2}{13440 L_{\max}^2}, \frac{L^2}{12L^2+84L_{\max}^2}\right\} = \mathcal{O}\left(\kappa_{\ell}^{-2}\right), & (\text{SONATA-L}).
\end{aligned}\right. 
 $$
Then,  Assumption~\ref{Lyapunov} holds with the following positions: \begin{itemize}
    \item Lyapunov function: 
 $$
        \cL(\vs) =  \frac{2}{\mu m}\summ \left(u(x_i)-u^\star\right) + \frac{\eta}{m}\left(4L_{\max}^2\|\vx_{\bot}\|^2+2\|\vy_{\bot}\|^2\right),
$$        
    where $\vs = [\vx,\vy],$  
$$
   {\vx_{\bot}} := \vx - 1_{m}(\bx)^{\top},\,\,  {\vy_{\bot}} := \vy - 1_{m}(\bar{y})^{\top},\,\, \text{and}\,\,  \eta = \left\{\begin{aligned}
    &  \frac{34}{\mu(16\beta+\mu)}, & (\text{SONATA-F}), \\
 & \frac{10}{\mu(4L+\mu)},    &(\text{SONATA-L}).
\end{aligned}\right. 
 $$
\item Linear convergence: the contraction property (\ref{assum3}) holds with $$c_{\mathcal L}=1\quad \text{and}\quad r_{\cM} = \left\{\begin{aligned}
 &  2+\frac{32\beta}{\mu} = \mathcal{O}\left(\frac{\beta}{\mu}\right), &(\text{SONATA-F}); \\
    & 2+\frac{8L}{\mu} = \mathcal{O}\left(\kappa_g\right), &(\text{SONATA-L}).
\end{aligned}\right.$$
\end{itemize}
   
\end{lemma}

\noindent  $\bullet$ {\bf On Assumption~\ref{assump:warm-start}:} When embedded in the inner loop of Algorithm~\ref{alg:framework} to minimize $u^k$,   and  initialized at the beginning of the $(k+1)$th inner loop  as  \begin{equation}\label{cata-sonata:ini}
    \begin{aligned} 
        & \vx^{k+1} = \vx^{k,T_k}, \\
        & \vy^{k+1} = \vy^{k,T_k}+\delta\left(\vz^{k+1}-\vz^{k}\right),
    \end{aligned}\quad \text{with}\quad \delta = \left\{\begin{aligned}
 &  \beta-\mu, &(\text{SONATA-F}); \\
    & L-\mu , &(\text{SONATA-L}),
\end{aligned}\right. 
\end{equation} 
 SONATA satisfies Assumption~\ref{assump:warm-start}, as proved next. %

\begin{lemma}\label{application:sonata-assump5}
   Consider   SONATA applied to the subproblem (\ref{eq:subproblem_DCat}), with initialization    (\ref{cata-sonata:ini})  and  tuning  as in Lemma~\ref{application:sonata-assump4} (applied here to $u^k$).  
 Assumption~\ref{assump:warm-start} holds with the following positions:    
    $$ c_{\cM}= 2\,\, \text{and}\,\, d_{\cM} =  \frac{2\delta^{2}}{(\delta+\mu)^2}+16  \eta \delta^{2}, \,\,\text{with}\,\, 
         \eta  := \left\{\begin{aligned}
    &  \frac{34}{(\mu+\delta)(16\beta+\mu+\delta)}, & (\text{SONATA-F}), \\
 & \frac{10}{(\mu+\delta)(4L+\mu+5\delta)},    &(\text{SONATA-L}).
\end{aligned}\right. $$

\end{lemma}\begin{proof}
    See Appendix~\ref{app:proof-lemma_app_sonata}.
\end{proof}
\noindent $\bullet$ \textbf{Convergence of DCatalyst-SONATA:} Since both   Assumptions~\ref{Lyapunov} and \ref{assump:warm-start} are satisfied,   convergence  of DCatalyst, equipped with  SONATA in the inner loop,  comes readily from    Theorems~\ref{theorem:total-iter-scvx} (for strongly convex $u$) and~\ref{theorem:total-iter-cvx} (for  convex $u$), as summarized below.

\begin{corollary}\label{coro:sonata-scvx} Given    Problem (\ref{eq:problem}) under Assumption~\ref{assump:class} with $\mu > 0$, let $\{\mathbf{x}^k\}$ be the sequence generated by Algorithm~\ref{alg:framework} in  the following setting:  {\bf (i)} in Step  \texttt{(S.1)}, either SONATA-L or SONATA-F is employed, with    initialization as in (\ref{cata-sonata:ini}) and  tuning  as in Lemma~\ref{application:sonata-assump4} (applied   to $u^k$); {\bf (ii)} in  Step $(\texttt{S.2})$,   $\alpha^{k}$  is chosen as   $\alpha^{k} =  \alpha := \sqrt{ \mu/(\mu+\delta)}$, for all $k\geq 0$; {\bf (iii) }     for each outer loop $k\geq 0$,  SONATA    is terminated after     $T_k$ inner iterations, given by    
\begin{equation*}
\begin{aligned}
 T_k \! = \!\left\lceil r_{\cM,\delta}\log\frac{ 2  + \frac{36}{(1-c\,\alpha)^{2}}\left( \frac{2\delta^{2}}{(\delta+\mu)^2}+16 \eta\delta^{2}\right)\left(2+\frac{\delta}{\mu}+\frac{(2m+1)(\mu+\delta)^2}{\mu^2(1-c)}+  \frac{2\sqrt{2000}(\mu+\delta)^2\sqrt{m}}{\mu^2(1-c)^2} \right)}{1-c\alpha}\right\rceil,
\end{aligned}
\end{equation*}
where $\eta$ is the same as in Lemma~\ref{application:sonata-assump5} and 
\begin{equation}\label{application:sonata-r-cM}
    \begin{aligned}
    r_{\cM,\delta} =  \left\{\begin{aligned}
 &  2+\frac{32\beta}{\mu+\delta} = 34,    &(\text{DCatalyst-SONATA-F}); \\
    & 2+\frac{8(L+\delta)}{\mu+\delta} \leq 18, &(\text{DCatalyst-SONATA-L}).
\end{aligned}\right.  
\end{aligned}
\end{equation}
Then,    
    $(1/m)\,\|\vx^{k}-\vx^{\star}\|^{2}   \leq \ep$    after a total number of  communication rounds and computations    given respectively by       
\begin{equation*}      \left\{\begin{aligned}
 &  \widetilde{\mathcal{O}}\left(\sqrt{\frac{\beta}{\mu(1-\rho)}}\log\frac{1}{\epsilon}\right)   \,\quad &\text{and}&\quad\,\,  \widetilde{\mathcal{O}}\left(\sqrt{\frac{L+\beta}{\beta}}\cdot\frac{\beta}{\mu}\log^2\frac{1}{\epsilon}\right),    \quad &(\text{DCatalyst-SONATA-F}); \\
    &    \widetilde{\mathcal{O}}\left(\sqrt{\frac{\kappa_g}{1-\rho}}\log\frac{1}{\epsilon}\right)   \,\quad &\text{and}&\quad\,\, \widetilde{\mathcal{O}}\left(\sqrt{\kappa_g}\log\frac{1}{\epsilon}\right), \quad &(\text{DCatalyst-SONATA-L}).
\end{aligned}\right.  
\end{equation*} 
  \end{corollary}
  \begin{proof}
      See Appendix~\ref{app:proof:coro:sonata-scvx}. 
  \end{proof}

\begin{corollary}\label{coro:sonata-cvx}  
Given    Problem (\ref{eq:problem}) under Assumption~\ref{assump:class} with $\mu = 0$, let $\{\mathbf{x}^k\}$ be the sequence generated by Algorithm~\ref{alg:framework} in  the following setting:  {\bf (i)} in Step  \texttt{(S.1)}, either SONATA-L or SONATA-F is employed, with    initialization as in (\ref{cata-sonata:ini}) and  tuning  as in Lemma~\ref{application:sonata-assump4} (applied   to $u^k$); {\bf (ii)} in  Step $(\texttt{S.2})$,  $\alpha^{k}$  is chosen recursively according to (\ref{eq:alpha-rec}); {\bf (iii)} the iterates are bounded, $\|x_i^k\|\leq B$ for any $i\in[m],$ $ k\geq 0$ and some $B \in (0,\infty)$; {\bf (iv) }     for each outer loop $k\geq 0$,  SONATA    is terminated after     $T_k$ inner iterations, given by
$$
\begin{aligned}
 T_k  =  \left\lceil  r_{\cM,\delta}\log \left(2^{5+2r_0}+1224(k+3)^{4+2r_0}\right)\right\rceil  
 = \widetilde{\mathcal{O}}(\log k). 
\end{aligned}
 $$ 
Then,     $   (1/(2\delta m))\cdot \summ \|\nabla \Mud(x_i^{K})\|^{2}  \leq \epsilon, $  after a total number of  communication rounds and computations respectively given by        
\begin{equation*}      \left\{\begin{aligned}
 &\widetilde{\mathcal{O}}\left(\sqrt{\frac{\beta}{(1-\rho)\epsilon}}\log\frac{1}{\epsilon} \right) \!\!\!\!\!\!\!\!\!\!\!\!  \quad &\text{and}&\quad\,\,\, \widetilde{\mathcal{O}}\left(\sqrt{\frac{L+\beta}{\epsilon}}\log^2\frac{1}{\epsilon} \right),    \quad &(\text{DCatalyst-SONATA-F}); \\
    &\widetilde{\mathcal{O}}\left( \sqrt{\frac{L}{(1-\rho)\epsilon}}\log\frac{1}{\epsilon}\right) \quad  &\text{and}& \quad\,\,\, \widetilde{\mathcal{O}}\left(\sqrt{\frac{L}{\epsilon}}\log\frac{1}{\epsilon} \right), \quad &(\text{DCatalyst-SONATA-L}).
\end{aligned}\right.  
\end{equation*}    
\end{corollary}

\begin{proof}
    See Appendix~\ref{app:proof:coro:sonata-cvx}.
\end{proof}

Corollary~\ref{coro:sonata-scvx} and~\ref{coro:sonata-cvx} demonstrate that DCatalyst-SONATA-L can achieve optimal performance in terms of both communication and computational complexities,   for composite functions $u$.  BY exploiting function similarity, DCatalyst-SONATA-F achieves (nearly) optimal communication complexity (Table~\ref{table:scvx}) $\widetilde{\mathcal{O}}((1/\sqrt{1-\rho})\cdot \sqrt{\beta/\mu}\cdot\log(1/\epsilon))$, which compare  favorably with existing algorithms, especially when $f$ is ill-conditioned.

\subsection{\texorpdfstring{Accelerating   PUDA \cite{alghunaim2020decentralized}}{Accelerating   PUDA}}
The PUDA algorithm applied to Problem~(\ref{eq:problem}) (with $\mu_{\min} > 0$) reads: for any $t\geq 0$ (with proper initialization of $\vx^0, \vy^0\in\mathbb{R}^{m\times d}$),
\begin{equation}\label{application:puda-update}
    \begin{aligned}
    &\vy^{t+1/2} =  (I-C)\vx^{t}-\eta\nabla F(\vx^t)-H\vy^t, \\
    & \vy^{t+1} = \vy^t+H\vy^{t+1/2},\\
    & \vx^{t+1} = \text{\bf prox}_{\eta r} (W\vy^{t+1/2}), 
    \end{aligned}
\end{equation}
where $\nabla F(\vx^{k}) \in \mathbb{R}^{m\times d}$ denotes $[\nabla f_1(x_1^k), \cdots, \nabla f_m(x_m^k)]^{\top}$, and the choice of $W, H, C \in\mathbb{R}^{m\times m}$ determines the specific algorithm for consideration. For instance,   EXTRA \cite{shi2015extra} and Prox-ED \cite{alghunaim2020decentralized} are special instances of PUDA--see \cite{alghunaim2020decentralized} for details. To ensure that DCatalyst is compatible with the diverse algorithms encompassed by PUDA,  we do not restrict the structure of these matrices but rather adhere to the  conditions for PUDA's convergence as outlined in \cite{alghunaim2020decentralized}, namely:  
(i) $W$ is symmetric and doubly stochastic; (ii)  $H$ and $C$ are symmetric, gossip matrices satisfying  $H\vx \,\,(\text{resp.} C\vx) = 0 \Leftrightarrow \left(I-\frac{1}{m}1_m 1_{m}^{\top}\right) {\vx} = 0$; and (iii)    $W, H,$ and $C$ satisfy
 $
    W^2 \preceq I-H^2,\quad 0\preceq C \preceq 2I.
 $

We show next that    Theorem~\ref{theorem:total-iter-scvx} and \ref{theorem:total-iter-cvx} are applicable to state convergence of DCatalyst when PUDA is used an decentralized algorithm to solve the inner   subproblems (\ref{eq:subproblem_DCat}).\smallskip

\noindent  $\bullet$ {\bf On Assumption~\ref{Lyapunov}:} In the setting above, and under a suitable tuning, PUDA applied to Problem~(\ref{eq:problem}) (under $ \mu_{\min}>0$)  satisfies  Assumption~\ref{Lyapunov}.  

\begin{lemma}\label{application-puda:assum4} 
Consider PUDA applied to Problem~(\ref{eq:problem}) under Assumption~\ref{assump:class} (with $\mu_{\min} > 0$), with   stepsize $\eta = (2-\sigma_{\max}(C))/(2 L_{\max} )$, and initialization: $x_i^0\in\texttt{dom}\,r$, and $y_i^0\in \texttt{range}(H)$,  for all $i\in[m]$. Then, Assumption~\ref{Lyapunov} holds with the following positions:
\begin{itemize}
    \item Lyapunov function: 
 
$$
            \cL(\vs) := \overm\summ\|x_i-x^\star\|^2 + \overm\summ\|y_i-y^\star\|^2,
       $$
where $\vs = [\vx, \vy]$, and $y^{\star}$ is the unique limit point of $\{y_i^{k}\}_{k\geq 0}$ (belonging to \texttt{range}(H)) 
\cite{alghunaim2020decentralized}).
    \item Linear convergence: the contraction property ($\ref{assum3}$) holds with:
    $$
  c_{\cL} = 1 \quad \text{and}\quad r_{\cM} = \max\left\{ \frac{4L_{\max}}{\left(2-\sigma_{\max}(C)\right)^2\mu_{\min}}, \frac{1}{\sigma_{\min}^{+}(H^2)}\right\}.
    $$
\end{itemize} 
\end{lemma}

\noindent  $\bullet$ {\bf On Assumption~\ref{assump:warm-start}:}    When embedding PUDA in the inner loop   (S.1) of Algorithm~\ref{alg:framework}, the selection of $\delta$ and the initialization of the $(k+1)$th inner loop are: 
\begin{equation}\label{cata-puda:ini}
        \begin{aligned}
        &\vx^{k+1} = \vx^{k,T_k},\\
        &\vy^{k+1} = \vy^{k,T_k},
    \end{aligned}\quad\text{with}\quad \delta = \frac{    4\sigma_{\min}^{+}(H^2)\left(L_{\max}-\mu_{\min}\right) }{  \left(2-\sigma_{\max}(C)\right)^2-4\sigma_{\min}^{+}(H^2) }-\mu_{\min},
\end{equation}
  then PUDA satisfies Assumption~\ref{assump:warm-start}, as shown below. \begin{lemma}\label{application-puda:assum5}
    With the initialization and selection of $\delta$ as in  (\ref{cata-puda:ini}), Assumption~\ref{assump:warm-start} holds with   
$$
          c_{\cM} = 2\quad\text{and}\quad d_{\cM} = 2+\frac{    \sigma_{\max}(H^2)\left(9+9\delta^2\eta^2    \right)}{\sigma_{\min}^{+}(H^2)},
  $$
  where $\eta$ denotes the stepsize of PUDA with the tuning as in Lemma~\ref{application-puda:assum4} (but applied to $u^k$).
\end{lemma}
\begin{proof}
    See Appendix~\ref{app:proof-puda-assum5}.
\end{proof}

\noindent $\bullet$ \textbf{Convergence of DCatalyst-PUDA:} Since both of Assumption~\ref{Lyapunov} and ~\ref{assump:warm-start} hold, we can apply Theorem~\ref{theorem:total-iter-scvx} to assess  the convergence property of DCatalyst equipped with PUDA.

\begin{corollary}\label{coro:puda-scvx} Given    Problem (\ref{eq:problem}) under Assumption~\ref{assump:class} with $\mu_{\min} > 0$, let $\{\mathbf{x}^k\}$ be the sequence generated by Algorithm~\ref{alg:framework} in  the following setting:  {\bf (i)} in Step  \texttt{(S.1)}, PUDA is employed, satisfying   (\ref{cata-puda:ini}) and  tuning  as in Lemma~\ref{application-puda:assum4} (applied   to $u^k$); {\bf (ii)} in  Step $(\texttt{S.2})$,   $\alpha^{k}$  is chosen as   $\alpha^{k} =  \alpha := \sqrt{\mu_{\min}/(\mu_{\min}+\delta)}$, for all $k\geq 0$; {\bf (iii) }     for each outer loop $k\geq 0$,  PUDA is terminated after     $T_k$ inner iterations, given by    
\begin{equation*}
\begin{aligned}
 T_k   = \!\!\left\lceil   r_{\cM,\delta}\log\frac{ 2\! + \!\frac{36}{(1-c\,\alpha)^{2}}\!\left(\!2\!+\!\frac{    \sigma_{\max}(H^2)\left(9+9\delta^2\eta^2   \right)\!}{\sigma_{\min}^{+}(H^2)}\right)\!\!\left(2\!+\!\frac{\delta}{\mu_{\min}}\!+\!\frac{(2m+1)(\mu_{\min}+\delta)^2}{\mu_{\min}^2(1-c)}\!+\!  \frac{2\sqrt{2000}(\mu_{\min}+\delta)^2\sqrt{m}}{\mu_{\min}^2(1-c)^2}  \!\right)}{1-c\alpha}\!\right\rceil,
\end{aligned}
\end{equation*}
where $\eta$ is the same as in Lemma~\ref{application-puda:assum5}, and $r_{\cM,\delta} =  1/\sigma_{\min}^{+}(H^2)$.
Then,    
    $(1/m)\,\norm{\vx^{k}-\vx^{\star}}^{2}   \leq \ep$    after a total number of  communication rounds and computations    given by   
$$\widetilde{\mathcal{O}}\left(\sqrt{\frac{\kappa_{\ell}}{\left(\left(2-\sigma_{\max}(C)\right)^2-4\sigma_{\min}^{+}(H^2)\right)\sigma_{\min}^{+}(H^2)}}\log\frac{1}{\epsilon}\right).
    $$
  \end{corollary}

  \begin{proof}
      See Appendix~\ref{app:proof:coro:puda-scvx}.
  \end{proof}

  To further comment the convergence results above,   let us consider a specific admissible choice of the weight matrices  within the PUDA framework, as used in the  Prox-ED algorithm \cite{alghunaim2020decentralized}, namely:   $W=(I+\tilde{W})/2$, $H^2=(I-\tilde{W})/2$, and $C =0$, where $\tilde{W}\in\mathbb{R}^{m\times m}$ is symmetric, doubly stochastic and primitive matrix matching the graph $\mathcal G$, with associated $\rho =\|\tilde{W}-(1/m) 1_{m}1_{m}^\top\|<1$. 
    The iteration and communication complexity of the plain Prox-ED reads
$$
 \mathcal{O}\left(\left(\kappa_{\ell}+\frac{1}{1-\rho}\right)\log\frac{1}{\epsilon}\right).
$$
By applying our acceleration framework,    DCatalyst-Prox-ED improves to
$$\widetilde{\mathcal{O}}\left(\sqrt{\frac{\kappa_{\ell}}{1-\rho}}\log\frac{1}{\epsilon}\right).
$$
This shows  a more favorable dependence on the agents' losses condition numbers $\kappa_\ell$ and network connectivity $\rho$.  Yet, the dependence on the condition number is more favorable in DCatalyst-SONATA-L (see Corollary~\ref{coro:sonata-scvx}), with the latter being the global condition number $\kappa_g$ rather than the local one $\kappa_\ell$.

\subsection{\texorpdfstring{ Accelerating    PMGT-LSVRG \cite{ye2021pmgtvr}}{ Accelerating    PMGT-LSVRG}}
We apply  the DCatalyst framework to the classes of (\ref{eq:problem})
where the agent losses have an additive separable structure, as specified in   (\ref{eq:f_ij}) (Assumption~\ref{assump:class-finite-sum}). We bring acceleration to the variance reduction decentralized algorithm  PMGT-LSVRG \cite{ye2021pmgtvr}. 

PMGT-LSVRG applied to Problem~(\ref{eq:problem}) reads: for any $t \geq 0$ (and proper initialization), 
 
\begin{equation}\label{application:pmgt-update}
\begin{aligned}
    &\vx^{t+1} = \texttt{FastMix}\left(\textbf{prox}_{\eta r}\left(\vx^{t}-\eta \vy^t\right), N_{\texttt{FM}}\right),\\
    & \text{Each agent $i\in[m]$ computes:}
    \\ 
    & \quad\quad v_i^{t+1} = \left\{\begin{aligned}
        &x_i^t,\quad \text{with probability}\,\,\,\frac{1}{n},\\
        & v_i^t, \quad\text{otherwise};
    \end{aligned}\right. \\ 
    &\quad\quad  \tg_i^{t+1} = \left\{\begin{aligned}
        &\nabla f_i(x_i^t),\,\,\, \text{if $v_i^{t+1} = x_i^t$}, \\
        & \tg_i^t, \quad\quad\quad\text{otherwise};
    \end{aligned}\right. \\  
    & \quad\quad\text{Pick $j_i\in[n]$ with probability}\,\,\,p_{ij_i} =  \frac{L_{ij_i}}{\sum_{j}^n L_{ij}},\,\,\text{and update} \\  
           &\quad\quad\,\,  g_i^{t+1} = \frac{1}{np_{ij_i}}\left(\nabla f_{ij_i}(x_i^{t+1})-\nabla f_{ij_i}(v_i^{t+1})\right)+\tg_i^{t+1},\\
        & \vy^{t+1} = \texttt{FastMix}\left(\vy^{t}+\vg^{t+1}-\vg^{t}, N_{\texttt{FM}}\right). 
 \end{aligned}  
 \end{equation}
Here,  \texttt{FastMix} is defined   in Algorithm~\ref{alg:fastmix} below, used for accelerating  the consensus steps   \cite{liu2011accelerated};   
$N_{\texttt{FM}}$
  is a parameter representing the number of communications per iteration, and $\rho = \norm{W-(1/m)
1_{m}1_{m}^{\top}}$, associated with the  gossip matrix $W\in\mathbb{R}^{m\times m}$ defined   as in Sec.~\ref{subsec:application-sonata}. In PMGT-LSVRG,  the $x$-and $y$-variables are the decision and tracking variables. The gradient tracking mechanism employs a variance reduction technique via the auxiliary variables $\tilde{g}_i$'s and $v_i$'s, whereby at each iteration the full batch gradient $\nabla f_i(x_i)$ is computed only with probability $1/n$ and $v_i$ represents the most recent iterate at which $\nabla f_i(\cdot)$ is evaluated.
  The number of communications per   iteration is counted as $N_{\texttt{com},\cM} = 2N_{\texttt{FM}}$.

\begin{algorithm}[!ht] 
\caption{\texttt{FastMix}$\left(\vx^0, N_{\texttt{FM}}\right)$}\label{alg:fastmix}
 {\bf Input}: $\vx^0 = (x_i^0)_{i\in [m]},\,\, N_{\texttt{FM}}$; $\rho$;    
{\bf Output}: $\vx^{N_{\texttt{pm}-1}}$;  
 {\bf Set}:  $\vx^{-1} := \vx^{0}$  and $\eta: = \frac{1-\sqrt{1-\rho^2}}{1+\sqrt{1-\rho^2}}$;\\
\textbf{for} {$t=0,1,2,\cdots, N_{\texttt{FM}}-1$} \textbf{do} \\ 
 $\vx^{t+1} = (1+\eta)W\vx^t-\eta\vx^{t-1}$;\\
\textbf{end for}
\end{algorithm} 

\begin{remark}
Differently from the original PMGT-LSVRG \cite{ye2021pmgtvr}, we introduced here a variant where the index $j_i$ in (\ref{application:pmgt-update}) is selected with  probability $p_{ij_i} = ({L_{ij_i}})/{\sum_{j}^n L_{ij}}$ rather than uniform. As it will be  showed in the convergence results below, this  improves the rate dependence of the algorithm the condition number of the problem, from   $\tilde{\kappa}_s$ to $\kappa_s$. \end{remark}

  \smallskip
\noindent  $\bullet$ {\bf On Assumption~\ref{Lyapunov}:} The result below summarizes the convergence   of PMGT-LSVRG, in agreement with  Assumption~\ref{Lyapunov}.  The proof builds on  \cite{ye2021pmgtvr} and is omitted.
 
\begin{lemma}\label{application:pmgt-assum4} 
Consider the PMGT-LSVRG   applied to Problem~(\ref{eq:problem}) under Assumption~\ref{assump:class} ($\mu > 0$) and Assumption~\ref{assump:class-finite-sum}, with the following initialization: for all $i\in[m]$, $x_i^0 = v_i^0 \in\texttt{dom}\, r$;    $g_i^0$  and  $\tilde{g}_i^0$ are   unbiased estimators of $\nabla f_i^0(x_i^0)$ and $\nabla f_i^0(v_i^0)$, respectively; and $(1/m)\summ y_i^0$ is an unbiased estimator of $(1/m)\summ \nabla f_i(x_i^0)$.      
Further,   $\eta = 1/(16L_{\max})$ and  $N_{\texttt{FM}} = (1/\sqrt{1-\rho})\cdot\log\left(36\max\{6\kappa_s, n\}\right)$.   
Then, Assumption~\ref{Lyapunov} holds with following positions:
\begin{itemize}
    \item Lyapunov function:
    $$
    \cL(\vs) := \frac{1}{c_{\texttt{pm}}}\left(\norm{\bar{x}-x^{\star}}^2+4n\eta^2\Delta_{f}\right)+\overm\left(\norm{\vx_{\bot}}^2+\eta^2\norm{\vy_{\bot}}^2\right),
    $$
where $\vs := [\vx, \vy, \vg, \vv, \tilde{\vg}]$, $c_{\texttt{pm}} = 20r_{\texttt{pm}}/(1-40r_{\texttt{pm}}\rho_{\texttt{pm}}^2)$, $r_{\texttt{pm}} = \max\{12\kappa_s, 2n\}$,
$\rho_{\texttt{pm}} = (1-\sqrt{1-\rho})^{N_{\texttt{FM}}}$, 
$$
\Delta_{f} \!:=\! \frac{1}{mn}\summ\sum_{j=1}^{n}\frac{1}{np_{ij}}\|\nabla f_{ij}(v_i)-\nabla f_{ij}(x^\star)\|^2,\,\, \vx_{\bot} \!:=\! \vx-1_m(\bar{x})^{\top},\,\, \vy_{\bot} \!:=\! \vy-1_m(\bar{y})^{\top}.
$$
    \item Linear convergence: the contraction (\ref{assum3}) holds with
    $
    c_{\cL} = 1$  and  $r_{\cM} = 4r_{\texttt{pm}}$.

\end{itemize}
 \end{lemma}

\noindent  $\bullet$ {\bf On Assumption~\ref{assump:warm-start}:} Suppose PMGT-LSVRG is embedded in the inner loop of Algorithm~\ref{alg:framework}, with the following choice of $\delta$ and   warm-start at   at the beginning of the $(k+1)$th inner loop:  \begin{equation}\label{cata-pmgt:ini} 
    \begin{aligned}
      &    \vx^{k+1} = \vx^{k,T_k},
        \,\,\vv^{k+1} = \vv^{k,T_k},   
         \vy^{k+1} = \vy^{k,T_k}+\delta(\vz^{k}-\vz^{k+1}),  \,\,
         \vg^{k+1} = \vg^{k,T_k}+\delta(\vz^{k}-\vz^{k+1}), \\&    \tilde{\vg}^{k+1} = \tilde{\vg}^{k,T_k}+\delta(\vz^{k}-\vz^{k+1}), \quad \text{and}\quad \delta = \frac{\bar{L}_{\max}}{n}-\mu. 
    \end{aligned} 
\end{equation}

Then, PMGT-LSVRG satisfies Assumption~\ref{assump:warm-start}, as proved below.

\begin{lemma}\label{application:pmgt-assum5}
Consider PMGT-LSVRG applied to subproblem~(\ref{eq:subproblem_DCat}), with the initialization (\ref{cata-pmgt:ini}) and tunning as in Lemma~\ref{application:pmgt-assum4} (applied to $u^k$). Assumption~\ref{assump:warm-start} holds with
$$
c_{\cM} = 2,\quad \text{and} \quad d_{\cM} = 2+8\eta^2\delta^2+\frac{8\eta^2(\bar{L}_{\max}+\delta)^2}{n^2}.
$$
\end{lemma}
\begin{proof}
    See Appendix~\ref{proof-pmgt-assum5}.
\end{proof}

\noindent $\bullet$ \textbf{Convergence of DCatalyst-PMGT-LSVRG:} Since both  Assumptions~\ref{Lyapunov} and ~\ref{assump:warm-start} hold, we can apply Theorem~\ref{theorem:total-iter-scvx} to assess  the convergence   of DCatalyst-PMGT-LSVRG.

\begin{corollary}\label{coro:pmgt-total} Given Problem (\ref{eq:problem}) under Assumption~\ref{assump:class} with $\mu  > 0$, let $\{\mathbf{x}^k\}$ be the sequence generated by Algorithm~\ref{alg:framework} in  the following setting:  {\bf (i)} in Step  \texttt{(S.1)}, PMGT-LSVRG is employed, satisfying    (\ref{cata-pmgt:ini})  and  tuning  as in Lemma~\ref{application:pmgt-assum4} (applied   to $u^k$); {\bf (ii)} in  Step $(\texttt{S.2})$,   $\alpha^{k}$  is chosen as   $\alpha^{k} =  \alpha := \sqrt{ \mu/(\mu+\delta)}$, for all $k\geq 0$; {\bf (iii) }     for each outer loop $k\geq 0$,  PMGT-LSVRG is terminated after     $T_k$ inner iterations, given by    
\begin{equation*}
\begin{aligned}
 T_k   =  \left\lceil r_{\cM,\delta}\log\frac{ 2  + \frac{90}{(1-c\,\alpha)^{2}}\left(2+\frac{\delta}{\mu}+\frac{(2m+1)(\mu+\delta)^2}{\mu^2(1-c)}+  \frac{2\sqrt{2000}(\mu+\delta)^2\sqrt{m}}{\mu^2(1-c)^2}  \right)}{1-c\alpha}\right\rceil, 
\end{aligned}
\end{equation*}
with  $r_{\cM,\delta} =  48n$.
Then,    
    $(1/m)\,\Eb[\|\vx^{k}-\vx^{\star}\|^{2}]   \leq \ep$    after a total number of  communication rounds and a number of computations  {(in expectation)}  given   respectively by       
    $$ \widetilde{\mathcal{O}}\left(\sqrt{\frac{\kappa_s n}{1-\rho}}\log\frac{1}{\epsilon}\right),\quad\text{and}\quad \widetilde{\mathcal{O}}\left(\sqrt{\kappa_s n}\log\frac{1}{\epsilon}\right).
    $$
  \end{corollary}

  \begin{proof}
      See Appendix~\ref{proof-pmgt-corollary}. 
  \end{proof}

Compared with    PMGT-LSVRG, (reported in Table~\ref{table:finite-sum}),  DCatalyst-PMGT-LSVRG exhibits better performance:  PMGT-LSVRG has communication  complexity $\mathcal{O}\left((1/(\sqrt{1-\rho}))\right.$ $\left.  (\tilde{\kappa}_s\log\tilde{\kappa}_s+n\log n )\cdot\log(1/\epsilon)\right)$   and computational complexity $\mathcal{O}\left((n+\tilde{\kappa}_s)\right.$ $\left.\log(1/\epsilon)\right)$, while DCatalyst-PMGT-LSVRG's reads $ \widetilde{\mathcal{O}}\left((1/\sqrt{1-\rho})\sqrt{\kappa_s n }\right.$ $\left. \log(1/\epsilon)\right)$ and $
\widetilde{\mathcal{O}}\left(\sqrt{\kappa_s n }\cdot \log(1/{\epsilon})\right)  
$, showing significant savings in communication/computation especially for ill-conditioned problems. This sets a new benchmark for composite optimization functions.  On the other hand, existing decentralized VR algorithms, such as Acc-VR-EXTRA-CA \cite{li2022variance}, ADFS \cite{hendrikx2021anOptimal} and DVR-Catalyst \cite{hendrikx2020dual}, have more favorable dependence on the condition number of $f$, namely: $\kappa_{\ell}$-dependence versus   $\kappa_{s}$ in DCatalyst-PMGT-LSVRG.  However, all these methods are applicable only to    {\it smooth} functions $u$.  DCatalyst-PMGT-LSVRG is the  first decentralized algorithm merging   acceleration and VR technique that is applicable to  composite functions $u$ (with $f_i$   additively separable).   

\section{Inexact Estimating Sequences 
}\label{sec:prelim-ES}

This section introduces  the {\it inexact estimating sequences}, a   novel  extension of the well-known estimating sequences \cite{nesterov2013introductory} to the decentralized setting. This framework is instrumental for designing and  analyzing    accelerated decentralized algorithms. The subsequent section  utilizes this machinery  to establish   convergence of Algorithm~\ref{alg:framework}.

Given Problem~(\ref{eq:problem}) (under Assumption~\ref{assump:class}) and the  the Moreau envelope  $\Mud:\mathbb{R}^{d}\rightarrow\mathbb{R}$, we recall that  $\Mud$  is $\delta$-smooth and $\muM$-strongly convex, with   
$\muM =  ({\delta\mu})/({\delta+\mu})$. Furthermore,  notice that  
$\Muds := \Mud(x^{\star}) = u^\star$, where $x^\star$ is a minimizer of $u$.  
 We are ready to state the definition of   general inexact estimating sequences. 
\begin{definition}[Inexact estimating sequence]\label{def:inexact-ES}
An {\it inexact estimating sequence} is a tuple of a sequence of  functions $\{(\psi_i^{k})_{i\in[m]}\}_{k=0}^{\infty}$,  a sequence of positive numbers $\{\alpha^k\}_{k=0}^{\infty}$, and a real-value error-sequence  $\{(\epsilon_i^{k})_{i\in[m]}\}_{k=0}^{\infty}$  such that   

{\bf (i)}  $\alpha^k\in(0,1)$, for all   $k\geq 0$, and   $\lim_{k\to \infty}\prod_{t=0}^{k}(1-\alpha^t) = 0$;   

{\bf (ii)} each $\psi_i^k:\mathbb{R}^d\rightarrow\mathbb{R}$ satisfies    
$$\psi_i^{k+1}(x^\star) \leq (1-\alpha^k)\psi_i^{k}(x^\star)+\alpha^k\left(\Muds+\epsilon_i^{k}\right),\quad \forall k\geq 0.$$
\end{definition}

The following result shows how inexact estimating sequences can be leveraged to guide the development  of appropriate  decentralized algorithms solving effectively  Problem~(\ref{eq:problem}).

\begin{lemma}\label{lemma:prelim} 
Let $\{(\psi_i^{k})_{i\in[m]}\}_{k=0}^{\infty}$,  $\{\alpha^k\}_{k=0}^{\infty}$,   $\{(\epsilon_i^{k})_{i\in[m]}\}_{k=0}^{\infty}$ be an inexact estimating sequence. Suppose there exists a sequence of iterates $\{(\tilde{x}^k_i)_{i\in[m]} \}_{k=0}^{\infty}$ such that, for all  $k\geq0$ and  $i\in[m]$,
\begin{equation}\label{eq:lower-bound-inexact} 
    \Mud(\tilde{x}_i^k) \leq \left[\psi_i^{k,\star}:=\min_{x\in\mathbb{R}^d} \psi_i^k(x)\right]  +\epsilon_{\psi,i}^k,   
\end{equation}
for some  $\{(\epsilon_{\psi,i}^{k})_{i\in[m]}\}_{k=0}^{\infty}$, with   $\epsilon_{\psi,i}^0  \equiv 0$. Then,  
\begin{equation}\label{eq:telescope-error-M}
    0 \leq \psi_i^k(x^\star) + \epsilon_{\psi,i}^k-\Muds \leq \lambda^k\left(\psi_i^0(x^\star)-\Muds+\sum_{t=0}^{k-1}\frac{\epsilon_{tot,i}^t}{\lambda^{t+1}}\right),
\end{equation}
 where \begin{equation}\label{def_lambdak}\lambda^{k}:= \prod_{t=0}^{k-1}(1-\alpha^t)\quad \text{and}\quad\epsilon_{tot,i}^k:= \ep_{\psi,i}^{k+1}-(1-\alpha^k)\ep_{\psi,i}^k+\alpha^k\ep_i^k.
\end{equation} 
\end{lemma}

\begin{proof}
    See Appendix~\ref{app:proof_lemma_preliminary}.
\end{proof} 
When  the algorithm generating  $\{(\tilde{x}^k_i)_{i\in[m]}\}_{k=0}^{\infty}$  is stochastic,     $\{(\tilde{x}^k_i)_{i\in[m]} \}_{k=0}^{\infty}$ is a random sequence, and Lemma~\ref{lemma:prelim} is understood to hold for any realization of the random variables.  

The subsequent proposition   elucidates the convergence rate of decentralized algorithms producing iterates    $\{(\tilde{x}^k_i)_{i\in[m]} \}_{k=0}^{\infty}$  conforming to  Lemma~\ref{lemma:prelim}. This convergence depends on the decay rate of the error sequences   $\mathbb{E}[|\ep_{tot,i}^{k}|]$.

\begin{proposition}\label{prop:preliminary} Let $\{(\tilde{x}^k_i)_{i\in[m]} \}_{k=0}^{\infty}$ be a sequence satisfying  Lemma~\ref{lemma:prelim}.  
\begin{itemize}
    \item[\bf (i)] $\mu > 0$:  Suppose    $(1/m)\summ\Eb\big[|\ep_{tot,i}^{k}|\big] \leq C_{scvx}\left(1-c\alpha\right)^{k}$, for all   $k\geq 0$ and  for some  problem-dependent constant $C_{scvx}$. Then,  
  $$\overm\Eb[\|\tilde{\vx}^{k}-\vx^{\star}\|^{2}]\leq c_{scvx}(1-c\alpha)^{k}, 
 $$
where $$c_{scvx} = \overm\summ\left(\frac{2}{\muM}\left(\psi_i^{0}(x^{\star})-\Muds+\frac{C_{scvx}}{\alpha(1-c)}\right)\right).$$
    \item[\bf (ii)]  $\mu = 0$:  Suppose   $(1/m)\summ\Eb\big[|\ep_{tot,i}^{k}|\big] \leq  C_{cvx}/(k+1)^{3+r_{0}}$, for all $k\geq 0$ and some   problem-dependent constant   $C_{cvx}$. Then,
\begin{equation}\label{eq:estimate-seq-sublinear-M}
\frac{1}{2\delta m}\summ\Eb[\|\nabla\Mud(\tx_i^{k})\|^{2} ]\leq  \frac{c_{cvx}}{(k+2)^{2}},  
\end{equation}
where  
$$c_{cvx} = \overm\summ\left(\psi_i^{0}(x^{\star})-\Muds+\frac{10C_{cvx}}{r_{0}} \right).
 $$

Moreover, if $r\equiv 0$, then 
\begin{equation}\label{eq:estimate-seq-sublinear-u}
    \frac{1}{2\delta m}\summ\Eb[\|\nabla u(\tx_i^{k})\|^{2}] \leq \frac{c_{cvx}(L+\delta)^2}{\delta^2(k+2)^{2}}.
\end{equation}
\end{itemize} 
\end{proposition}
\begin{proof}
    See Appendix \ref{app:proof_Cor_preliminary}.
\end{proof}  

\subsection{A constructive approach to a family of inexact estimating sequences   }\label{subsec:family-inexact-ES}
 We  provide a constructive approach to derive an explicit   class of inexact estimating sequences. which  leads to a wide range of decentralized designs. Our method unfolds in two steps:   {\bf (i)} we first specify a family of inexact  estimating sequence, parametrized by certain free quantities  (see Lemma~\ref{lm:est-seq-family} below);   {\bf (ii)} we then exploit the degrees of freedom offered by   this family  to satisfy the remaining condition (\ref{eq:lower-bound-inexact}).

\begin{lemma} \label{lm:est-seq-family}
   Assume,  for any $i\in [m]$, and  $k\geq0$,
    \begin{itemize}
        \item[{\bf (i)}]
    $\psi_i^0(\bullet):\mathbb{R}^d\rightarrow\mathbb{R}$ is an arbitrary  convex function on $\mathbb{R}^d$;
        \item[{\bf (ii)}]
            $\{(\tz_i^k)_{i\in[m]}\}_{k=0}^{\infty}$ and  $\{(e_i^k)_{i\in[m]}\}_{k=0}^{\infty}$ are arbitrary sequences in $\mathbb{R}^d$; 
        \item[{\bf (iii)}]
         $\{\alpha^k\}_{k=0}^{\infty}$ is chosen according to Definition~\ref{def:inexact-ES}(ii).  
         \item[{\bf (iv)}]   $\{(\psi_i^k)_{i\in[m]}\}_{k=0}^{\infty}$ is defined recursively  as   \begin{equation}\label{Moreau_local:ES}
 \hspace{-0.4cm}   \psi_i^{k+1}(x) := (1-\alpha^k)\psi_i^k(x)+\alpha^k \left(\Mud(\tilde{z}_i^k)+\langle\nabla\Mud(\tz_i^k)+e_i^k,x-\tz_i^k\rangle+\frac{\muM}{2}\|x-\tz_i^k\|^2\right). 
    \end{equation}
    \end{itemize}
    Then, the tuple $\{(\psi_i^k)_{i\in[m]}\}_{k=0}^{\infty}$, $\{\alpha^k\}_{k=0}^{\infty}$, $\{(\epsilon_i^k)_{i\in[m]}\}_{k=0}^{\infty}$, with $\epsilon_i^k := \inn{e_i^k}{x^{\star}-\tz_i^k}$   
    forms an inexact estimating sequence.
\end{lemma}
 
Given  the above family of inexact estimating sequences, we have two control sequences to choose,  $\{(\tz_i^k)_{i\in[m]}\}_{k}$ and   $\{(\tx_i^k)_{i\in[m]}\}_{k}$, to maintain recursively the relation (\ref{eq:lower-bound-inexact}), for a suitable sequence $\{\epsilon_{\psi,i}^k\}_k$. 
According to Lemma~\ref{lm:est-seq-family}, we are also free in the choice of $\psi_i^0$. Following \cite{nesterov2013introductory}, we choose it as a simple quadratic function:
\begin{equation}\label{def:psi_0}
    \psi_i^0(x) = \psi_i^{0,*}+\frac{\zeta^0}{2}\|x-v_i^0\|^2,
\end{equation}
for some $\psi_i^{0,*} \in\mathbb{R}$, $\zeta^0\in(0,\infty)\in \mathbb{R}^d$, and  $v_i^0\in\mathbb{R}^d$. We can obtain $\psi_i^{k,*}$  to be used in (\ref{eq:lower-bound-inexact}) in a closed form recurrence, as derived next. The proof is standard, hence it  is omitted.

\begin{lemma}[Canonical form]\label{lemma:canoi} 
Let $ \psi_i^0(x)$ be given as in (\ref{def:psi_0}). The process (\ref{Moreau_local:ES}) preserves the canonical form of functions   $\left\{(\psi_i^{k})_{i\in [m]}\right\}_{ k}$: 
 $$
        \psi_i^{k}(x) = \psi_i^{k,\star}+\frac{\zeta^{k}}{2}\|x-v_i^{k}\|^{2},\quad \forall i\in[m], \,\,k\geq0,
 $$
where the sequences $\{\zeta^{k}\}$, $\{v_i^{k}\}$, and $\{\psi_i^{k,\star}\}$ are defined as follows: 
\begin{equation}\label{recur:zeta}
\begin{aligned}
     \zeta^{k+1}  = & \,\,(1-\alpha^{k})\zeta^{k}+\alpha^{k}\muM, \\
     v_i^{k+1} = &\,\, \frac{(1-\alpha^{k})\zeta^{k}}{\zeta^{k+1}}v_i^{k}+\frac{\alpha^{k}\muM}{\zeta^{k+1}}\tz_i^{k}-\frac{\alpha^{k}}{\zeta^{k+1}}\left(\nabla \Mud(\tilde{z}_i^k)+e_i^k\right),\\
      \psi_i^{k+1,\star}   = & \,\, (1-\alpha^{k})\psi_i^{k,\star}+\alpha^{k}\Mud(\tz_i^{k})-\frac{(\alpha^{k})^{2}}{2\zeta^{k+1}}\|\nabla\Mud(\tz_i^k)+e_i^{k}\|^{2}\\
        & +\frac{\alpha^{k}(1-\alpha^{k})\zeta^{k}}{\zeta^{k+1}}\left(\inn{\nabla \Mud(\tz_i^k)+e_i^k}{v_i^k-\tz_i^k}+\frac{\muM}{2}\|\tz_i^k-v_i^k\|^{2}\right).  
\end{aligned}
\end{equation}
\end{lemma}

We are left to  to enforce (\ref{eq:lower-bound-inexact}), which holds under the following positions.

\begin{lemma}\label{lemma-lower-bound}[On condition (\ref{eq:lower-bound-inexact})] In the setting above,   let for each $i\in [m]$, \begin{itemize}
        \item[{\bf (i)}] $\psi_i^{0,\star} = \Mud(\tx_i^0)$ and $v_i^0 = \tx_i^0$;
       \item[{\bf (ii)}]   $\tz_i^0 = \tx_i^0$, and  
\begin{equation}\label{eq:specify-tz}
\tz_i^{k+1} := \tx_i^{t+1}+\frac{\alpha^k(1-\alpha^k)}{(\alpha^k)^2+\alpha^{k+1}}(\tx_i^{k+1}-\tx_i^k); 
\end{equation}
\item[\bf (iii)] and 
\begin{equation}\label{eq:specify-ei}
e_i^k := \delta(\tz_i^k-\tx_i^{k+1})-\nabla\Mud(\tz_i^k).
\end{equation} \end{itemize}
Then,  (\ref{eq:lower-bound-inexact}) holds, with
\begin{equation}\label{def:ep-psi}
    \ep_{\psi,i}^{k+1} := (1-\alpha^{k})\ep_{\psi,i}^{k}+(1-\alpha^{k})\inn{e_i^{k}}{\tx_i^{k}-\tz_i^{k}}+\inn{e_i^k+\nabla\Mud(\tz_i^k)}{\frac{1}{\delta}e_i^{k}},\,\text{and}\,\,\ep_{\psi,i}^{0} = 0.
\end{equation}  
\end{lemma} 
\begin{proof}
   The lemma is proved by induction.
Notice that  (\ref{eq:lower-bound-inexact}) holds for $k=0$, and  $\ep_{\psi,i}^{0} = 0$.  
  Suppose   (\ref{eq:lower-bound-inexact})
    hold for a given  $k>0$. 
    Then, in view of  (\ref{recur:zeta})  and $\Mud(\tx_i^{k})\geq \langle\nabla\Mud(\tz_i^k), \tx_i^k-\tz_i^k\rangle+\Mud(\tz_i^k),$
    we have
    \begin{equation*}\label{ineq:psi*}
    \begin{aligned}
        \psi_i^{k+1,\star}&\geq -(1-\alpha^{k})\ep_{\psi,i}^{k}+
        \Mud(\tz_i^{k})+  
        (1-\alpha^{k})\inn{\nabla\Mud(\tz_i^k)}{\tx_i^k-\tz_i^k}\\
        &-\frac{(\alpha^{k})^{2}}{2\zeta^{k+1}}\|\nabla\Mud(\tz_i^k)+e_i^k\|^{2}+(1-\alpha^{k})\inn{\nabla\Mud(\tz_i^k)+e_i^k}{\frac{\alpha^{k}\zeta^{k}}{\zeta^{k+1}}\left(v_i^k-\tz_i^k\right)}\\
        & \overset{(a)}{=} -(1-\alpha^{k})\ep_{\psi,i}^{k}+
        \Mud(\tz_i^{k})+  
        (1-\alpha^{k})\inn{-e_i^{k}}{\tx_i^k-\tz_i^k}\\
        &-\frac{1}{2\delta}\|\nabla\Mud(\tz_i^k)+e_i^k\|^{2}+(1-\alpha^{k})\inn{\nabla\Mud(\tz_i^k)+e_i^k}{\frac{\alpha^{k}\zeta^{k}}{\zeta^{k+1}}\left(v_i^k-\tz_i^k\right)+\tx_i^k-\tz_i^k},
    \end{aligned}
\end{equation*} 
where in (a) we used    $\zeta^{k+1} = \delta(\alpha^{k})^{2}$ for all $k\geq 0$, which is a consequence of the recursion on $\zeta^{k+1}$ in 
  (\ref{recur:zeta}), given $\zeta^{0} = ({\delta(\alpha^{0})^{2}-\alpha^{0}\muM})/({1-\alpha^{0}})$.

 Chaining the above inequality with the following (due to  (\ref{eq:specify-ei}) and  $\delta$-smoothness of $\Mud$),
$$
    \Mud(\tx_i^{k+1})\leq \Mud(z_i^k)-\frac{1}{2\delta}\|\nabla\Mud(\tz_i^k)+e_i^k\|^{2}+\frac{1}{\delta}\inn{e_i^k}{\nabla\Mud(\tz_i^k)+e_i^k},
$$
   and using (\ref{eq:specify-tz}), 
  yields  (\ref{eq:lower-bound-inexact}), with    $\ep_{\psi,i}^{k+1}$ as  
in~(\ref{def:ep-psi}). \end{proof}

\section{Convergence Analysis of DCatalyst}\label{section:analysis} 
 We are now ready to prove convergence of Algorithm~\ref{alg:framework}, leveraging  the   framework of inexact estimating sequences. Specifically, Sec.~\ref{sec:outer-loop-convergence} establishes convergence of the outer-loop,   showing that   
    $\{(z_i^k)_{i\in [m]}\}_{k=0}^{\infty}$ and $\{(x_i^k)_{i\in [m]}\}_{k=0}^{\infty}$     are an instance of the sequences  $\{(\tz_i^k)_{i\in [m]}\}_{k=0}^{\infty}$ 
 and $\{(\tx_i^k)_{i\in [m]}\}_{k=0}^{\infty}$  
 The proof of Propositions~\ref{prop:outer-loop-scvx} and~\ref{prop:outer-loop-cvx} will then follow by checking the conditions on the decay of the error sequence $(1/m)\summ\Eb\big[|\ep_{tot,i}^{k}|\big] $  as in  (i) and (ii) of Proposition~\ref{prop:preliminary}, respectively.    Subsequently, Sec.~\ref{Subsec:proof-inner-loop} is dedicated to the convergence analysis of the inner loop,  proving     Propositions~\ref{prop:inner-loop-scvx} and \ref{prop:inner-loop-cvx}.   
 Finally, the proofs of  Theorems~\ref{theorem:total-iter-scvx} and~\ref{theorem:total-iter-cvx} follow readily combining the  results from both the outer- and inner-loop analyses.

\subsection{Analysis of the outer-loop of Algorithm~\ref{alg:framework}: Proof of Propositions~\ref{prop:outer-loop-scvx} and ~\ref{prop:outer-loop-cvx}} \label{sec:outer-loop-convergence}

Setting $\tz_i^k=z_i^k$ and  $\tx_i^k=x_i^k$,  for all $i\in[m]$ and $k\geq 0$, with $\{(z_i^k)_{i\in [m]}\}_{k=0}^{\infty}$ and  $\{(x_i^k)_{i\in [m]}\}_{k=0}^{\infty}$ being the sequences generated by   Algorithm~\ref{alg:framework}, Algorithm~\ref{alg:framework} is an instance of the  the setting of Lemma~\ref{lemma-lower-bound}.  Hence, one can apply  Proposition~\ref{prop:preliminary}. Next,   we study   $(1/m)\summ\Eb[|\epsilon_{tot,i}^k|]$.   

\subsubsection{\texorpdfstring{Bounding $\epsilon_{tot,i}^k$}{Bounding epsilon\_tot,i\^k}}

Using (\ref{def_lambdak}) and  (\ref{def:ep-psi}),  
\begin{equation*}
    \begin{aligned}
        \ep_{tot,i}^{k}& =\ep_{\psi,i}^{k+1}-(1-\alpha^{k})\ep_{\psi,i}^{k}+\alpha^{k}\ep_{i}^{k }\\ &=\alpha^{k}\inn{e_i^{k}}{x^{\star}-z_i^{k}}+(1-\alpha^{k})\inn{e_i^{k}}{x_i^k-z_i^{k}}+\frac{1}{\delta}\inn{e_i^{k}}{\nabla\Mud(z_i^k)+e_i^k}\\
        &=\inn{e_i^{k}}{\alpha^{k} x^{\star}+(1-\alpha^{k})x_i^{k}-x_i^{k+1}} =\alpha^{k}\inn{e_i^{k}}{x^{\star}-v_i^{k+1}}.
    \end{aligned}
\end{equation*}
Therefore,  
\begin{equation}\label{ineq:ep-tot}
\begin{aligned}
\overm\summ\Eb[|\ep_{tot,i}^{k}|] &\leq \alpha^{k}\overm\summ\Eb[\|e_i^{k}\|\|x^{\star}-v_i^{k+1}\|]\\& \leq \alpha^{k}\left(\overm\summ\Eb[\|e_i^k\|^2]\right)^{1/2}\left(\overm\summ\Eb[\|x^{\star}-v_i^{k+1}\|^2]\right)^{1/2}.
\end{aligned}
\end{equation}

We proceed bounding the two terms above. \smallskip

\noindent   {\bf 1) Bounding $\left(\overm\summ\Eb [\|e_i^k\|^2]\right)^{1/2}$:} For all $i\in[m]$ and $ k \geq 0$, using    (\ref{eq:specify-ei}) and $x^{k,\star} = \bz^k-(1/\delta)\nabla\Mud(\bz^k)$,  
$$
e_i^k = \nabla\Mud(\bz^{k})-\nabla\Mud(z_i^{k})+\delta(\bz^{k}-z_i^{k})+\delta(\bx^{k+1}-x_i^{k+1})+\delta(x^{k,\star}-\bx^{k+1}).
$$
  Then, for any $k\geq0$, setting  $x_i^{-1} = 0$ (as a dummy variable), we have
\begin{equation}\label{local_grad_err}
\begin{aligned}
       &\|e_i^{k}\|^{2} = \|\nabla\Mud(\bz^{k})-\nabla\Mud(z_i^{k})+\delta(\bz^{k}-z_i^{k})+\delta(\bx^{k+1}-x_i^{k+1})+\delta(x^{k,\star}-\bx^{k+1})\|^{2} \\
       & \leq 12\delta^{2}\|z_i^{k}-\bz^{k}\|^{2}+3\delta^{2}\|\bx^{k+1}-x_i^{k+1}\|^{2}+3\delta^{2}\|x^{k,\star}-\bx^{k+1}\|^{2}\\
       &\leq 72\delta^{2}\|x_i^{k}-\bx^{k}\|^{2}+48\delta^{2}\|x_i^{k-1}-\bx^{k-1}\|^{2}
+3\delta^{2}\|\bx^{k+1}-x_i^{k+1}\|^{2}+3\delta^{2}\|x^{k,\star}-\bx^{k+1}\|^{2}.
\end{aligned}    
\end{equation}

We proceed distinguishing  the two cases of   $\mu > 0$ and  $\mu =0$.

\begin{itemize}
    \item[\bf (i)]  $\mu>0$:  
By Assumption~\ref{Lyapunov} and  (\ref{ineq:prop-out-scvx-Tk}), 
 $$
\begin{aligned} \overm\Eb[\|\vx^{k+1}-\vx^{k,\star}\|^{2}] \leq  \Eb[\cL^k(\vs^{k,T_{k}})]  \leq \ep_{0}(1-c\alpha)^{k+1},\quad\forall k \geq 0.
\end{aligned}
 $$
When $k \geq 0$, using (\ref{local_grad_err}), yields 
 $$
\begin{aligned}
 \Eb[\|e_i^{k}\|^{2}]
       &\leq  
        72\delta^2(1-c\alpha)^{k} m\epsilon_0 +  48\delta^2(1-c\alpha)^{k-1} m\epsilon_0 +  3\delta^2(1-c\alpha)^{k+1} m\epsilon_0.
\end{aligned}    
 $$
 Then
\begin{equation}\label{eq:final_bound_ei_scvx}\Eb[\|e_i^{k}\|^{2}] \leq  125\delta^2(1-c\alpha)^{k-1}m\epsilon_0,\quad \text{for any}\,\,  i\in[m],\,\, k \geq 0.
\end{equation}

\item[\bf (ii)]   $\mu = 0$, Invoking  Assumption~\ref{Lyapunov} and (\ref{ineq:prop-out-cvx-Tk}), and following similar steps as in the case   $\mu > 0$, we obtain:   for any $i\in[m]$ and $k\geq 0$,
\begin{equation}\label{eq:final_bound_ei_cvx}\Eb[\|e_i^{k}\|^{2}] \leq 125\delta^2m\epsilon_0\left(\frac{1}{k+1}\right)^{4+2r_0}.
 \end{equation}
\end{itemize}

\noindent  {\bf 2) Bounding $\alpha^{k}\left(\overm\summ\Eb[\|x^{\star}-v_i^{k+1}\|^2]\right)^{1/2}$}:
By the canonical form of $\psi_i^k$ and Lemma~\ref{lemma-lower-bound},   
$$
\begin{aligned}
    &\frac{\zeta^{k+1}}{2}\|x^{\star}-v_i^{k+1}\|^{2}+\Mud(x_i^{k+1})-\Muds  
    & \leq 
    \frac{\zeta^{k+1}}{2}\|x^{\star}-v_i^{k+1}\|^{2}+\psi_i^{k+1,\star}+\epsilon_{\psi,i}^{k+1}-\Muds \\
    & = \psi_i^{k+1}(x^{\star})+\ep_{\psi,i}^{k+1}-\Muds.
\end{aligned}
$$
 
By Lemma~\ref{lemma:prelim} 
, we have   $$
\begin{aligned}
   &\frac{1}{\lambda^{k+1}}\left(\frac{\zeta^{k+1}}{2}\|x^{\star}-v_i^{k+1}\|^{2}+\Mud(x_i^{k+1})-\Muds\right)   
   \leq \psi_i^{0}(x^{\star})-\Muds+\sum_{t=0}^{k}\frac{ \alpha^{t}\|e_i^{t}\|\|x^{\star}-v_i^{t+1}\|}{\lambda^{t+1}} \\
    & \leq \frac{\zeta^{0}+\delta}{2}\|x^\star-x_i^0\|^2 +\sum_{t=0}^{k}\frac{ \sqrt{\zeta^{t+1}}\|e_i^{t}\|\|x^{\star}-v_i^{t+1}\|}{\sqrt{\delta}\lambda^{t+1}},
    \end{aligned}
 $$
where in the inequality  we   used  $\zeta^{t+1} = \delta(\alpha^{t})^2$. 
Then,
\begin{equation}\label{ineq:zeta-lambda}
    \begin{aligned}
&\Eb\left[\frac{\zeta^{k+1}}{2\lambda^{k+1}}\|x^{\star}-v_i^{k+1}\|^{2}\right]\\
&\leq \frac{\zeta^{0}+\delta}{2}\|x^\star-x_i^0\|^2+\sum_{t=0}^{k}\Eb\left[\left(\sqrt{\frac{2}{\delta\lambda^{t+1}}}\|e_i^{t}\|\right)\left(\sqrt{\frac{\zeta_{t+1}}{2\lambda^{t+1}}}\|x^{\star}-v_i^{t+1}\|\right)\right] \\
& \leq \frac{\zeta^{0}+\delta}{2}\norm{x^\star-x_i^0}^2+\sum_{t=0}^{k}\left(\Eb\left[\frac{2}{\delta\lambda^{t+1}}\|e_i^{t}\|^2\right]\right)^{1/2}\left(\Eb\left[ \frac{\zeta_{t+1}}{2\lambda^{t+1}}\|x^{\star}-v_i^{t+1}\|^2\right]\right)^{1/2}.
    \end{aligned}
\end{equation}

We apply now \cite[Lemma 1]{schmidt2011convergence} to (\ref{ineq:zeta-lambda}), and obtain
\begin{equation}\label{ineq:boundvx}
    \begin{aligned}
        &  
       \Eb\left[\frac{\zeta^{k+1}}{2\lambda^{k+1}}\|x^{\star}-v_i^{k+1}\|^{2}\right] 
        \leq  \left(\sqrt{\frac{\zeta^{0}+\delta}{2}\|x^\star-x_i^0\|^2}+\sum_{t=0}^{k}\left( \Eb\left[\frac{2}{\delta\lambda^{t+1}} \|e_i^{t}\|^2\right]\right)^{1/2}\right)^{2}. \\
    \end{aligned}
\end{equation}

Next, we bound $\Eb[({\zeta^{k}}/{2})\|v_i^{k}-x^{\star}\|^{2}]$  separately  for the two cases $\mu > 0$ and $\mu = 0$, utilizing the bounds of $\Eb[\|e_i^{k}\|^2]$ derived   in  (\ref{eq:final_bound_ei_scvx}) and (\ref{eq:final_bound_ei_cvx}), respectively. 
 
    \begin{itemize}
        \item[\bf(i)]  $\mu > 0$:  Using $\alpha^{k}\equiv \alpha$,   $\zeta^{k} \equiv \muM$, and 
     
          (\ref{ineq:boundvx}), we have 
  
\begin{equation}\label{eq:final-bound-v-scvx}
\begin{aligned}
   & \overm\summ\Eb\left[\frac{\muM}{2}\|v_i^{k}-x^{\star}\|^{2}\right]\leq \left(\sqrt{\frac{\muM+\delta}{2}\|x^\star-x_i^0\|^2}+\sum_{t=0}^{k-1}\Eb\left[\left(\frac{\alpha\|e_i^{t}\|}{\sqrt{\lambda^{t+1}}}\sqrt{\frac{2}{\muM}}\right)^2\right]^{1/2}\right)^{2} \\
&\overset{(\ref{eq:final_bound_ei_scvx}) }{\leq} (1-c\alpha)^{k}\overm\summ\left(\sqrt{\frac{\muM+\delta}{2}\|x^\star-x_i^0\|^2}+\frac{\alpha\sqrt{125\delta^2 m\ep_{0}}}{ 1-c\alpha}\sqrt{\frac{2}{\muM}}\sum_{t=0}^{\infty}\left(\sqrt{\frac{1-\alpha}{1-c\alpha}}\right)^{t}\right)^{2}\\
    & = (1-c\alpha)^{k}\overm\summ \left(\sqrt{\frac{\muM+\delta}{2}\|x^\star-x_i^0\|^2}+\sqrt{\frac{2}{\muM(1-c\alpha)}}\frac{\alpha\sqrt{125\delta^2 m\ep_{0}}}{\sqrt{1-c\alpha}-\sqrt{1-\alpha}}\right)^{2} \\
    & \leq (1-c\alpha)^{k}\overm\summ\left((\muM+\delta)\|x^{\star}-x_i^0\|^2+\frac{2000\delta^2 m\epsilon_{0} }{\muM(1-c)^2}\right).
\end{aligned}
 \end{equation}
Denote  $c_{v,scvx} := \overm\summ\left((\muM+\delta)\|x^{\star}-x_i^0\|^2+\frac{2000\delta^2 m\epsilon_{0} }{\muM(1-c)^2}\right)$.

\item[\bf (ii)] $\mu=0$:    
    Using (\ref{ineq:boundvx}) and (\ref{eq:final_bound_ei_cvx}), yields 
\begin{equation}\label{eq:final-bound-v-cvx}
    \begin{aligned}
    &\overm\summ\Eb\left[\frac{\zeta^{k}}{2}\|v_i^{k}-x^{\star}\|^{2}\right] \\
    & \leq \lambda^{k}\overm \summ \left( \sqrt{\frac{\zeta^{0}+\delta}{2}\|x^\star-x_i^0\|^2} + \sum_{t=0}^{k-1} \left(\Eb\left[\frac{2}{\delta\lambda^{t+1}}\|e_i^{t}\|^2\right]\right)^{1/2}
 \right)^2  \\
& \leq \lambda^{k}\overm \summ \left( \sqrt{\frac{\zeta^{0}+\delta}{2}\|x^\star-x_i^0\|^2} + \sum_{t=0}^{k-1} \sqrt{\frac{250\delta m \epsilon_0}{\lambda^{t+1}}}\left(\frac{1}{t+1}\right)^{2+r_0}\right)^2 \\
    & \leq  \lambda^{k}\overm\summ\left( \sqrt{\frac{\zeta^{0}+\delta}{2}\|x^\star-x_i^0\|^2}+\sum_{t=0}^{k-1}  3\sqrt{125\delta m\epsilon_0} \left(\frac{1}{t+1}\right)^{1+r_0}\right)^2\\ 
    & \leq \frac{4}{(k+2)^2}\frac{1}{m}\summ\left(\sqrt{\frac{\zeta^{0}+\delta}{2}\|x^\star-x_i^0\|^2}+ 3\sqrt{125\delta m\epsilon_0} \left(\int_{0}^{\infty}\left(\frac{1}{t+1}\right)^{1+r_0}{\rm d}t+1\right)\right)^2  \\
    & = \frac{4}{(k+2)^2}\frac{1}{m}\summ\left(\sqrt{\frac{\zeta^{0}+\delta}{2}\|x^\star-x_i^0\|^2}+ 3\sqrt{125\delta m\epsilon_0} \left(\frac{1}{r_0}+1\right)  \right)^2  \\
    & \leq \frac{1}{(k+2)^2}\frac{1}{m}\summ\left(8 \delta \|x^\star-x_i^0\|^2 +9000 \delta m\epsilon_0\left(\frac{1}{r_0}+1\right)^2 \right),
    \end{aligned}
 \end{equation}
where in the last inequality we used   $\zeta^{0} = \delta$. 
Let  $c_{v,cvx} :=   ({1}/{m})\summ(8 \delta \|x^\star-x_i^0\|^2 +9000 \delta m\epsilon_0({1}/{r_0}+1)^2)$. 
    \end{itemize}
    
 Equipped with the bound on    
$(1/m)\summ\Eb[|\epsilon_{tot,i}^{k}|]$,   we are ready to prove convergence of  the outer loop  of Algorithm 1 by applying Proposition~\ref{prop:preliminary}.

\subsubsection{Application of Proposition~\ref{prop:preliminary}}

We separate the two cases $\mu > 0$   and $\mu = 0$.  

\begin{itemize}
    \item[\bf (i)]  $\mu > 0$: 
 Using (\ref{ineq:ep-tot}), (\ref{eq:final_bound_ei_scvx}) and (\ref{eq:final-bound-v-scvx}), we have
 $$
\begin{aligned}
\overm\summ\Eb[|\ep_{tot,i}^{k}|]  \leq \alpha\sqrt{\frac{\delta^2 m \epsilon_0 c_{v,scvx}}{\muM}}(1-c\alpha)^{k-1}.
\end{aligned}
 $$
 Denote $C_{scvx} : = \alpha\sqrt{\frac{\delta^2 m \epsilon_0 c_{v,scvx}}{\muM}}$.
By Corollary~\ref{prop:preliminary}, 
 $(1/m)\Eb[\|\vx^{k+1}-\vx^{\star}\|^{2}]\leq c_{scvx}(1-c\alpha)^{k+1}, 
 $
 where  
\begin{equation}\label{cal:c_scvx}
 \begin{aligned}
     & c_{scvx} = \overm\summ\left(\frac{2}{\muM}\left(\psi_i^{0}(x^{\star})-\Muds+\frac{C_{scvx}}{\alpha(1-c)}\right)\right) \\ 
     & \leq \overm\summ\left(\|x^\star-x_i^0\|^2+\frac{2}{\muM}(\Mud(x_i^0)-\Muds) + \frac{2}{\muM(1-c)}\sqrt{ \frac{\delta^2 m\epsilon_0 c_{v,scvx}}{\muM}}\right) \\ 
    & \leq \overm\summ\left(2+\frac{\delta}{\mu}\right)\|x^\star-x_i^0\|^2\!+\!\frac{2}{\muM(1-c)}\!\!\left(\!\sqrt{\!\frac{\delta^2\epsilon_0(\muM+\delta)}{\muM} \summ\|x-x_i^0\|^2}+\!\sqrt{\frac{2000\delta^4 m\epsilon_0^2}{\muM^2(1-c)^2}}\right)\\
    & \leq \overm\summ\left(\!2+\frac{\delta}{\mu}\!\right)\|x^\star-x_i^0\|^2+\frac{\delta(\muM+\delta)}{\muM^2(1-c)}\summ\|x^\star-x_i^0\|^2+\! \left(\!\frac{\delta}{\muM(1-c)}+\!\frac{2\sqrt{2000m}\delta^2}{\muM^2(1-c)^2}\right)\epsilon_0. 
 \end{aligned}
 \end{equation}
    \item[\bf (ii)]   $\mu = 0$: Following similar steps as for $\mu > 0$ and using  (\ref{ineq:ep-tot}), (\ref{eq:final_bound_ei_cvx}), and (\ref{eq:final-bound-v-cvx}), we have

 $$
\begin{aligned}
    \Eb[|\ep_{tot,i}^{k}|]& \leq  
\alpha^{k}\sqrt{\frac{2}{\zeta^{k+1}}}\left(\overm\summ\Eb[\|e_i^k\|^2]\right)^{1/2}\left(\frac{\zeta^{k+1}}{2m}\summ\Eb[\|x^{\star}-v_i^{k+1}\|^2]\right)^{1/2} 
\\
& \leq \sqrt{250\delta m\epsilon_0 c_{v,cvx}}\left(\frac{1}{k+1}\right)^{3+r_0}. \\ 
\end{aligned}        
 $$
Denote $C_{cvx} := \sqrt{250\delta m\epsilon_0 c_{v,cvx}}$. Applying Proposition~\ref{prop:preliminary},
we obtain  $$   
\!\!\!\!\!\!\!\!\!    \frac{1}{2\delta m}\summ\Eb[\|\nabla \Mud(x_i^{k})\|^{2}] \leq \frac{c_{cvx}}{(k+1)^{2}}\,\, \text{and}\,\,  \frac{1}{2\delta m}\summ\Eb[\|\nabla u(x_i^{k})\|^{2}] \leq \frac{c_{cvx}(L+\delta)^2}{\delta^2(k+2)^{2}} \,\,(\text{if }  r\equiv 0),$$
 where
\begin{equation}\label{cal:c_cvx}
\begin{aligned}
c_{cvx}  &= \frac{4}{m}\summ\left(\psi_i^{0}(x^{\star})-\Muds+\frac{10C_{cvx}}{r_{0}}\right)  \leq \frac{4}{m}\summ\left(\delta\|x^\star-x_i^0\|^2+\frac{10C_{cvx}}{r_{0}}\right).  
\end{aligned}
\end{equation}
\end{itemize}

\subsection{Analysis of the inner-loop of Algorithm~\ref{alg:framework}: Proof of Propositions~\ref{prop:inner-loop-scvx} and~\ref{prop:inner-loop-cvx}} 
\label{Subsec:proof-inner-loop}

We prove the two propositions  showing that the   upper bounds on $T_k$ as given therein  make (\ref{ineq:prop-out-scvx-Tk}) and (\ref{ineq:prop-out-cvx-Tk}) hold, under  Assumption~\ref{assump:warm-start} and the setting of  Propositions~\ref{prop:outer-loop-scvx} and~\ref{prop:outer-loop-cvx}, respectively.  We study separately the case $\mu>0$ and $\mu=0$.

 \subsubsection{\texorpdfstring{Inner loop analysis ($\mu>0$):  Proof of Proposition~\ref{prop:inner-loop-scvx}}{Inner loop analysis for mu-strongly convex case: Proof of Proposition}}
 \label{subsubsec:inner-loop-scvx} 

For simplicity, denote 
$$
c_{T} := \left\lceil r_{\cM,\delta}\cdot \max\left\{   \log\frac{c_{\cL}\cL^{0}(\vs^{0,0})}{\epsilon_{0}(1-c\alpha)},  
 \,\, \log\frac{  c_{\cL} c_{\cM}  + 36 d_{\cM}c_{scvx}\frac{1}{\epsilon_{0}(1-c\alpha)^{2}} }{1-c\alpha}\right\}\right\rceil.
$$
We proceed by induction, with the following  induction hypothesis: given $k\geq 0$ and  any   $\ell\in[k]$, there exists $T_{\ell}\in (0,c_T]$ such that   (\ref{ineq:prop-out-scvx-Tk}) holds for the $\ell$th outer loop with such $T_{\ell}$. 

\begin{itemize}
    \item  $k = 0:$ It is sufficient to  pick $T_0 = \left\lceil r_{\cM,\delta}\log\left((c_{\cL}\cL^{0}(\vs^{0,0}))/(\epsilon_{0}(1-c\alpha))\right)\right\rceil \leq c_{T}$, for  (\ref{ineq:prop-out-scvx-Tk}) to hold, since
 $$  \Eb[\cL^{0}(\vs^{0,T_0})|\cF_0] \leq  c_{\cL}\left(1-\frac{1}{r_{\cM,\delta}}\right)^{T_0}\cL^{0}(\vs^{0,0})  \leq  \epsilon_{0}(1-c\alpha).    $$
 
    \item $k > 0:$ Suppose  the induction hypothesis holds at given $k$. we leverage  Assumption~\ref{assump:warm-start} to bound $\Eb[\cL^{k+1}(\vs^{k+1, 0})]$. 
    
    We preliminary bound 
\begin{equation}\label{boundZk}
\begin{aligned}
 \|\vz^{k+1}-\vz^{k}\|  
&\leq \left(1+\frac{1-\alpha}{1+\alpha}\right)\|\vx^{k+1}-\vx^{k}\|+\frac{1-\alpha}{1+\alpha}\|\vx^{k}-\vx^{k-1}\|\\
&\leq 3\max\{\|\vx^{k+1}-\vx^{k}\|,\|\vx^{k}-\vx^{k-1}\|\},\\
\end{aligned}
\end{equation}
Using  $\|\vx^{k+1}-\vx^{k}\|$ by 
$
\|\vx^{k+1}-\vx^{k}\| \leq \|\vx^{k+1}-\vx^{\star}\| + \|\vx^{\star} -\vx^{k}\|
$, we proceed  bounding $\|\vx^{\star}-\vx^{k}\|$. According to the induction hypothesis, (\ref{ineq:prop-out-scvx-Tk}) holds   for all $\ell \in [k]$. Therefore, invoking    Proposition~\ref{prop:outer-loop-scvx}   we have
$$
    \begin{aligned}
    &\Eb[\|\vx^{\star}-\vx^{\ell}\|^{2}]   \overset{(\ref{Prop1:result-scvx})}{\leq} mc_{scvx}(1-c\alpha)^{\ell},\quad\forall \ell\in[k+1]. 
      \end{aligned}
$$
This yields 
 $$
\begin{aligned}
    \Eb[\cL^{k+1}(\vs^{k+1,0})] &\leq c_{\cM} \epsilon_0(1-c\alpha)^{k+1} + 36 d_{\cM}c_{scvx}(1-c\alpha)^{k-1}.\\ 
\end{aligned}
 $$ 
Choose $$T_{k+1} = \left\lceil r_{\cM,\delta}     
 \log\frac{  c_{\cL} c_{\cM}  + 36 d_{\cM}c_{scvx}\frac{1}{\epsilon_{0}(1-c\alpha)^{2}} }{1-c\alpha}\right\rceil \leq c_{T}.$$
 Clearly,  (\ref{ineq:prop-out-scvx-Tk}) holds with
$$
    \begin{aligned}
        & \Eb[\cL^{k+1}(\vs^{k+1,T_{k+1}})] \leq  c_{\cL}\left(1-\frac{1}{r_{\cM,\delta}}\right)^{T_{k+1}}\Eb[\cL^{k+1}(\vs^{k+1,0})] \leq \epsilon_0(1-c\alpha)^{k+2}. 
    \end{aligned}
 $$

We proved that   the minimal $T_k$ for (\ref{ineq:prop-out-scvx-Tk}) to hold is upper bounded by $c_{T}$,  for all  $k\geq 0$.  
\end{itemize}

\hfill$\square$

\subsubsection{Inner Loop Analysis \texorpdfstring{$(\mu=0)$}{(mu=0)}: proof of Proposition~\ref{prop:inner-loop-cvx}}\label{subsubsec:inner-loop-cvx}

The procedure is similar to that in Sec.~\ref{subsubsec:inner-loop-cvx}. Let 
$$
c_{T_k} := \left\lceil r_{\cM,\delta}\cdot\max\left\{  
         \log\frac{ c_{\cL}9^{r_0+2}\cL^{0}(\vs^{0,0})}{\epsilon_0}, 
        \log \left(c_{\cL} c_{\cM}2^{4+2r_{0}} +\frac{36c_{\cL} d_{\cM}B^2(k+3)^{4+2r_{0}}}{\epsilon_{0}}\right)\right\}\right\rceil. 
$$
The induction hypothesis reads: given $k\geq 0$ and  $\ell\in[k]$, there exists $T_{\ell} \in (0,c_{T_{\ell}}]$ such that  (\ref{ineq:prop-out-cvx-Tk}) holds for the $\ell$th outer loop with such $T_{\ell}$.

\begin{itemize}
    \item $k=0$: Choose   $T_0 = \left\lceil r_{\cM,\delta}\log\left(( c_{\cL}9^{r_0+2}\cL^{0}(\vs^{0,0}))/\epsilon_0\right)\right\rceil \leq c_{T_0}$. It holds 
$$
\Eb[\cL^{0}(\vs^{0,T_0})|\cF_0] \leq c_{\cL}\left(1-\frac{1}{r_{\cM,\delta}}\right)^{T_0}\cL^{0}(\vs^{0,0}) \leq  \frac{1}{3^{4+2r_0}}\epsilon_0.
$$
 \item $k>0:$ Using   (\ref{boundZk}) and $\norm{x_i^k}\leq B$,  we have  
 \begin{equation}\label{ineq:zk diff}
     \Eb\left[\frac{1}{m}\|\vz^{k+1}-\vz^{k}\|^{2}\right] \leq 36B^{2},\quad \forall k\geq 0. 
 \end{equation}
  Then, by Assumption~\ref{assump:warm-start} and  (\ref{ineq:zk diff}), 
 $$
    \begin{aligned}
          \Eb[\cL^{k+1}(\vs^{k+1,0})] &\leq c_{\cM}\Eb[\cL^{k}(\vs^{k,T_{k}})]+\frac{d_{\cM}}{m}\Eb[\|\vz^{k}-\vz^{k+1}\|^2] \\
        & \leq c_{\cM} \left(\frac{1}{k+3}\right)^{4+2r_0}\epsilon_0 + 36 d_{\cM} B^2.
    \end{aligned}
 $$
Choosing $$T_{k+1} = \left\lceil r_{\cM,\delta}\log \left(c_{\cL} c_{\cM}2^{4+2r_{0}} +\frac{36c_{\cL} d_{\cM}B^2(k+4)^{4+2r_{0}}}{\epsilon_{0}}\right)\right\rceil \leq c_{T_{k+1}},$$ 
makes  (\ref{ineq:prop-out-cvx-Tk}) hold, with
 $$
 \begin{aligned}
     &\Eb[\cL^{k+1}(\vs^{k+1,T_{k+1}})] \leq c_{\cL} \Eb[\cL^{k+1}(\vs^{k+1,0})]\left(1-\frac{1}{r_{\cM,\delta}}\right)^{T_{k+1}}  \leq \left(\frac{1}{k+4}\right)^{4+2r_0}\epsilon_0.
 \end{aligned}
 $$
\end{itemize}
 
\hfill$\square$

\section{Numerical Experiments}\label{section:simulations}

This section presents experiments conducted on  
real datasets across three instances of Problem~(\ref{eq:problem}): strongly convex 
$u$, convex (non-strongly convex) 
$u$, and losses 
$f_i$  with a finite-sum structure. Given that DCatalyst introduces acceleration for minimizing {\it composite} objective functions (i.e., $r\neq 0$) across all these classes for the first time,  our experiments predominantly focus on these types of functions.  To offer a comprehensive comparison against existing accelerated decentralized methods, which are generally designed for{\it  smooth}, unconstrained optimization problems, we have included additional experiments for smooth instances of Problem~\eqref{eq:problem} in the arXiv version of this paper.  

Unless otherwise specified, the setup for all experiments is as follows: we simulate an     Erdos-Renyi graph with $m=30$ nodes (agents) and an edge probability of $p=0.5$. The gossip weight matrix used in all the algorithms is Metropolis-Hasting weight matrix. For strongly convex problems, the optimality gap at iteration  $k$ is defined as $\frac{1}{m}\sum_{i=1}^{m}\|x_i^k-x_{\text{opt}}\|^2$, while for (non strongly) convex functions,the gap is measured using   $\frac{1}{m}\sum_{i=1}^{m}u(x_i^k)-u_{\text{opt}}$.

\subsection{Strongly Convex objectives}
\label{subsec:strongly cvx}

Our first experiment concerns  the logistic regression model with the elastic net sparsity-inducing regularization. Beyond solely $\ell_1$ penalty (LASSO), the elastic net regularization is proposed for a better selection of groups of correlated variables in statistical learning problems, which can help improve the prediction accuracy  on many real datasets--see, e.g.,~\cite{zou2005regularization,tay2023elastic}.  The formulation of interest corresponds to Problem~\eqref{eq:problem} with 
\begin{equation}\label{eq:logistic-scvx}f_i(x)=\frac{1}{n}\sum_{j=1}^{n}\log(1+\exp(-b_{ij}\cdot\langle x,a_{ij}\rangle))+\frac{\gamma}{2}\|x\|_2^2\quad\text{and}\quad r(x)=\lambda\|x\|_1,\end{equation}
where  $a_{ij}\in \mathbb{R}^{d}$ and $b_{ij}\in\{-1,1\}$. Here, the data set  $\{(a_{ij},b_{ij})\}_{j=1}^n$ is assumed to be private and  owned only by  agent $i$. We use MNIST dataset from LIBSVM \cite{refer_libsvm} of size $N=60000$ and feature dimension $d=784$.  

We contrast the proposed  DCatalyst-SONATA-L and DCatalyst-SONATA-F with DPAG \cite{ye2020decentralized}, which is the   only available decentralized algorithm  applicable to nonsmooth objective functions.  All algorithms were implemented according to their theoretical guidelines. For DCatalyst-SONATA-L (resp.  DCatalyst-SONATA-F), we set    $\delta=L-\mu$ (resp. $\delta=\beta-\mu$), the momentum parameter $\alpha=\sqrt{\mu/(\mu+\delta)}$, and the number of inner-loop iterations $T_k=\lceil\log(L/\mu)\rceil$ (resp. $T_k=\lceil\log(\beta/\mu)\rceil$), for all $k$. In  DCatalyst-SONATA-L, the local agents' prox-updates are computed in closed-form while in   DCatalyst-SONATA-F,   
the  accelerated proximal gradient method is employed, rub  up  to a tolerance of   $10^{-8}$. For the DPAG, we followed the tuning recommendation as in~\cite{ye2020decentralized}.

\begin{figure}[htbp] 
\centering
\begin{minipage}[t]{0.48\textwidth} 
\centering
\includegraphics[width=6cm]{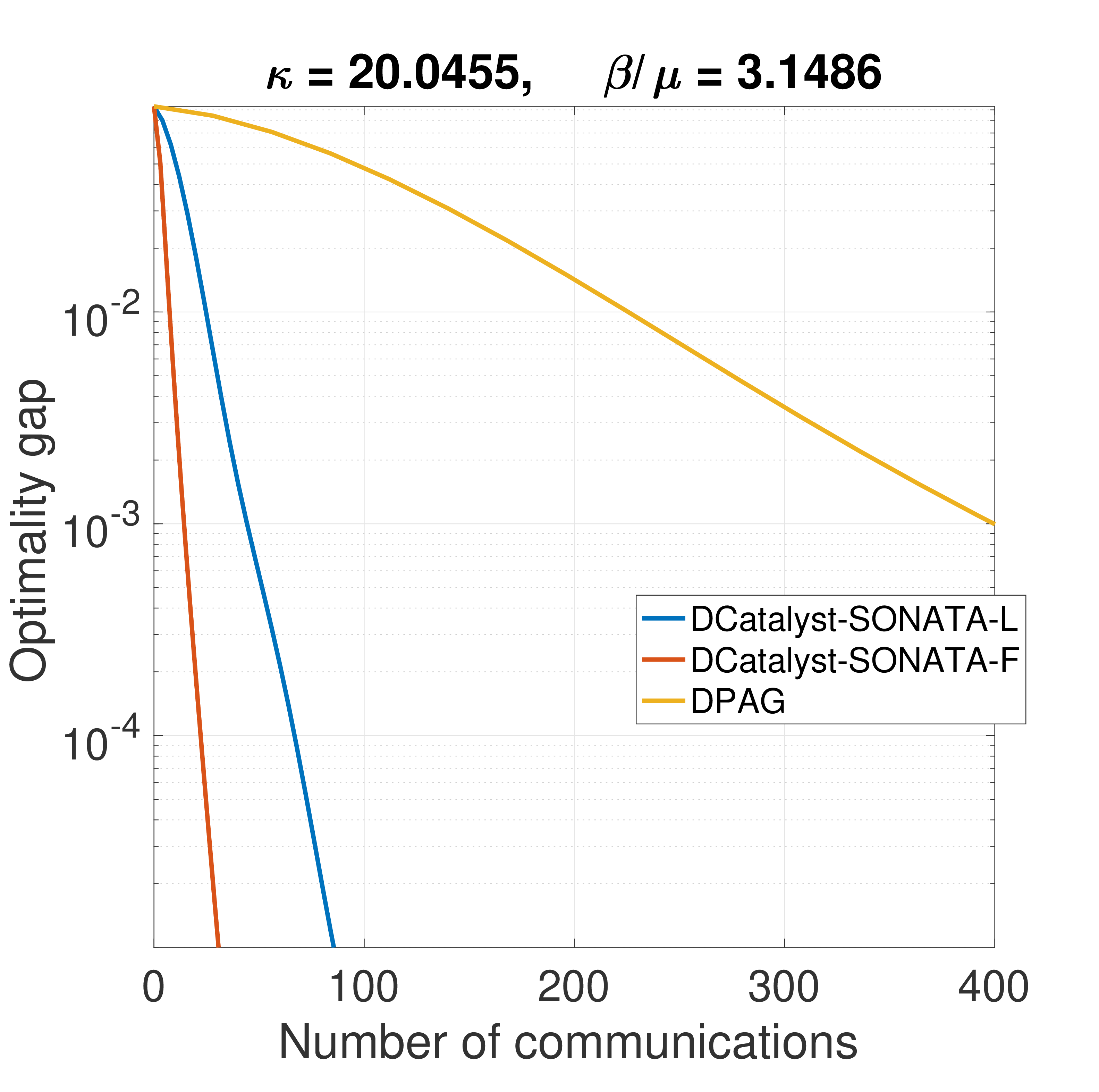}
\end{minipage}
\begin{minipage}[t]{0.48\textwidth}
\centering
\includegraphics[width=6cm]{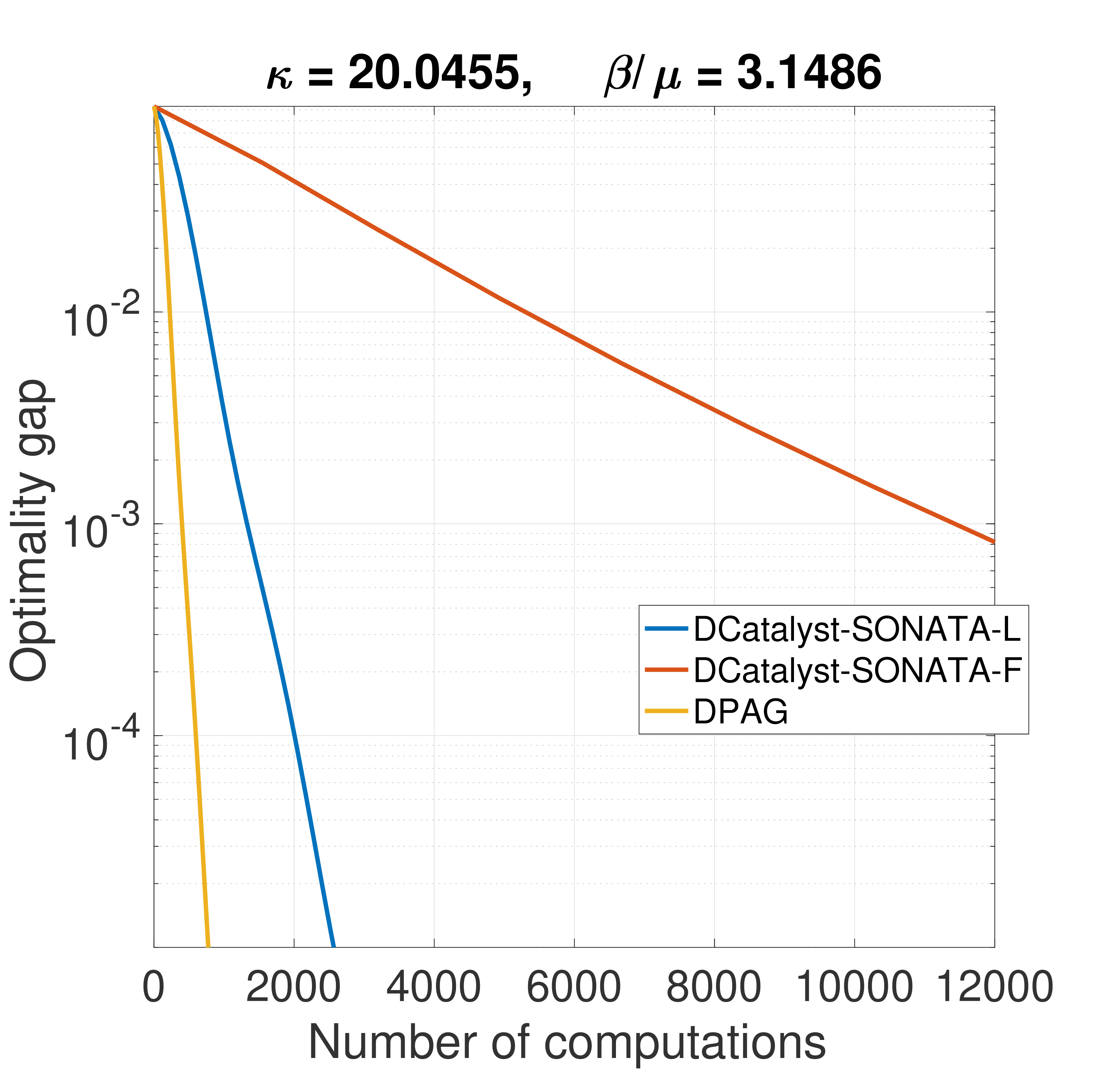} 
\end{minipage}
    \caption{ Comparison of distributed algorithms under strongly convex and non-smooth setting in \textbf{(a):} communication cost (left) and \textbf{(b):} computation cost (right).}
    \label{'6_1_2_1'} 
\end{figure}

$\bullet$ {\bf On the linear convergence:} The initial experiment aims to validate the linear convergence rate predicted by our theoretical findings for this class of problems.  In these experiments, we set    in (\ref{eq:logistic-scvx}),  $\lambda=0.01$ and  $\gamma=0.5$, resulting   in $\beta/\mu=3.1$ and $\kappa_g=20$.   Figure~\ref{'6_1_2_1'} 
summarizes the comparison, plotting the optimality gap of each algorithm  versus the number of total communications (subplot (a)) and the number of inner-plus-outer iterations (subplot (a)) (for DCatalyst-SONATA-F, this includes also the number of iterations run by the   accelerated proximal gradient method employed locally to compute the local agents' prox-updates).  All
the schemes achieve linear convergence. Consistent with our theoretical predictions (see Corollary~\ref{coro:sonata-scvx}), DCatalyst-SONATA-F excels particularly when the local functions $f_i$ exhibit some similarity, quantified by      $1+\beta/\mu<\kappa_g$,    outperforming both   DCatalyst-SONATA-L and DPAG~\cite{ye2020decentralized}, which do not leverage function similarity.    Notably,  DCatalyst-SONAT-L significantly outperforms    DPAG~\cite{ye2020decentralized} in terms of communication rounds while maintaining comparability in the number of iterations. This is quite remarkable considering that DPAG is a {\it single-loop} scheme, specifically designed for this class of problems while DCatalyst-SONAT-L is the result of a general unified framework applicable to a much larger class of objective functions. 

\begin{figure}[htbp] 
\centering
\begin{minipage}[t]{0.3\textwidth}
\centering
\includegraphics[width=5cm]{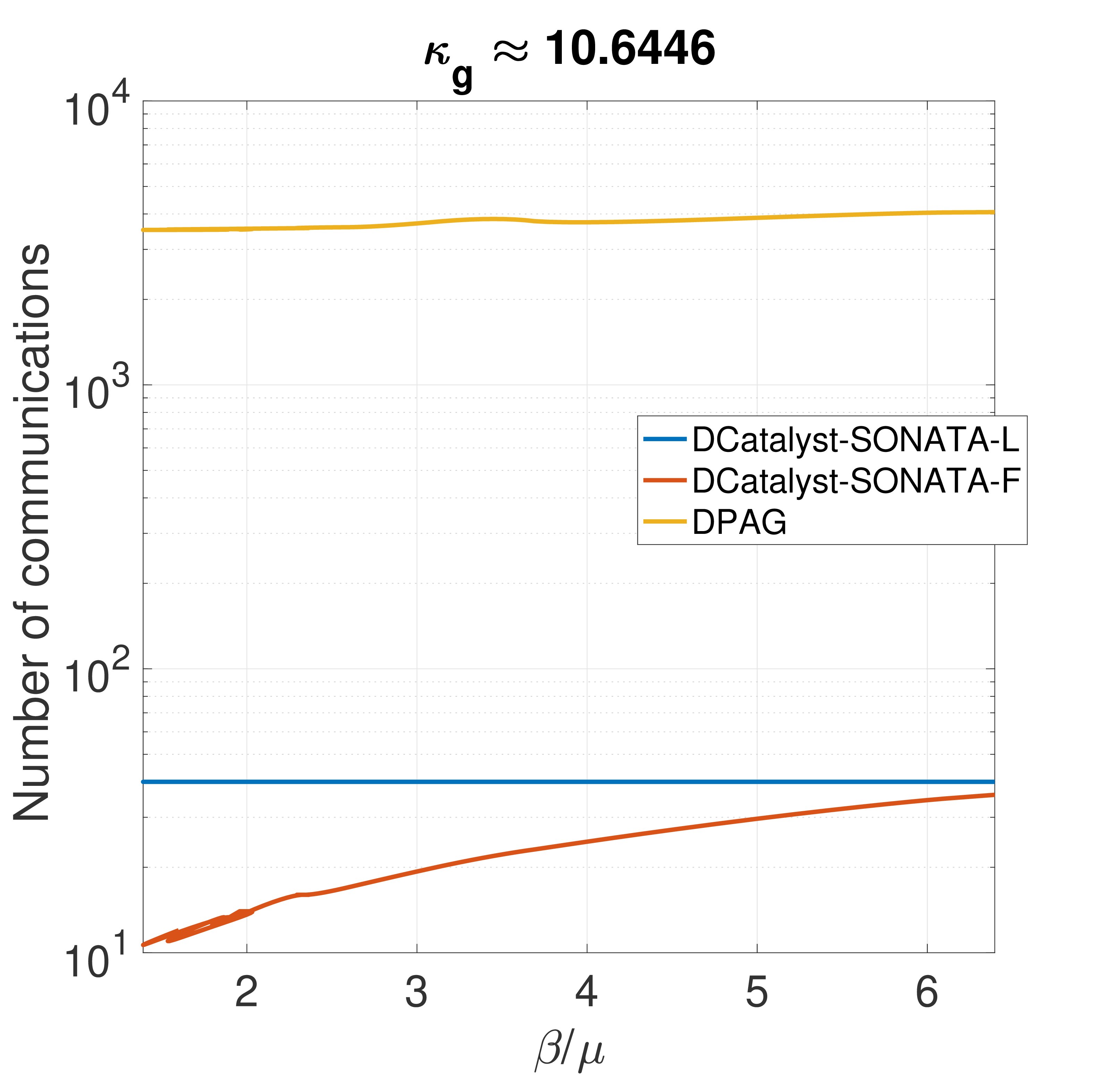}
\end{minipage}
\begin{minipage}[t]{0.3\textwidth}
\centering
\includegraphics[width=5cm]{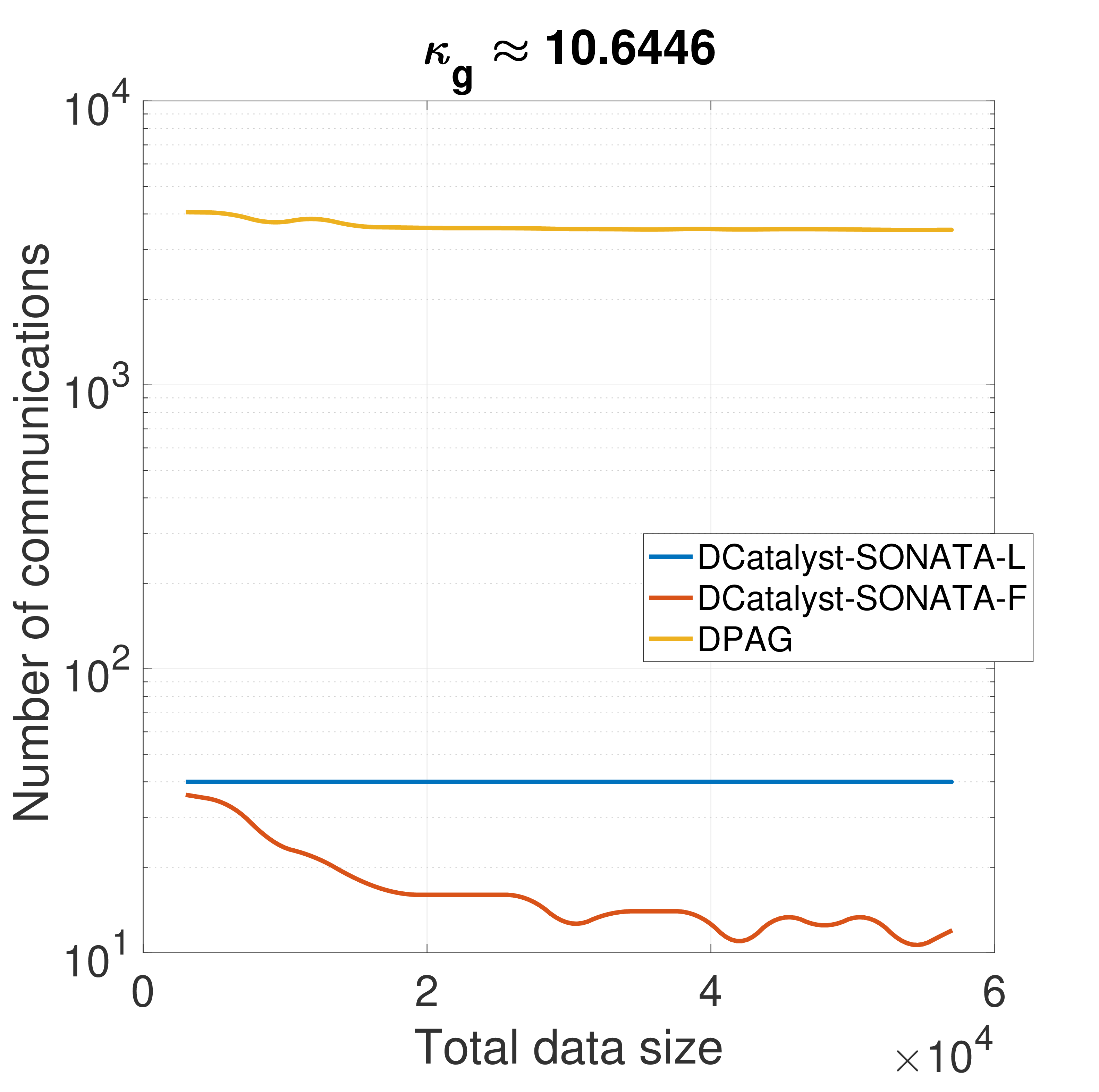}
\end{minipage}
\centering
\begin{minipage}[t]{0.3\textwidth}
\centering
\includegraphics[width=5cm]{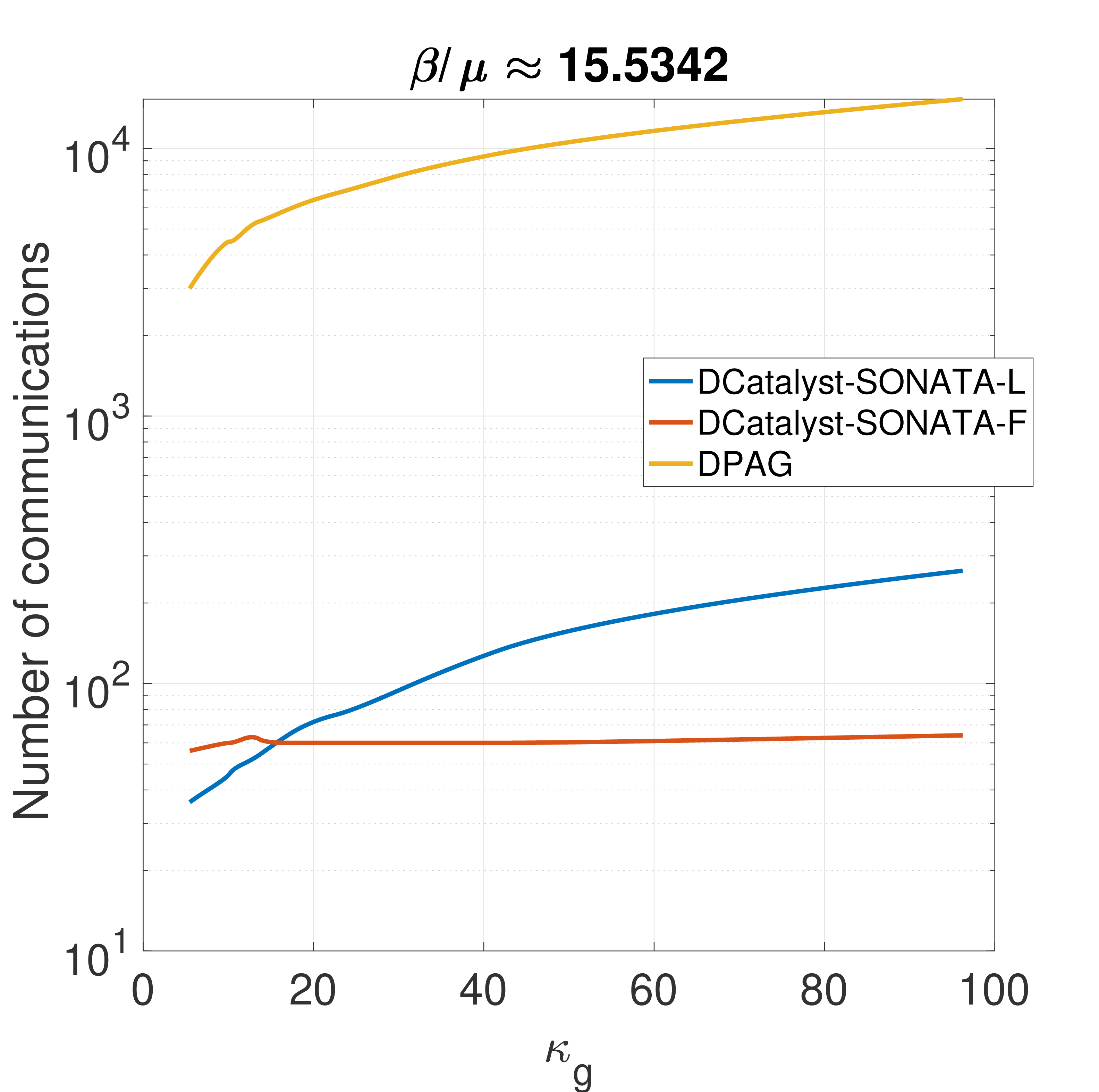}
\end{minipage} \caption{Comparison of distributed algorithms under strongly convex and non-smooth setting on the influence of parameters \textbf{(a):} similarity $\beta/\mu$ (left), \textbf{(b):} total sample size $N$ (middle) and  \textbf{(c):} global condition number $\kappa_g$  to the number communication rounds needed to reach a precision of $10^{-4}$ (right).}
    \label{'6_1_2_2'} 
\end{figure}

$\bullet$ {\bf  On the impact of $\kappa_g$, $\beta/\mu$, and  and total sample size.} We consider the following two scenarios. \textbf{(i)} {\it Changing $\beta/\mu$ with (almost) fixed $\kappa_g$}:  We generate instances of logistic  regression problem  with decreasing $\beta$ and (almost) fixed $\kappa_g$, increasing the local sample size $n$ (starting from $n=20$) while keeping  $\gamma=0.5$ fixed (and $\lambda=0.01$), resulting in a $\kappa_g\approx 10.646$. 
Fig.~\ref{'6_1_2_2'}  (left-panel) captures this scenario; we plot
the number of communications to drive the optimality gap below $10^{-4}$. In the mid-panel we report the same number of communications versus the total sample size $N$.  
 \textbf{(ii)} {\it Changing $\kappa_g$ with fixed $\beta/\mu$}:  We generate instances of logistic  regression problem  with varying   $\kappa_g$, acting on $\gamma$, while keeping $\beta/\mu$ constant by changing the local sample size to compensate for the variation of $\mu$ via $\gamma$, resulting in  $\beta/\mu\approx 15.53$. 
Fig.~\ref{'6_1_2_2'}  (right-panel) captures this scenario; we plot
the number of communications to drive the optimality gap below $10^{-4}$ versus $\kappa_g$, for fixed $\beta/\mu$.  

   The following comments are in order. The left panel confirms what predicted by Corollary~\ref{coro:sonata-scvx}: the convergence rate of DCatalist-SONATA-F  scales with $\sqrt{\beta/\mu}$ while that of the other reported algorithms is almost invariant with ${\beta/\mu}$. This is because  those other methods use only gradient information  and hence  cannot benefit from statistical similarity. On the other hand, the right-panel shows that DCatalist-SONATA-L and DPAG exhibit  a communication complexity that  deteriorates when  $\kappa_g$ grows (of the order of $\sqrt{\kappa_g}$) whereas that of DCatalist-SONATA-F  remains almost invariant. This is exactly what our theoretical results predicted. Notice also that both instances of the proposed framework uniformly outperform DPAG, in any simulated setting.

\subsection{(Non-strongly) Convex Setting}
\label{subsec:cvx}
As (non strongly) convex, nonsmooth instance of Problem~\eqref{eq:problem}, we consider the decentralized logistic regression problem with $\ell_1$-regularization. This corresponds to the formulation in (\ref{eq:logistic-scvx}), with $\gamma=0$. We set $\lambda=10^{-4}$.   
 
\begin{figure}[htbp] 
\centering
\begin{minipage}[t]{0.48\textwidth}
\centering
\includegraphics[width=6cm]{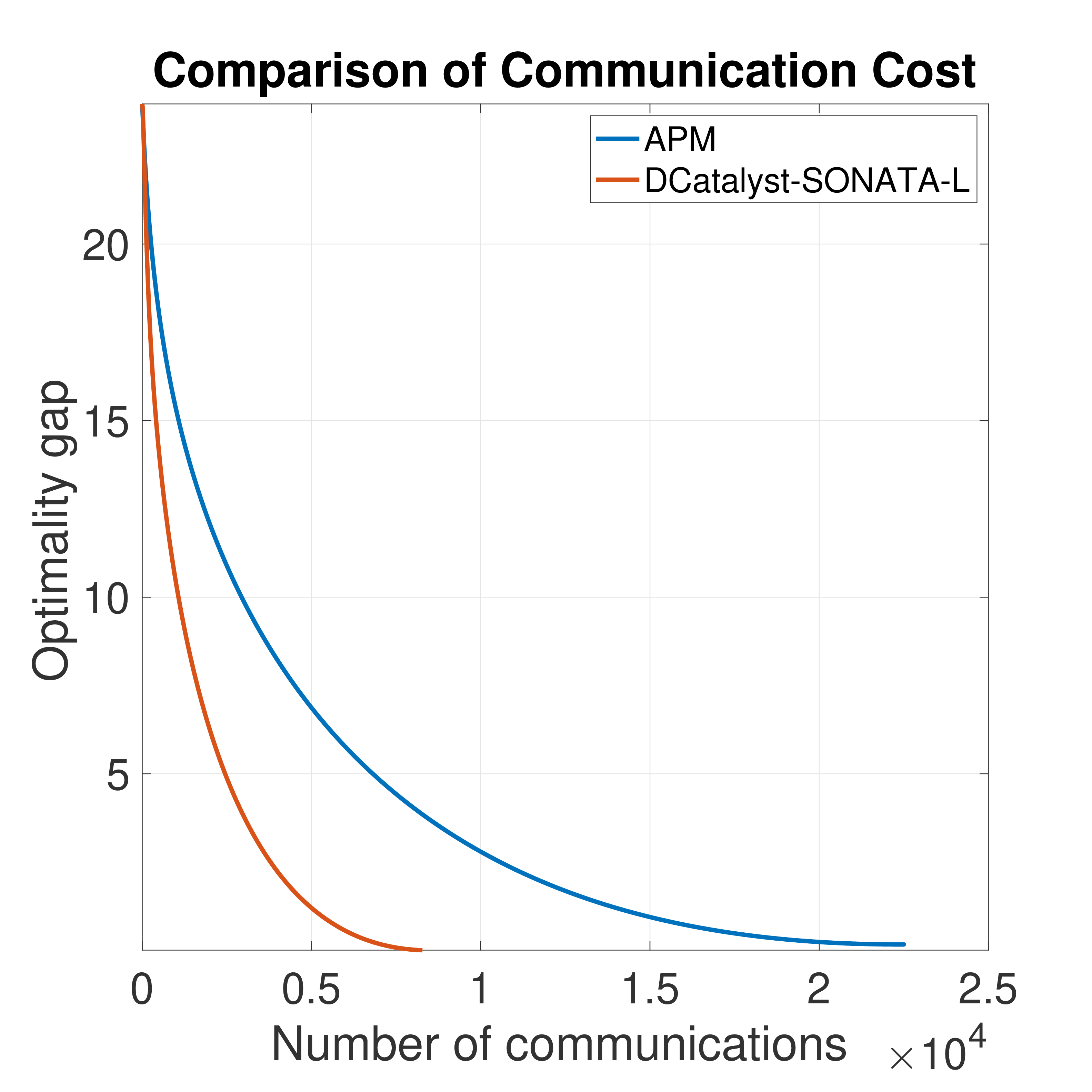}
\end{minipage}
\begin{minipage}[t]{0.48\textwidth}
\centering
\includegraphics[width=6cm]{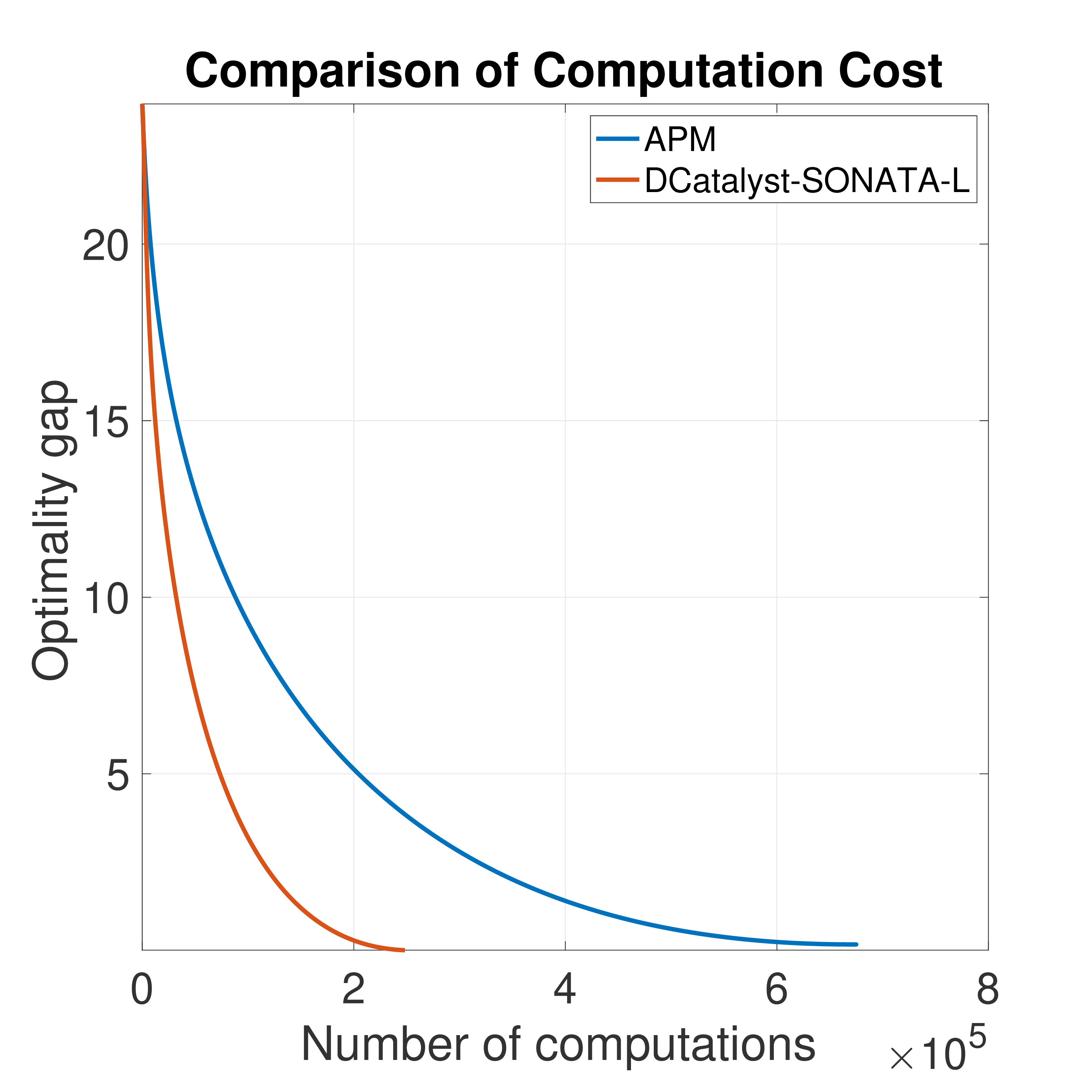}
\end{minipage}
\caption{The comparison between APM\cite{li2020decentralized} and DCatalyst-SONATA-L under convex and non-smooth setting in \textbf{(a):} communication cost (left) and \textbf{(b):} computation cost (right).}
\label{'6_2_2'} 
\end{figure} 

We contrast the proposed DCatalyst-SONATA-L with the APM algorithm~\cite{li2020decentralized},  which to our knowledge is the only algorithm available in the literature applicable to such a class of problems.  
The tuning of the free parameters in  APM   follows the theoretical guidelines in \cite{li2020decentralized}. For the DCatalyst-SONATA-L algorithm, the parameters are chosen according to the theory developed in Sec.~\ref{sec:convergence-cvx}: we set $\delta=L$, the momentum parameter is calculated using~\eqref{eq:alpha-rec}, and the number of inner-loop iterations is set as $\lceil\log(k+1)\rceil$ for all $k$ [see~\eqref{inner-lp-cvx}]. 

In Figure \ref{'6_2_2'}, we plot the optimality gap achieved by the two algorithms versus the number of communications (left-panel) and computations (right-panel). The figure  confirms the sublinear convergence of the algorithms, with DCatalyst-SONATA-L outperforming APM \cite{li2020decentralized} both on  communication and computation costs.

\subsection{Minimizing finite-sum functions via variance-reduction methods} 
\label{subsec:fns} 
The last set of experiments concerns the application of acceleration to decentralized variance reduction methods, testing the first scheme of this kind applicable to composite functions $u$ (whose $f$-part has a finite-sum structure).    To offer a comprehensive comparison against existing accelerated decentralized variance reduction  methods, which are available only for {\it  smooth}, unconstrained optimization problems, we also present in Sec.~\ref{subsubsec:fns+sm},  some  experiments for smooth functions $u$.  
\subsubsection{Logistic regression  with    elastic net  regularization}
We consider the logistic regression problem introduced in Sec.~\ref{subsec:strongly cvx} (see \eqref{eq:logistic-scvx}), but now exploiting algorithmically the finite-sum structure of each agent's loos $f_i$, that is,  
$$f_i(x)=\sum_{j=1}^n  f_{ij}(x), \quad \text{with}\quad f_{ij}(x)=\log(1+\exp(-b_{ij}\cdot\langle x,a_{ij}\rangle))+\frac{\gamma}{2}\|x\|_2^2,\quad r(x)=\lambda\|x\|_1.$$ 
In the experiments, we set $\lambda=10^{-7}$ and $\gamma=10^{-4}$, and use only the first   $N=6000$  data points of the data set MNIST. This  ensures   $n=200<<\kappa_s\approx2.73\times10^5$.  

\begin{figure}[htbp]
\centering
\begin{minipage}[t]{0.48\textwidth}
\centering
\includegraphics[width=6cm]{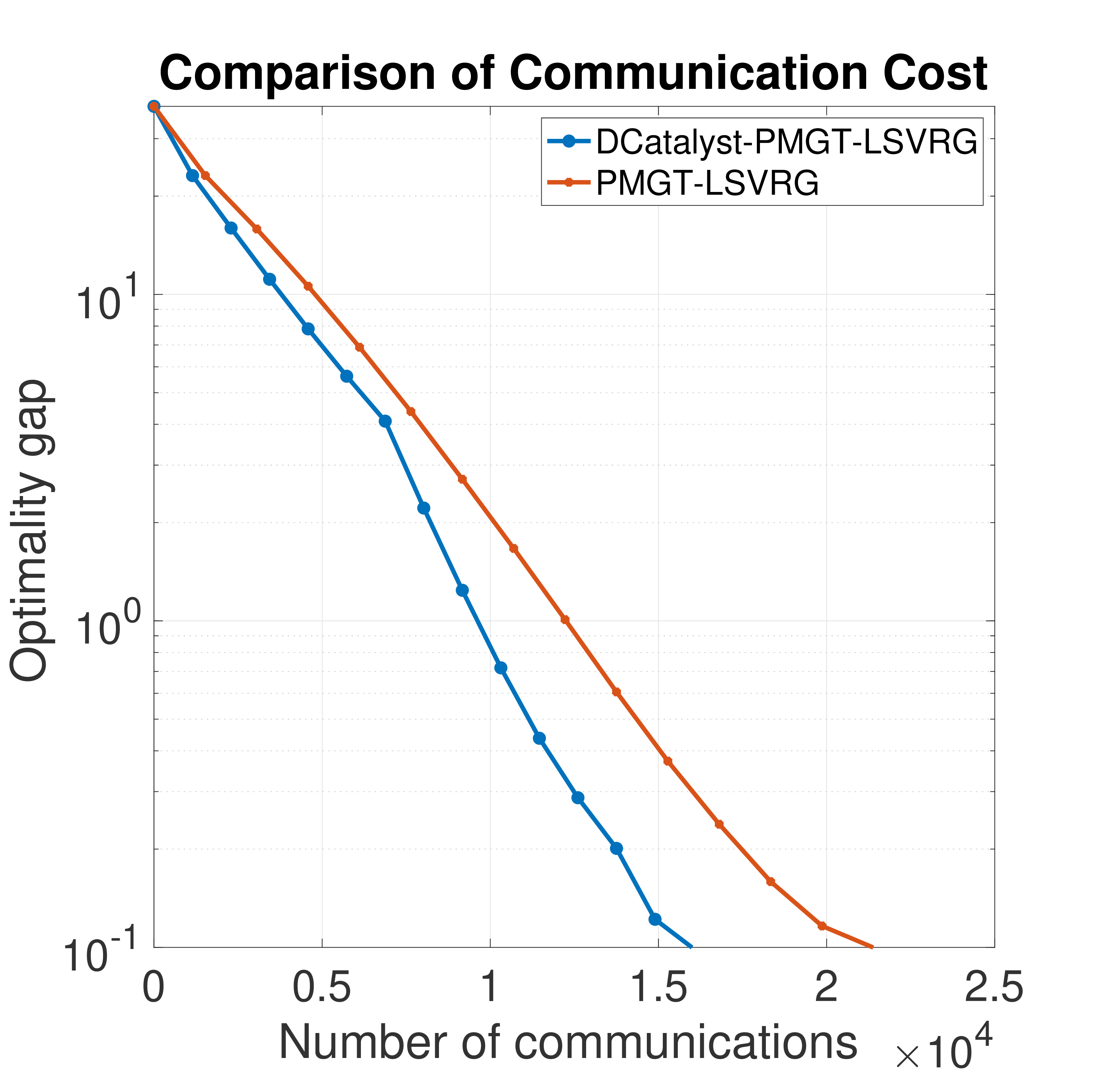}
\end{minipage}
\begin{minipage}[t]{0.48\textwidth}
\centering
\includegraphics[width=6cm]{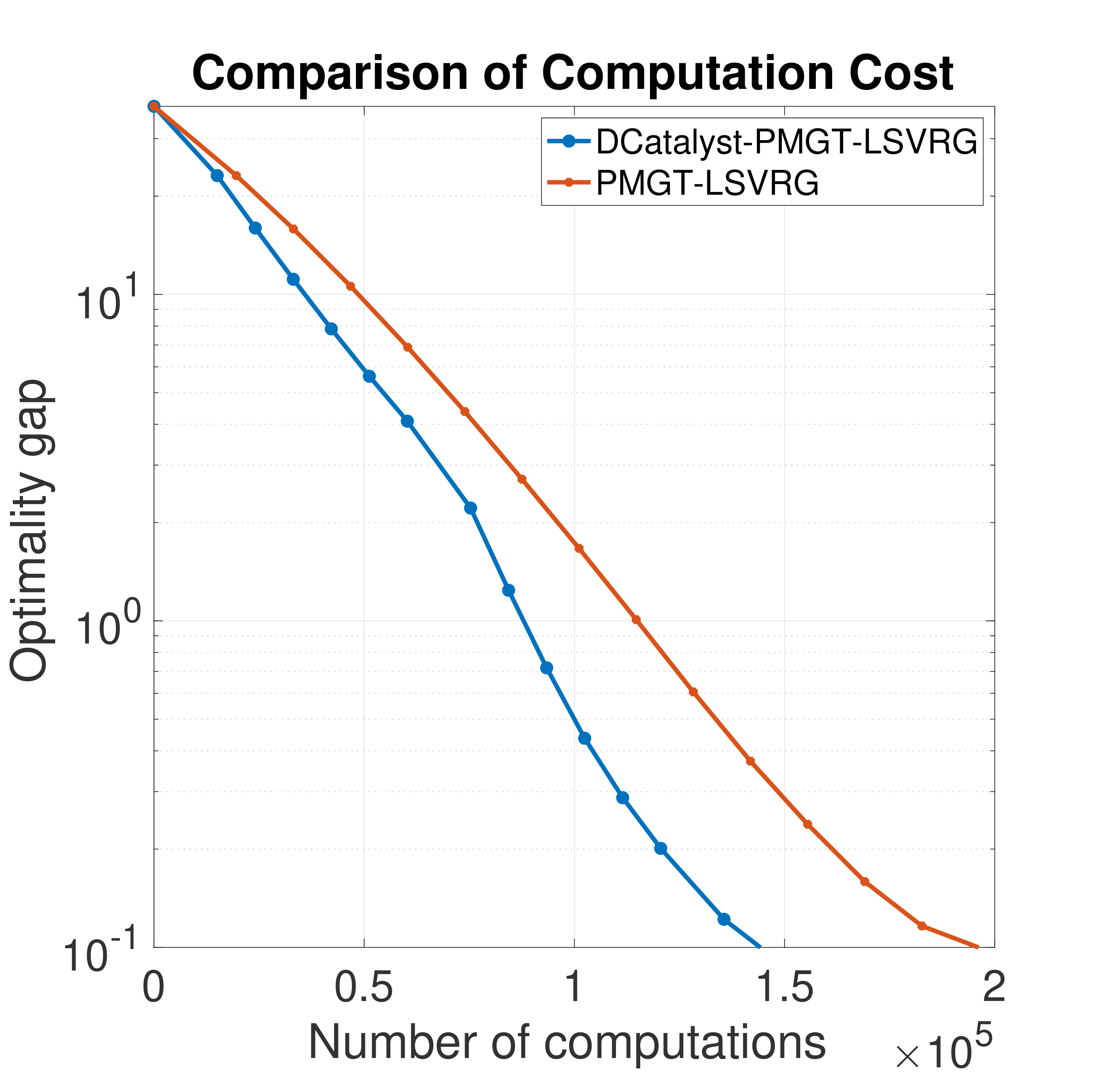}
\end{minipage}
\caption{Comparison between PMGT-LSVRG\cite{ye2021pmgtvr} and DCatalyst-PMGT-LSVRG, under the finite-sum setting, with   strongly convex non-smooth $u$ and  $n<<\kappa_s$.   \textbf{(a):} communication cost (left);  \textbf{(b):} computation cost  (right).}
\label{'6_3_1'} 
\end{figure}

Since for such classes of functions there is no accelerated variance reduction decentralized methods in the literature,  we compare the non-accelerted decentralized algorithm  PMGT-LSVRG~\cite{ye2021pmgtvr} with its accelerated counterpart obtained applying our DCatalyst framework, termed DCatalyst-PMGT-LSVRG,  The tuning parameters of PMGT-LSVRG are set  as recommended in~\cite{ye2021pmgtvr}. For DCatalyst-PMGT-LSVRG, similar to our previous experiments, we set  $\delta=L-\mu$, the momentum parameter 
$\alpha=\sqrt{\mu/(\mu+\delta)}$,  and the number of inner loop iterations 
$T_k=\lceil\log(L/\mu)\rceil$ Figure~\ref{'6_3_1'} plots the optimality gap versus the number of communications and iterations (evaluated in terms of single gradient computation $\nabla f_{ij}$) produced by the two algorithms. The plots show that,  when $n<<\kappa_s$, the acceleration introduced by   DCatalyst  improves both   computation and communication complexities of the plain PMGT-LSVRG algorithm.  

\subsubsection{Ridge regression }
\label{subsubsec:fns+sm} 
As instance of  strongly convex and smooth objective function with sum-structure we consider here   the ridge regression model, which corresponds to Problem~\eqref{eq:problem} with  
 $$f_i(x)=\sum_{j=1}^n f_{ij}(x),\quad \text{with}\quad  f_{ij}(x)=\frac{1}{2}\|a_{ij}^Tx-b_{ij}\|^2+\frac{0.1}{2}\|x\|^2,\quad r(x)=0. $$
Here,  $\{(a_{i,j},b_{ij})\}_{j=1}^n$, $a_{ij}\in\mathbb{R}^d$, $b_{ij}\in\mathbb{R}$, is the local data set accessible only to agent $i$.  

We simulate two notable decentralized algorithms that utilize variance reduction techniques and applicable to smooth functions: the  Acc-VR-EXTRA-CA  algorithm and the Acc-VR-DIGing-CA algorithm~\cite{li2022variance}. We compare these with our DCatalyst framework applied to the (non-accelerated) VR-EXTRA~\cite{li2022variance}, which we term DCatalyst-VR-EXTRA. For Acc-VR-EXTRA-CA and Acc-VR-DIGing-CA, we adhere to the tuning recommendations provided in~\cite{li2022variance}.       For DCatalyst-VR-EXTRA, similar to our previous experiments, we set  $\delta=L-\mu$, the momentum parameter 
$\alpha=\sqrt{\mu/(\mu+\delta)}$,  and the number of inner loop iterations 
$T_k=\lceil0.5\log(L/\mu)\rceil$. Additionally, to illustrate the balance between communication and computation costs, we experiment with different batch sizes:   $b=1$, $b=100$, and $b=1000$.

\begin{figure}[htbp] 
\centering
\begin{minipage}[t]{0.48\textwidth}
\centering
\includegraphics[width=6cm]{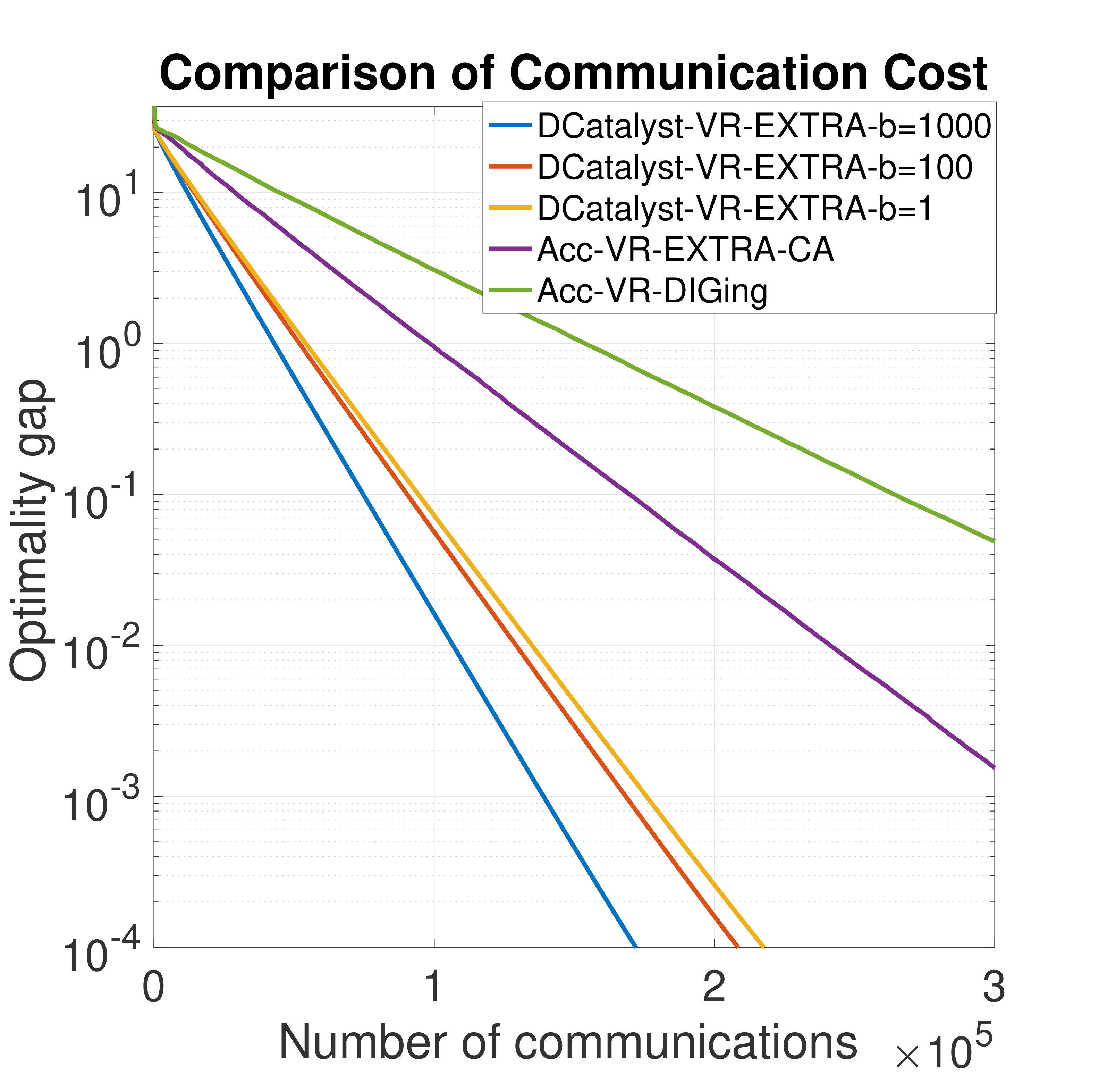}
\end{minipage}
\begin{minipage}[t]{0.48\textwidth}
\centering
\includegraphics[width=6cm]{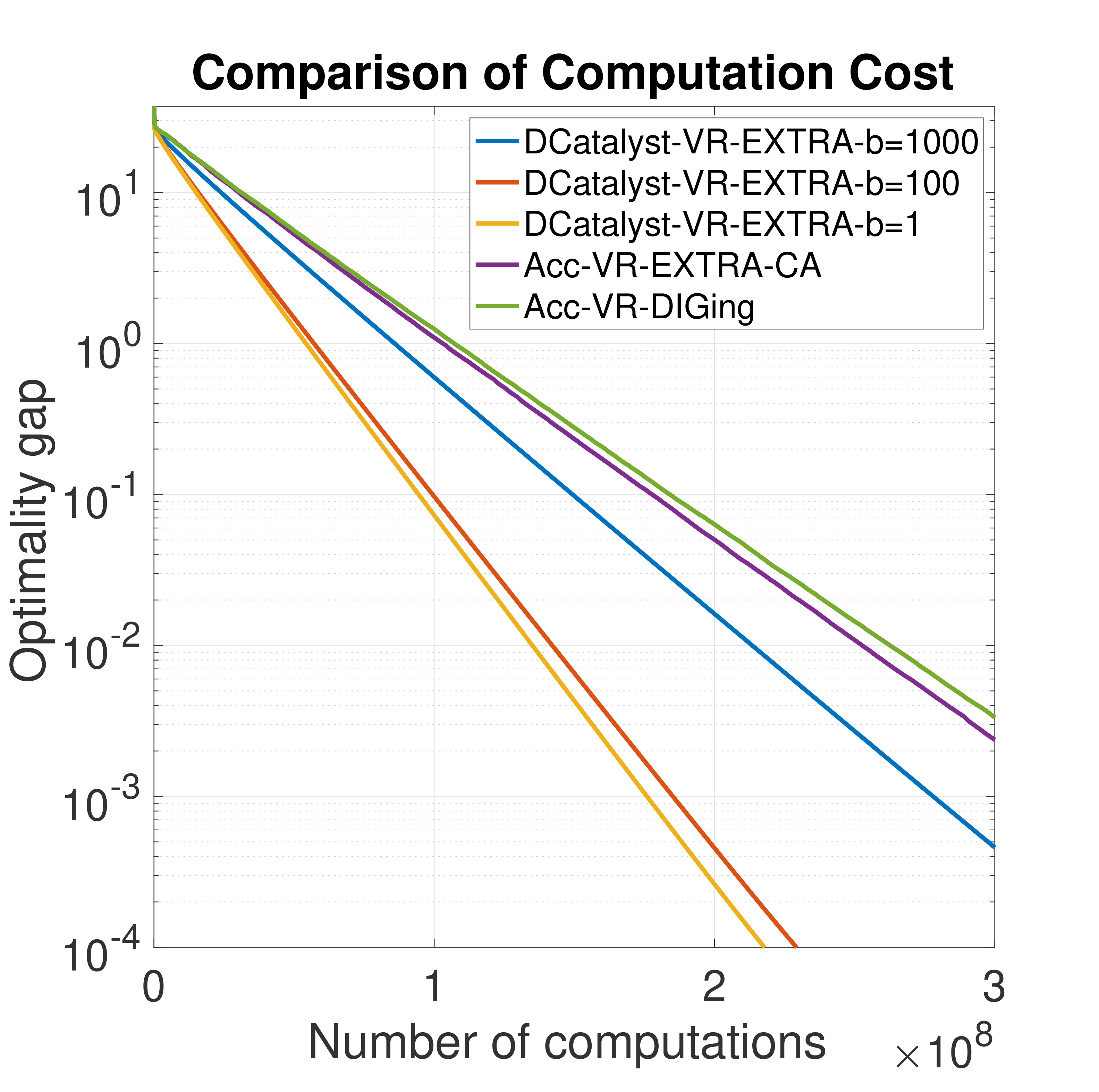}
\end{minipage}
\caption{  Comparison of distributed algorithms under finite-sum setting with a strongly convex smooth objective function in \textbf{(a):} communication cost (left) and \textbf{(b):} Computation cost (right).}
\label{'6_3_2'} 
\end{figure}

 Figure~\ref{'6_3_2'}  presents the main results. It plots the optimality gap achieved by the aforementioned algorithms as a function of the number of communication rounds (left panel) and the number of iterations (right panel), where iterations are measured as local gradient $\nabla f_{ij}$ evaluations.
 The results indicate that DCatalyst-VR-EXTRA surpasses both Acc-VR-EXTRA-CA and Acc-VR-DIGing-CA in terms of both computation and communication costs. Additionally, the figure highlights an interesting trade-off in DCatalyst-VR-EXTRA when varying the batch size   $b$: as $b$ grows, the communication cost reduces and the computation cost increases, increases, the communication cost decreases while the computation cost increases, demonstrating the efficiency gains from adjusting batch sizes within our DCatalyst framework.

\section*{Acknowledgements}
This work has been supported by   the Office of Naval Research (ONR Grant N. N000142412751).

\appendix 

\renewcommand{\thesection}{Appendix \Alph{section}}
\renewcommand{\thesubsection}{\Alph{section}.\arabic{subsection}}
\section{}
\label{appendix:a}

\subsection{Proof of Lemma~\ref{lemma:prelim}}\label{app:proof_lemma_preliminary} 
Since $\Mud(\tx_i^k) \leq \psi_i^{k,\star}+\epsilon_{\psi,i}^k$, we have
$
 \psi_i^{k,\star}+\ep_{\psi,i}^{k}-\Muds \geq 0.
$
Then
\begin{equation*}
    \begin{aligned}
     0&\leq    \ep_{\psi,i}^{k}+\psi_i^{k,\star}-\Muds \leq\ep_{\psi,i}^{k}+\psi_i^{k}(x^{\star})-\Muds\\    
          &  \leq  (1-\alpha^{k-1})(\psi_i^{k-1}(x^{\star} )+\ep_{\psi,i}^{k-1} -\Muds)+\ep_{\psi,i}^{k} +\alpha^{k-1}\ep_{i}^{k-1} -(1-\alpha^{k-1})\ep_{\psi,i}^{k-1},
    \end{aligned}
\end{equation*}
where the last inequality is due to (ii) of Definition~\ref{def:inexact-ES}.
Dividing  by $\lambda^{k}$ both sides of the last inequality, yields  
\begin{equation*}
  \!\!\!\!\!  \frac{\psi_i^{k}(x^{\star} )+\ep_{\psi,i}^{k} -\Muds}{\lambda^{k}}\leq \frac{\psi_i^{k-1}(x^{\star} )+\ep_{\psi,i}^{k-1} -\Muds}{\lambda^{k-1}}+\frac{\ep_{\psi,i}^{k} -(1-\alpha^{k-1})\ep_{\psi,i}^{k-1} +\alpha^{k-1}\ep_{i}^{k-1 } }{\lambda^{k}}. \qquad \hfill \square
\end{equation*}

\subsection{Proof of Proposition~\ref{prop:preliminary}}\label{app:proof_Cor_preliminary}
 
It follows from \eqref{eq:telescope-error-M}, 
$ \psi_i^{k}(x^{\star} )+\ep_{\psi,i}^{k} -\Muds  \geq \psi_i^{k,\star} +\ep_{\psi,i}^{k} -\Muds \geq \Mud(\tx_i^{k} )-\Muds.$
    We have  
\begin{equation}\label{eq:upper-bound-M-tele}
     \Mud(\tx_i^{k} )-\Muds  \leq \psi_i^{k}(x^{\star}  )+\ep_{\psi,i}^{k} -\Muds \leq \lambda^{k}\left( \psi_i^{0}(x^{\star})-\Muds+\sum_{j=0}^{k-1}\frac{\epsilon_{tot,i}^{j} }{\lambda^{j+1}}\right),
\end{equation} We separate the discussion into the two cases $\mu > 0$ and $\mu = 0$.

\noindent   {\bf (i) $\mu > 0$:}  Summing \eqref{eq:upper-bound-M-tele} over $i$ while taking the expected value, and using  
    $\lambda^k = (1-\alpha)^k$, we have 
    \begin{equation*}
    \begin{aligned}
          &\overm\summ\lambda^{k}\left(\psi_i^{0}(x^{\star})-\Muds+\sum_{j=0}^{k-1}\frac{\Eb[\ep_{tot,i}^{j}]}{\lambda^{j+1}}\right) \\
          &\leq(1-c\alpha)^{k} \overm\summ\left(\psi_i^{0}(x^{\star})-\Muds+\frac{C_{scvx}}{1-c\alpha}\sum_{j=0}^{k-1}\frac{(1-\alpha)^{k-1-j}}{(1-c\alpha)^{k-1-j}}\right)\\
          &\leq (1-c\alpha)^{k}\overm\summ\left(\psi_i^{0}(x^{\star})-\Muds+\frac{C_{scvx}}{\alpha(1-c)}\right). 
    \end{aligned}
\end{equation*}
Invoking the    $\muM$-strong convexity of $\Mud$,  we deduce
\begin{equation*}
    \begin{aligned}
             \overm\Eb[\|\tilde{\vx}^{k}-\vx^{\star}\|^{2}]&\leq \frac{2}{m\muM}\summ\left(\Eb[\Mud(\tx_i^{k})]-\Muds\right) \leq (1-c\alpha)^{k}\overm\summ\left(\frac{2}{\muM}\left(\psi_i^{0}(x^{\star})-\Muds+\frac{C_{scvx}}{\alpha(1-c)}\right)\right)\\
             &= c_{scvx}\left(1-c\alpha\right)^{k}. 
    \end{aligned}
\end{equation*}
 
  \noindent   {\bf (ii) $\mu = 0$:}  We hinge on the   lemma below,  which follows from \cite{lin2015universal} applied to our setting.

\begin{lemma}\label{lemma:lambdak}
 Let $\{\lambda^{k}\}_{k }$ be defined in \eqref{def_lambdak}  where $\{\alpha^{k}\}_{k }$ is defined as in Algorithm~\ref{alg:framework} for the case $\mu = 0$. Then, for all $k\geq0$,  
$$
     \frac{4}{(k+2)^{2}} \geq \lambda^{k} \geq \frac{2}{(k+2)^{2}}.
$$
 \end{lemma}
 Using Lemma~\ref{lemma:lambdak}, we have  
\begin{equation*}
    \begin{aligned} &\overm\summ\left(\Eb[\Mud(\tx_i^{k})]-\Muds\right) \leq  \lambda^{k}\overm\summ\left(\psi_i^{0}(x^{\star})-\Muds+\sum_{j=0}^{k-1}\frac{\Eb[\ep_{tot,i}^{j}]}{\lambda^{j+1}}\right) \\
    &\leq \lambda^{k}\overm\summ\left(\psi_i^{0}(x^{\star})-\Muds+\sum_{j=0}^{\infty} \frac{C_{cvx}(j+3)^{2}}{2(j+1)^{3+r_{0}}}\right) \\
          &\leq \lambda^{k}\overm\summ\left(\psi_i^{0}(x^{\star})-\Muds+\sum_{j=1}^{\infty} \frac{C_{cvx}((j+1)^{2}+4)}{(j+1)^{3+r_{0}}}+\frac{9C_{cvx}}{2}\right) \\
    &\leq \frac{4}{(k+2)^{2}}\overm\summ\left(\psi_i^{0}(x^{\star})-\Muds+\sum_{j=1}^{\infty} \frac{C_{cvx}((j+1)^{2}+4)}{(j+1)^{3+r_{0}}}++\frac{9C_{cvx}}{2}\right) \\ 
    & \leq   \frac{4}{(k+2)^{2}}\overm\summ\left(\psi_i^{0}(x^{\star})-\Muds+\frac{10C_{cvx}}{r_{0}} \right).
    \end{aligned}
\end{equation*}
The desired result \eqref{eq:estimate-seq-sublinear-M} follows from the above inequality and the fact 
$$  
    \Mud(\tilde{x}_i^k)-\Muds \geq \frac{1}{2\delta}\|\nabla\Mud(\tx_i^{k})\|^{2}.
 $$  
If $r\equiv 0$,   $u(x)$ is $L$-smooth. Define $\tx_i^{k,\star} := \tx_i^{k} -(1/\delta)\nabla\Mud(\tx_i^{k} )$. Then
\begin{equation*}
\begin{aligned}
    \|\nabla\Mud(\tx_i^{k} )\| &= \|\nabla u(\tx_i^{k,\star} )\|  
     \geq \|\nabla u(\tx_i^{k} )\|-\frac{L}{\delta}\|\nabla \Mud(\tx_i^{k} )\|.
\end{aligned}
\end{equation*}
Using this lower bound in \eqref{eq:estimate-seq-sublinear-M} yields \eqref{eq:estimate-seq-sublinear-u}. \hfill $\square$

\renewcommand{\thesection}{Appendix \Alph{section}}
\renewcommand{\thesubsection}{\Alph{section}.\arabic{subsection}}
\section{}
\label{appendix:b}

\subsection{Proof of Lemma~\ref{application:sonata-assump5}}\label{app:proof-lemma_app_sonata}

  For any variables $ x_1,  x_{0} \in \mathbb{R}^d$,  
 $$
\begin{aligned}
&u^{k+1}(x_1)-u^{k}(x_1)-(u^{k+1}(x_0)-u^{k}(x_0))\\
& = \frac{\delta}{2m} \summ\left(\|x_1-z_{i}^{k+1}\|^2-\|x_1-z_{i}^{k}\|^2-\|x_0-z_{i}^{k+1}\|^{2}+\|x_0-z_{i}^{k}\|^{2}\right).
\end{aligned}    
 $$
Hence, 
 $$
\begin{aligned}
&u^{k+1}(x_{i}^{k+1})-u^{k+1}(x^{k+1,\star})\\ 
& = u^{k}(x_{i}^{k+1})-u^{k}(x^{k,\star})+u^{k}(x^{k,\star})-u^{k}(x^{k+1,\star})+\delta\langle\bar{z}^{k}-\bar{z}^{k+1},x_{i}^{k+1}-x^{k+1,\star}\rangle\\
&\leq u^{k}(x_{i}^{k+1})-u^{k}(x^{k,\star})-\frac{\delta+\mu}{2}\|x^{k,\star}-x^{k+1,\star}\|^2+\delta\langle\bar{z}^{k}-\bar{z}^{k+1},x_{i}^{k+1}-x^{k+1,\star}\rangle.
\end{aligned} 
 $$
Then, 
\begin{equation}\label{ineq1:proof-sonata-assum4}
    \begin{aligned}
    &\frac{2}{(\mu+\delta) m}\summ(u^{k+1}(x_{i}^{k+1})-u^{k+1}(x^{k+1,\star}))\\
    &\leq \frac{2}{(\mu+\delta)m}\summ (u^k(x_i^{k+1})-u^{k}(x^{k,*}) )-\|x^{k,\star}-x^{k+1,\star}\|^2
    +\frac{2\delta}{\mu+\delta}\langle\bar{z}^{k}-\bar{z}^{k+1}, \bar{x}^{k+1}-x^{k+1,\star}\rangle.
    \end{aligned}
\end{equation}
We bound the last term on the RHS of the inequality as follows:  
\begin{equation}\label{ineq2:proof-sonata-assum4}
    \begin{aligned}
    &\frac{2\delta}{\mu+\delta}\langle\bar{z}^k-\bar{z}^{k+1},x^{k,\star}-x^{k+1,\star}\rangle \leq  \|x^{k,\star}-x^{k+1,\star}\|^2+\frac{\delta^2}{(\delta+\mu)^2}\|\bar{z}^{k}-\bar{z}^{k+1}\|^2,\\
    &\frac{2\delta}{\mu+\delta}\langle\bar{z}^k-\bar{z}^{k+1}, \bar{x}^{k+1}-x^{k,\star}\rangle \leq  \|x^{k,\star}-\bar{x}^{k+1}\|^2+\frac{\delta^2}{ (\delta+\mu)^2}\|\bar{z}^{k}-\bar{z}^{k+1}\|^2.
    \end{aligned}
\end{equation}
Combine  \eqref{ineq1:proof-sonata-assum4} and  \eqref{ineq2:proof-sonata-assum4},
\begin{equation}\label{part1}
    \begin{aligned}
    &\frac{2}{(\mu+\delta)m}\summ (u^{k+1}(x_{i}^{k+1})-u^{k+1}(x^{k+1,\star}) )\\
    &\leq \frac{2}{(\mu+\delta)m}\summ (u^k(x_i^{k+1})-u^k(x^{k,*}) )+\frac{2\delta^2}{(\delta+\mu)^2} \|\bar{z}^{k}-\bar{z}^{k+1}\|^2+\|\bar{x}^{k+1}-x^{k,\star}\|^2\\
    &\leq  \frac{4}{(\mu+\delta)m}\summ (u^k(x_i^{k+1})-u^k(x^{k,*}) )+\frac{2\delta^2}{(\delta+\mu)^2}\|\bar{z}^k-\bar{z}^{k+1}\|^2.
    \end{aligned}
\end{equation}

Using the  initialization \eqref{cata-sonata:ini},
\begin{equation}\label{Consensus}
\begin{aligned}
& \frac{\eta}{m} \left(4 (L_{\max}+\delta)^2 \|\vx_{\bot}^{k+1}\|^2 + 2\|\vy_{\bot}^{k+1}\|^2\right)\\
& =\frac{\eta}{m} \left(4 (L_{\max}+\delta)^2 \|\vx_{\bot}^{k,T_k}\|^2 + 2\|\vy_{\bot}^{k,T_k}+\delta(I-\overm\mathbf{1}_{m}\mathbf{1}_{m}^{T})(\vz^{k}-\vz^{k+1})\|^2\right) \\
& \leq \frac{\eta}{m}\left(4(L_{\max}+\delta)^2\|\vx_{\bot}^{k,T_k}\|^2 + 4  \|\vy_{\bot}^{k,T_k}\|^2+16\delta^{2}\|\vz^{k+1}-\vz^{k}\|^2\right).\\ 
\end{aligned}
\end{equation}

Combine \eqref{part1} and \eqref{Consensus}, yields  
\begin{equation*} 
    \begin{aligned}
     \cL^{k+1}(\vs^{k+1,0})  
    &\leq 2\cL^{k}(\vs^{k,T_k})+\frac{2\delta^2}{(\delta+\mu)^2}\|\bar{z}^{k+1}-\bar{z}^{k}\|^2+\frac{16\eta\delta^{2}}{m}\|\vz^{k+1}-\vz^{k}\|^2\\
    &\leq 2L^{k}(\vs^{k,T_k})+\left(\frac{2\delta^{2}}{m(\delta+\mu)^2}+\frac{16 \eta\delta^{2}}{m}\right)\|\vz^{k+1}-\vz^{k}\|^2.\qquad \square
    \end{aligned}
\end{equation*}

\subsection{Proof of  Corollary~\ref{coro:sonata-scvx}}\label{app:proof:coro:sonata-scvx}

By Lemma~\ref{application:sonata-assump4} and~\ref{application:sonata-assump5}, SONATA algorithm (under the initialization \eqref{cata-sonata:ini})  satisfies Assumptions~\ref{Lyapunov} and~\ref{assump:warm-start}. 

We provide next a suitable value for  $T_k$. It is trivial that check that  \eqref{application:sonata-r-cM}  holds and $r_{\cM,\delta} = \mathcal{O}(1)$. Let us  choose $\epsilon_0 = \max\left\{(1/m)\|\vx^{\star}-\vx^0\|^2, \cL^{0}(\vs^{0,0})\right\}.$ 
 
By \eqref{cal:c_scvx}, it follows 
$$
\begin{aligned}
    & c_{scvx} \\
    &\leq 
     \overm\left(2+\frac{\delta}{\mu}\right)\|\vx^{\star}-\vx^0\|^2+\frac{2(\mu+\delta)^2}{\mu^2(1-c)}\|\vx^{\star}-\vx^0\|^2+ \left(\frac{\mu+\delta}{\mu(1-c)}+\frac{2\sqrt{2000}(\mu+\delta)^2\sqrt{m}}{\mu^2(1-c)^2}\right)\epsilon_0\\
     & = \left(2+\frac{\delta}{\mu}+\frac{(2m+1)(\mu+\delta)^2}{\mu^2(1-c)}+  \frac{2\sqrt{2000}(\mu+\delta)^2\sqrt{m}}{\mu^2(1-c)^2} \right)\epsilon_0,
\end{aligned}
 $$
and, according to  Proposition~\ref{prop:inner-loop-scvx},  \eqref{ineq:prop-out-scvx-Tk}  holds under the following bound on $T_k$: 
$$
\begin{aligned}
    & r_{\cM,\delta}\log\frac{  c_{\cL}  \,c_{\cM}  + 36 d_{\cM}\,c_{scvx}\dfrac{1}{\epsilon_{0}(1-c\,\alpha)^{2}} }{1-c\alpha}   = r_{\cM,\delta}\log\frac{ 2  + 36\left( \frac{2\delta^{2}}{(\delta+\mu)^2}+16 \eta\delta^{2}\right)\,\dfrac{c_{scvx}}{\epsilon_{0}(1-c\,\alpha)^{2}} }{1-c\alpha} \\
    & \leq   r_{\cM,\delta}\log\frac{ 2  + \frac{36}{(1-c\,\alpha)^{2}}\left( \frac{2\delta^{2}}{(\delta+\mu)^2}+16 \eta\delta^{2}\right)\left(2+\frac{\delta}{\mu}+\frac{(2m+1)(\mu+\delta)^2}{\mu^2(1-c)}+  \frac{2\sqrt{2000}(\mu+\delta)^2\sqrt{m}}{\mu^2(1-c)^2}  \right)}{1-c\alpha}. 
\end{aligned}
$$

We can now   apply Theorem~\ref{theorem:total-iter-scvx}. The total numbers of inner-plus-outer iterations $N_{\texttt{it}}$ of DCatalyst-SONATA-F and DCatalyst-SONATA-L are respectively
$$
 \widetilde{\mathcal{O}}\left(\sqrt{\frac{\beta}{\mu}}\log\frac{1}{\epsilon}\right)\quad\text{and}\quad \widetilde{\mathcal{O}}\left(\sqrt{\kappa_g}\log\frac{1}{\epsilon}\right).   
$$
When solving subproblem~\eqref{eq:subproblem_DCat}, the network connectivity $\rho$ is required to  satisfy (Lemma~\ref{application:sonata-assump4}) 
$$
\rho = \left\{\begin{aligned}
 &   \mathcal{O}\left(\left(\frac{\beta}{L+\beta-\mu}+1\right)^{2}\right), &(\text{DCatalyst-SONATA-F});\\
    &   \mathcal{O}\left(\frac{L}{L_{\max}+L-\mu}\right), & (\text{DCatalyst-SONATA-L}).
\end{aligned}\right. 
$$ 
This conditons can be eforced by running multiple communication rounds  per algorithm iteration.   Specifically, employing Chebyshev acceleration \cite{scaman2017optimal} yields  $N_{\texttt{com}, \cM} = \widetilde{\mathcal{O}}(\sqrt{{1}/({1-\rho})})$. 
 Then, the total number of communications $N_{\texttt{com}}$ for DCatalyst-SONATA-F and DCatalyst-SONATA-L reads  respectively
 $$\widetilde{\mathcal{O}}\left(\sqrt{\frac{\beta}{\mu(1-\rho)}}\log\frac{1}{\epsilon}\right)    
\quad\text{and}\quad       \widetilde{\mathcal{O}}\left(\sqrt{\frac{\kappa_g}{1-\rho}}\log\frac{1}{\epsilon}\right).
    $$
 
For the sake of completeness, we provide also the computational complexities of DCatalyst-SONATA-F (the one of DCatalyst-SONATA-L coincides with $N_{\texttt{it}}$).     
Suppose  that, at every iteration     $t$ of SONATA-F,  each agent solves its own subproblem inexactly  
\begin{equation*}
   x_i^{k,t+1} = \arg\min_{x\in\mathbb{R}^d} \underbrace{f_i(x)+\frac{\beta-\mu}{2}\|x-z_i^k\|^2+\frac{\beta}{2}\|x-x_i^{k,t}\|^2+r(x)}_{u_i^{k,t}(x)},
\end{equation*}
employing the  accelerated proximal gradient method, such that the resulting inexact solution $\hat{x}_i^{k,t+1}$ satisfies
$$\max\left\{\frac{1}{2\beta}( u_i^{k,t}(\hat{x}_i^{k,t+1})-u_i^{k,t}(x_i^{k,t+1})), \|\hat{x}_i^{k,t+1}-x_i^{k,t+1}\|^2\right\} = 
\mathcal{O}\left(\epsilon_{0}(1-c\alpha)^{k+1}\left(1-\frac{1}{r_{\cM,\delta}}\right)\right).$$   
Following similar steps as those in  \cite[Appendix E]{similarity2022}, the
 number of computations for DCatalyst-SONATA-F  to meet such an inexact termination is
$$   
\widetilde{\mathcal{O}}\left(\sqrt{\frac{L+2\beta-\mu}{\beta-\mu}}\cdot\frac{\beta}{\mu}\log^2\frac{1}{\epsilon}\right) = \widetilde{\mathcal{O}}\left(\sqrt{\frac{L+\beta}{\beta}}\cdot\frac{\beta}{\mu}\log^2\frac{1}{\epsilon}\right) \qquad \square. $$

\subsection{Proof of  Corollary~\ref{coro:sonata-cvx}}\label{app:proof:coro:sonata-cvx}
Lemma~\ref{application:sonata-assump4} and  Lemma~\ref{application:sonata-assump5}  hold.  
To found a suitable bound of  $T_k$, we let $\epsilon_0 := \max\{B^2, \cL^{0}(\vs^{0,0})\}$, then 
$$
\begin{aligned}
    & r_{\cM,\delta}\log \left(c_{\cL} c_{\cM}2^{4+2r_{0}} +\frac{36c_{\cL} d_{\cM}B^2(k+3)^{4+2r_{0}}}{\epsilon_{0}} \right)  
      \leq r_{\cM,\delta}\log \left(2^{5+2r_0}+1224(k+3)^{4+2r_0}\right) .
\end{aligned}
$$
According to Proposition~\ref{prop:inner-loop-cvx},   \eqref{ineq:prop-out-cvx-Tk}  holds if $T_k \geq r_{\cM,\delta}\log \left(2^{5+2r_0}+1224(k+3)^{4+2r_0}\right)$.

We can  apply Theorem~\ref{theorem:total-iter-cvx} and, following  similar steps as  in Sec.~\ref{app:proof:coro:sonata-scvx},   we obtain the communication and computational complexities for DCatalyst-SONATA-F and DCatalyst-SONATA-L as stated in the corollary.
\hfill$\square$

\subsection{Proof of Lemma~\ref{application-puda:assum5}}\label{app:proof-puda-assum5} 
When applying PUDA to solve the subproblem~\eqref{eq:subproblem_DCat}, denote the fixed point of $\{y_i^{k,t}\}_{t\geq 0}$ as $y^{k,\star}$ and $\vy^{k,\star} := 1_m(y^{k,\star})^{\top}$. 
We have
    \begin{equation} \label{ineq:puda-assum5-ineq1}
        \begin{aligned} &\|\vx^{k,T_k}-\vx^{k+1,\star}\|^2 \leq 2\|\vx^{k,T_k}-\vx^{k,\star}\|^2 + 2\|\vx^{k+1,\star}-\vx^{k,\star}\|^2 \leq 2\|\vx^{k,T_k}-\vx^{k,\star}\|^2 + 2\|\vz^{k+1}-\vz^{k}\|^2, \\ 
            & \|\vy^{k,T_k}-\vy^{k+1,\star}\|^2 \leq 2\|\vy^{k,T_k}-\vy^{k,\star}\|^2+2\|\vy^{k,\star}-\vy^{k+1,\star}\|^2. 
        \end{aligned}
    \end{equation}
    Given $k\geq 0$, according to   \cite[Lemma 1 ]{alghunaim2020decentralized}, the following holds  
     \begin{equation*}
        H\vx^{k,\star} - \eta  H\nabla F^{k}(\vx^{k,\star})-H^2 \vy^{k,\star} = 0,
    \end{equation*}
    where
    $\nabla F^k(\vx^{k,\star})$ denotes $[\nabla f_1^{k}(x^{k,\star}), \cdots, \nabla f_m^k(x^{k,\star})]^{\top}$. 
    
    Then
    \begin{equation}\label{ineq:puda-assum5-ineq2}
        \begin{aligned}
            &\sigma_{\min}^{+}(H^2)\|\vy^{k+1,\star}-\vy^{k,\star}\|^2  \leq \| H\vx^{k,\star} - \eta  H\nabla F^{k}(\vx^{k,\star})- H\vx^{k+1,\star} + \eta  H\nabla F^{k+1}(\vx^{k+1,\star})\|^2 \\ 
            & \leq \left((1+\delta\eta )\|H(\vx^{k,\star}-\vx^{k+1,\star})\|+\eta \|H(\nabla F(\vx^{k+1,\star})-\nabla F(\vx^{k,\star}))\|+\delta\eta \|H(\vz^{k}-\vz^{k+1})\|\right)^2\\  
               & \leq 
            \sigma_{\max}(H^2)\left(9+9\delta^2\eta^2     \right)\|\vz^{k}-\vz^{k+1}\|^2.\\ 
        \end{aligned}
    \end{equation}
   Combine  \eqref{ineq:puda-assum5-ineq1}  and  \eqref{ineq:puda-assum5-ineq2}, we obtain the desired expression of $c_{\cM}$ and $d_{\cM}$. 

 \hfill$\square$

\subsection{Proof of Corollary~\ref{coro:puda-scvx}}\label{app:proof:coro:puda-scvx}

By Lemma~\ref{application-puda:assum4} and~\ref{application-puda:assum5}, Assumption~\ref{Lyapunov} and~\ref{assump:warm-start} hold when using  PUDA. 
With  $\delta$ chosen as in  \eqref{cata-puda:ini}, the convergence rate reads
$$
r_{\cM, \delta} = \max\left\{\frac{1}{\sigma_{\min}^{+}(H^2)}, \frac{4(L_{\max}+\delta)}{\left(2-\sigma_{\max}(C)\right)^2(\mu_{\min}+\delta)}\right\} = \frac{1}{\sigma_{\min}^{+}(H^2)}.
$$
We can now provide a valid bound on   $T_k$. Similar to the proof of Corollary~\ref{coro:sonata-scvx}, we pick $\epsilon_0 = \max\{(1/m)\|\vx^{\star}-\vx^0\|^2, \cL^{0}(\vs^{0,0})\}$, yielding
$$ c_{scvx} \leq \left(2+\frac{\delta}{\mu_{\min}}+\frac{(2m+1)(\mu_{\min}+\delta)^2}{\mu_{\min}^2(1-c)}+  \frac{2\sqrt{2000}(\mu_{\min}+\delta)^2\sqrt{m}}{\mu_{\min}^2(1-c)^2} \right)\epsilon_0,
$$
 and thus the following bound can be used as   valid value of $T_k$: 
$$
  \begin{aligned} &r_{\cM,\delta}\log\frac{  c_{\cL}  \,c_{\cM}  + 36 d_{\cM}\,c_{scvx}\dfrac{1}{\epsilon_{0}(1-c\,\alpha)^{2}} }{1-c\alpha} \\ 
&\leq r_{\cM,\delta}\log\frac{ 2\! + \!\frac{36}{(1-c\,\alpha)^{2}}\!\left(\!2\!+\!\frac{    \sigma_{\max}(H^2)\left(9+9\delta^2\eta^2    \right)\!}{\sigma_{\min}^{+}(H^2)}\right)\!\!\left(2\!+\!\frac{\delta}{\mu_{\min}}\!+\!\frac{(2m+1)(\mu_{\min}+\delta)^2}{\mu_{\min}^2(1-c)}\!+\!  \frac{2\sqrt{2000}(\mu_{\min}+\delta)^2\sqrt{m}}{\mu_{\min}^2(1-c)^2}  \!\right)}{1-c\alpha}. 
\end{aligned}
$$
 
Under the above setting, we can call   Theorem~\ref{theorem:total-iter-scvx} and obtain the total number $N_{\texttt{it}}$ of inner-plus-outer iterations of PUDA as given in Corollary~\ref{coro:puda-scvx}.  
According to \eqref{application:puda-update},  
both   communication and computational complexities are of the order of   $N_{\texttt{it}}$.
\hfill$\square$

\subsection{Proof of Lemma~\ref{application:pmgt-assum5}}\label{proof-pmgt-assum5} 

Denote the smoothness parameter for $f^k$ of (S.1) in Algorithm~\ref{alg:framework} as    $L_{\delta}$,  $p_{ij}$ in  \eqref{application:pmgt-update}  becomes $ {L_{\delta,ij} }/({\sum_{j=1}^n L_{\delta,ij}})$. We have   

\begin{equation}\label{pmgt-assum5-bx}
    \begin{aligned} 
        &\|\bx^{k,T_k}-x^{k+1,\star}\|^2 
         \leq 2\|\bx^{k,T_k}-x^{k,\star}\|^2+\frac{2}{m}\|\vz^k-\vz^{k+1}\|^2, \\
        &\frac{1}{m}\|\vx_{\bot}^{k+1}\|^2+\frac{\eta^2}{m}\|\vy_{\bot}^{k+1} \|^2   
          \leq \overm\left(\|\vx_{\bot}^{k,T_k}\|^2+2\eta^2\|\vy_{\bot}^{k,T_k}\|^2+8\eta^2\delta^2\|\vz^{k}-\vz^{k+1}\|^2\right),
    \end{aligned}
\end{equation} and 
\begin{equation}\label{pmgt-assum5-Deltaf}
    \begin{aligned}
         \Delta_{f}^{k+1,0} 
        &  \leq \frac{1}{mn}\summ\sum_{j=1}^{n}\frac{1}{np_{ij}}\left(2\|\nabla f^{k+1}_{ij}(v_i^{k,T_k})-\nabla f^{k+1}_{ij}(x^{k,\star})\|^2+ 2\|\nabla f^{k+1}_{ij}(x^{k+1,\star})-\nabla f^{k+1}_{ij}(x^{k,\star})\|^2\right) \\
        & \leq 2\Delta_{f}^{k,T_k}+\frac{2\bar{L}_{\delta,\max}}{(mn)^2}\summ\sum_{j=1}^n L_{\delta,ij}\|\vz^k-\vz^{k+1}\|^2,
    \end{aligned}
\end{equation}
where   $L_{\delta,ij} = L_{ij}+\delta$ and  $\bar{L}_{\delta,\max} = \bar{L}_{\max}+\delta$.
The proof of the lemma follows readily 
combining \eqref{pmgt-assum5-bx}  and \eqref{pmgt-assum5-Deltaf}. 
\hfill$\square$

\subsection{Proof of Corollary~\ref{coro:pmgt-total}
}\label{proof-pmgt-corollary}

    By Lemma~\ref{application:pmgt-assum4} and ~\ref{application:pmgt-assum5}, PMGT-LSVRG  (under the  initialization  \eqref{cata-pmgt:ini}) satisfies Assumption~\ref{Lyapunov} and~\ref{assump:warm-start}. 
    The convergence rate is
    $r_{\cM,\delta} = 48n$. 
    
   We proceed bounding $T_k$. Similarly to the procedure in the analysis of Corollary~\ref{coro:sonata-scvx}, we set $\epsilon_0 = \max\{(1/m)\norm{\vx^0-\vx^{\star}}^2, \cL^{0}(\vs^{0,0})\}$. Then,
$$
\begin{aligned}
     c_{scvx}  
      \leq \left(2+\frac{\delta}{\mu}+\frac{(2m+1)(\mu+\delta)^2}{\mu^2(1-c)}+  \frac{2\sqrt{2000}(\mu+\delta)^2\sqrt{m}}{\mu^2(1-c)^2} \right)\epsilon_0.
\end{aligned}
 $$ 
 Invoking Proposition~\ref{prop:inner-loop-scvx}, the following bound holds for $T_k$:
$$
\begin{aligned}
    &   r_{\cM,\delta}\log\frac{ 2  + 36\left(2+8\eta^2\delta^2+\frac{8\eta^2(\bar{L}_{\max}+\delta)^2}{n^2}\right)\,\dfrac{c_{scvx}}{\epsilon_{0}(1-c\,\alpha)^{2}} }{1-c\alpha} \\
    & \leq r_{\cM,\delta}\log\frac{ 2  + \frac{36}{(1-c\,\alpha)^{2}}\left(2+8\eta^2\delta^2+\frac{8\eta^2(\bar{L}_{\max}+\delta)^2}{n^2}\right)\left(2+\frac{\delta}{\mu}+\frac{(2m+1)(\mu+\delta)^2}{\mu^2(1-c)}+  \frac{2\sqrt{2000}(\mu+\delta)^2\sqrt{m}}{\mu^2(1-c)^2}  \right)}{1-c\alpha} \\
    & \leq r_{\cM,\delta}\log\frac{ 2  + \frac{90}{(1-c\,\alpha)^{2}}\left(2+\frac{\delta}{\mu}+\frac{(2m+1)(\mu+\delta)^2}{\mu^2(1-c)}+  \frac{2\sqrt{2000}(\mu+\delta)^2\sqrt{m}}{\mu^2(1-c)^2}  \right)}{1-c\alpha}. 
\end{aligned}
$$

    The total number of inner-plus-outer iterations of DCatalyst-PMGT-LSVRG is then 
    $$ N_{\texttt{it}} = \widetilde{\mathcal{O}}\left(\sqrt{\kappa_s n}\log\frac{1}{\epsilon}\right).
    $$
    For each algorithm iteration, there are $N_{\texttt{com},\cM}  = 2N_{\texttt{FM}} = \mathcal{O}\left((1/\sqrt{1-\rho})\cdot\log n\right)$ numbers of communications. Hence,   the total number of communication   is
    $$
    N_{\texttt{com}}  = N_{\texttt{it}}N_{\texttt{com},\cM} = \widetilde{\mathcal{O}}\left(\sqrt{\frac{\kappa_s n}{1-\rho}}\log\frac{1}{\epsilon}\right).
    $$

The computation of full-batch gradient $\nabla f_i$ requires computing all  $(\nabla f_{ij})_{j\in[n]}$. In addition, at each iteration, one sample $j_i$ is randomly selected and $\nabla f_{ij_i}$ is computed. Therefore,  
the total    number of computations is $2N_{\texttt{it}}$ in expectation.  \hfill
 $\square$
  
\subsection{Additional numerical results}
We complement the main section of numerical results, with additional experiments. Specifically,   we include numerical results for smooth instances of Problem~\eqref{eq:problem} ($r\equiv0$) with strongly convex or (non strongly) convex $f_i$. This offers a comprehensive comparison against most representative existing accelerated decentralized methods on smooth instances of Problem~\eqref{eq:problem} while further corroborating  our theoretical results in Sec.~\ref{sec:convergence-scvx} and Sec.~\ref{sec:convergence-cvx}.  
 
\subsubsection{\texorpdfstring{$\ell_2$}--regularized Logistic Regression} 
\label{ssububsec:strongly cvx+sm}
As a strongly convex and smooth instance of Problem~\eqref{eq:problem} ($r\equiv0$), we consider the decentralized logistic regression model with $\ell_2$ regularization, which corresponds to the formulation  
\begin{equation*}
    f_i(x)=\frac{1}{n}\sum_{j=1}^{n}\log(1+\exp(-b_{ij}\cdot\langle x,a_{ij}\rangle))+\frac{\gamma_i}{2}\|x\|_2^2,\qquad r(x)=0,
\end{equation*} 
where $a_{ij}\in \mathbb{R}^{d}$ and $b_{ij}\in\{-1,1\}$. Here, the data set  $\{(a_{ij},b_{ij})\}_{j=1}^n$ is assumed to be private and  owned only by  agent $i$. We use MNIST dataset from LIBSVM \cite{refer_libsvm} of size $N=60000$ and feature dimension $d=784$. To control the values of $\kappa_g$ and $\kappa_{\ell}$, we set  each $\gamma_i=0.05$, for $i=1,2,\cdots 29$ and $\gamma_{30}=0.005$. Under this setting, we have $\kappa_g\approx 10.8494<<\kappa_l\approx 2318.6468$ and $\beta\approx 1.7035<L\approx 10.4896$.

We compare the proposed  DCatalyst-SONATA-L and DCatalyst-SONATA-F with  the representative existing accelerated decentralized algorithms, namely:  Mudag~\cite{ye2023multi}, and OPAPC~\cite{kovalev2020optimal}. All algorithms were implemented according to their theoretical guidelines~\cite{ye2023multi,kovalev2020optimal}. For DCatalyst-SONATA-L (resp.  DCatalyst-SONATA-F), we set $\delta=L-\mu$ (resp. $\delta=\beta-\mu$), the momentum parameter $\alpha=\mu/(\mu+\delta)$ (resp.  $\alpha=\mu/(\mu+\delta)$), and the number of inner-loop iterations $T_k=\lceil\log(L/\mu)\rceil$ (resp. $T_k=\lceil\log(\beta/\mu)\rceil$), for all $k$. In  DCatalyst-SONATA-L, the local agents' prox-updates are computed in closed-form while in   DCatalyst-SONATA-F,
the  accelerated proximal gradient method is employed, rub  up  to a tolerance of   $10^{-8}$.

\begin{figure}[htbp]
\centering
\begin{minipage}[t]{0.48\textwidth}
\centering
\includegraphics[width=6cm]{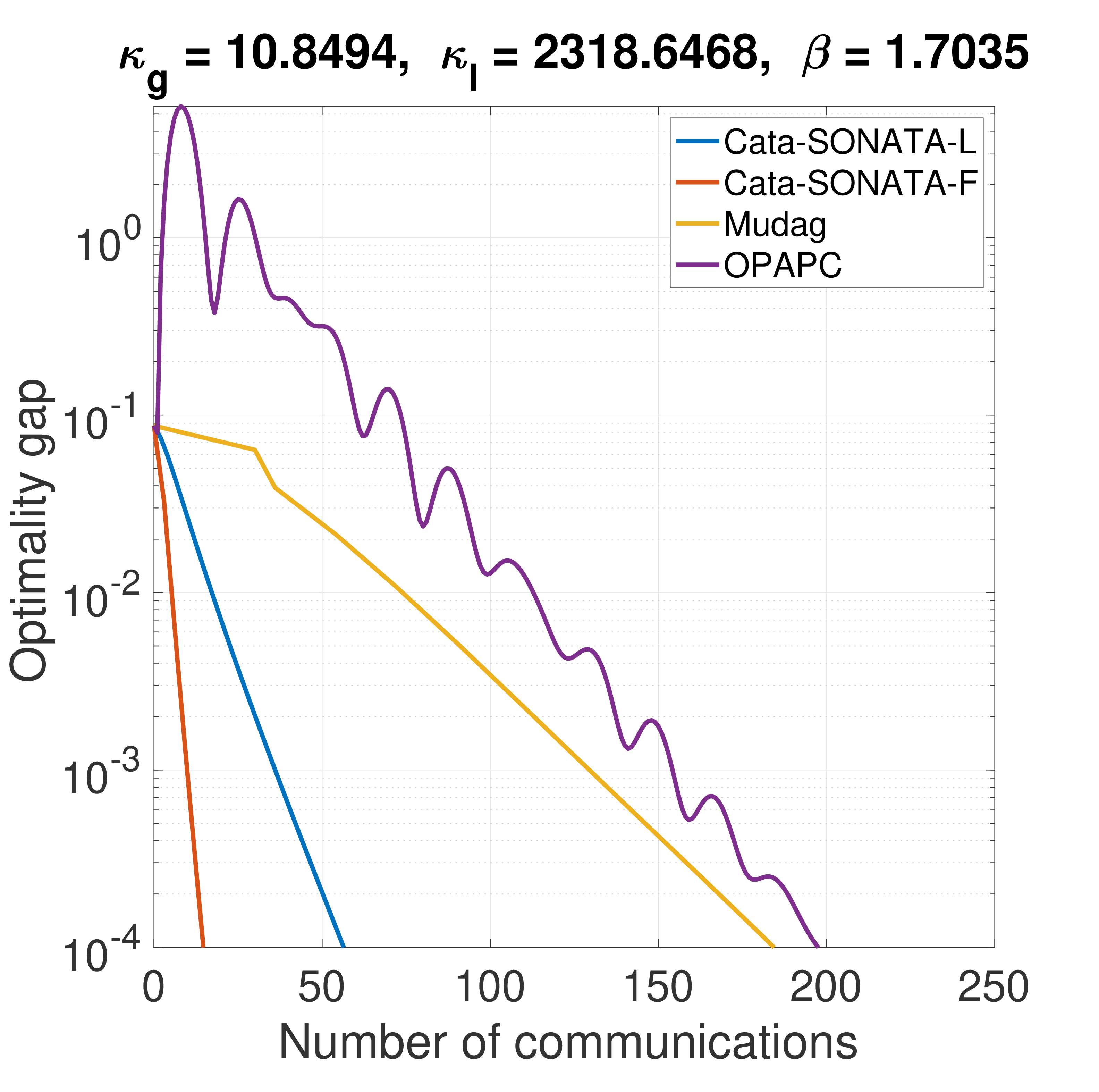}
\end{minipage}
\begin{minipage}[t]{0.48\textwidth}
\centering
\includegraphics[width=6cm]{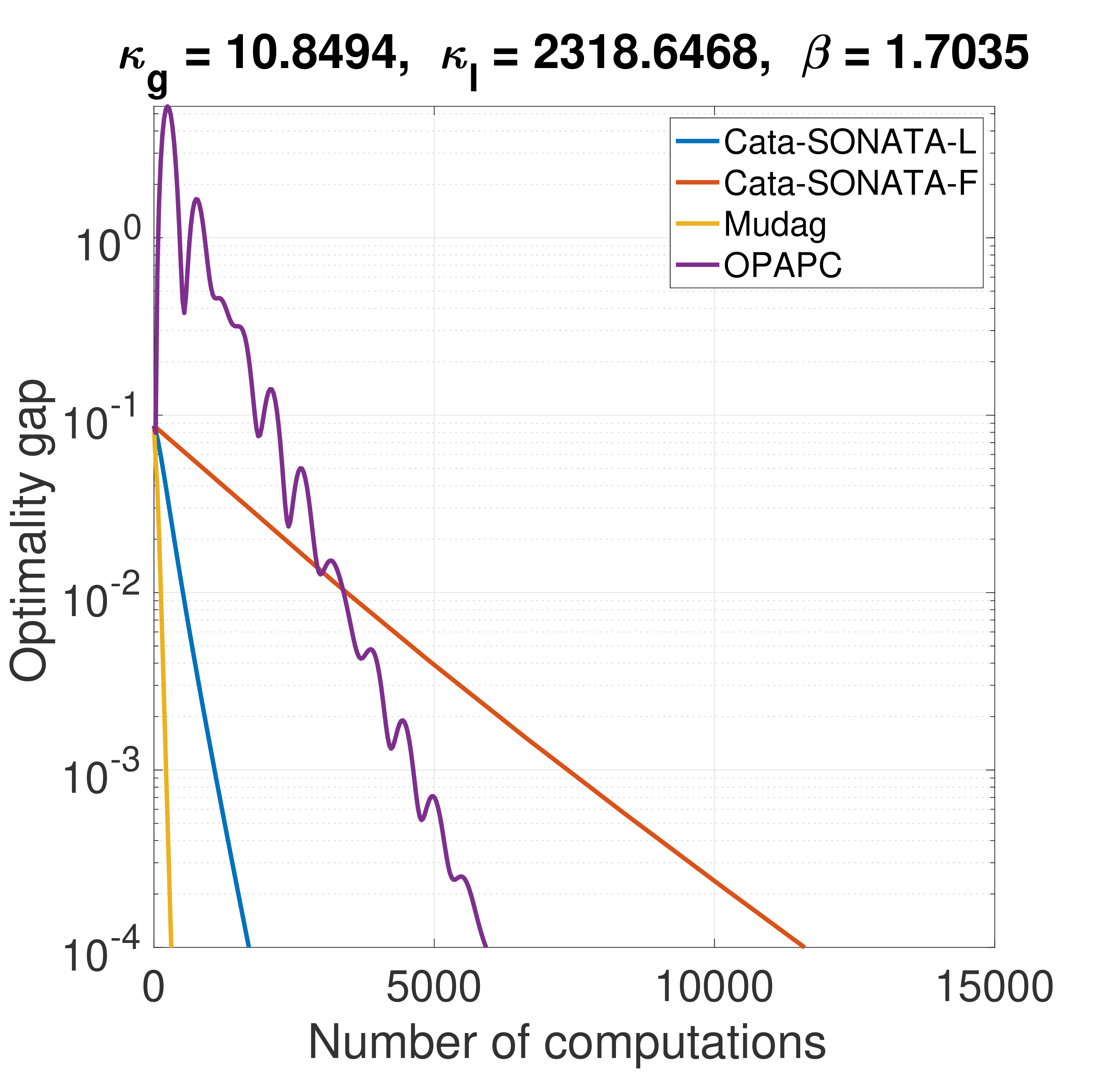}
\end{minipage}
    \caption{Comparison of distributed algorithms under strongly convex and smooth setting in \textbf{(a):} communication cost (left) and \textbf{(b):} computation cost (right).}
    \label{'6_1_1'}
\end{figure}

Figure \ref{'6_1_1'} summarizes the comparison, plotting the optimality gap  versus the number of total communications (left panel) and the number of inner-plus-outer iterations (right panel). For DCatalyst-SONATA-F, this include also the number of iterations run by the   accelerated proximal gradient method employed locally to compute the local agents' prox-updates.   The figure confirms linear convergence of the algorithms.    DCatalyst-SONATA-F is the only accelerated algorithm exploiting function similarity, yielding   much less communications, at the cost of more (albeit limited) local computations.   Notably, both DCatalyst-SONATA-L and DCatalyst-SONATA-F outperform OPAPC in communication and computation efficiency. This is because their convergence depends on the  global condition number  $\kappa_g$ rather than the local one   $\kappa_l$.  DCatalyst-SONATA-L   outperforms   Mudag   in term of communications but requires more local computations. This is compatible with our theoretical findings, predicting  an extra poly-log factor in the convergence of  DCatalyst-SONATA-L with respect to the complexity of Mudag.   

\subsubsection{  Linear Regression (\texorpdfstring{$N<d$}-) with Huber loss regularization}
\label{subsubsec:cvx+sm}
\begin{figure}[htbp] 
    \centering
\begin{minipage}[t]{0.48\textwidth}
\centering
\includegraphics[width=6cm]{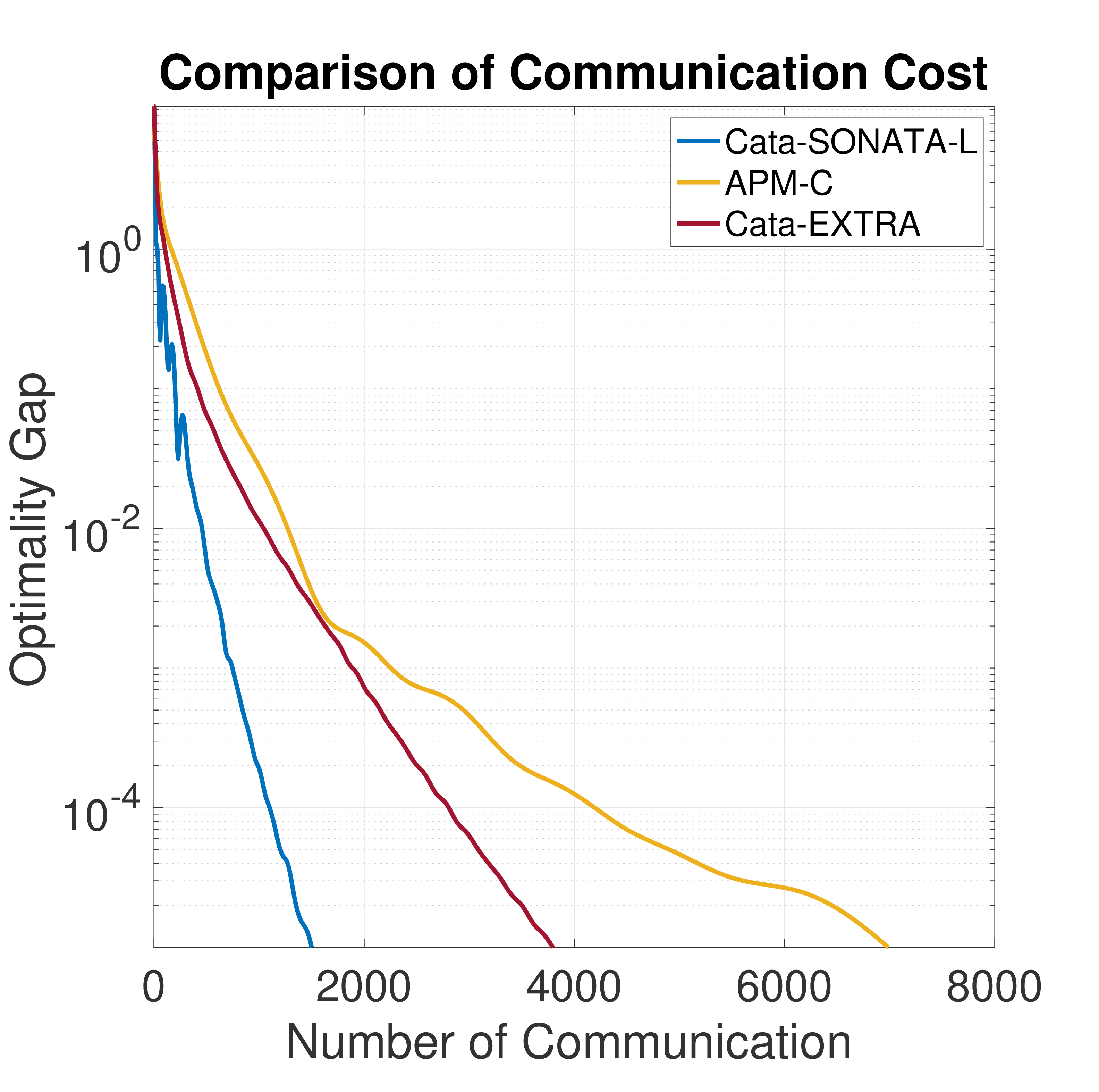}
\end{minipage}
\begin{minipage}[t]{0.48\textwidth}
\centering
\includegraphics[width=6cm]{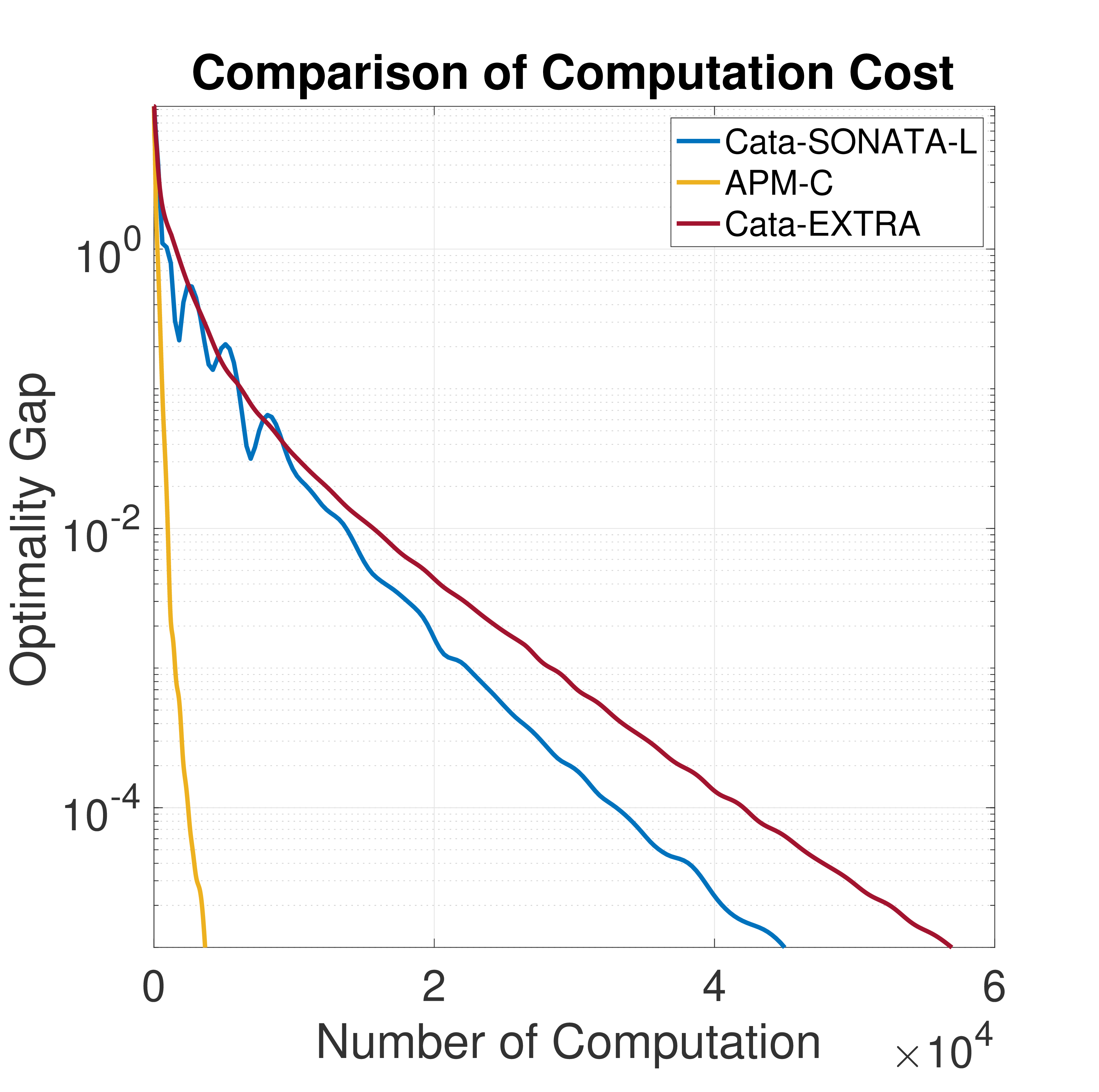}
\end{minipage}
    \caption{Comparison of distributed algorithms under convex and smooth setting in \textbf{(a):} communication cost (left) and \textbf{(b):} computation cost (right).}
    \label{'6_2_1'}
\end{figure}
As (non strongly) convex, smooth instance of Problem~\eqref{eq:problem} ($r\equiv 0$), we consider the under-determined linear regression model, with Huber loss regularization~\cite{huber2011robust}. Proposed for robust statistical procedures, the Huber norm has numerous applications in statistics and engineering-see, e.g.,~\cite{zoubir2018robust}. Specifically, we consider 
\begin{equation*}
    f_i(x)=\sum_{j=1}^n(a_{ij}^{\top}x-b_{ij})^2+\sum_{k=1}^d r_H^i(x_k),\qquad r(x)=0,
\end{equation*}
where $x_k$ denotes the $k-$th element of $x\in\mathbb{R}^d$ and the Huber norm is given by
\begin{equation*}
    r_H^i(x_k)=\left\{\begin{aligned}
    & \lambda_i(|x_k|-\frac{\lambda_i}{4\gamma_i}),& |x_k|\geq \frac{\lambda_i}{2\gamma_i},\\
    &\gamma_ix_k^2, & |x_k|<\frac{\lambda_i}{2\gamma_i}.
\end{aligned}\right.
\end{equation*}
Note that $a_{ij}\in \mathbb{R}^{d}$ and $b_{ij}\in\{-1,1\}$. Here, the data set  $\{(a_{ij},b_{ij})\}_{j=1}^n$ is assumed to be private and  owned only by  agent $i$. We use the MNIST from LIBSVM \cite{refer_libsvm} as  dataset, and   pick only $N=600$ data with feature dimension $d=784$, so that $N<d$. In the experiment, we set $\lambda_i=0.5$ for all $i$. We set $\gamma_i=0.5$ for all $i$ except $i=30$ which is set to be $50$ so that $L_{\text{max}}\approx137.72>L\approx 39.78$ and a significant difference between parameters $L_{\max}$ and $L$ is ensured.  

We compare the proposed DCatalyst-SONATA-L with the  APM-C \cite{li2020decentralized}  and Cata-EXTRA \cite{li2020revisiting}.  All the parameters are tuned according to their theoretical recommended values~\cite{li2020decentralized,li2020revisiting}. For the DCatalyst-SONATA-L algorithm,   we set $\delta=L$, the momentum parameter is calculated using~\eqref{eq:alpha-rec}, and the number of inner-loop iterations is set as $\lceil\log(k+1)\rceil$, for all $k$ [see~\eqref{inner-lp-cvx}]. 

In Figure \ref{'6_2_1'}, we plot the optimality gap achieved by the above algorithms versus the number of communications (left panel) and computations (right panel) respectively.   The figure  confirms the sublinear convergence of the algorithms. We notice that DCatalyst-SONATA-L outperforms all other algorithms in term of communication. As predicted by our theory,  DCatalyst-SONATA-L loses a bit to the APM-C in term of computations; this is  due to an extra poly-log factor and $\log\frac{1}{\epsilon}$ term appearing  in its complexity. This gap reduces when   $L_{\text{max}}$ grows with respect to   $L$.

\bibliographystyle{plain}
\bibliography{reference_}

\end{document}